\documentclass[12pt,leqno]{article}

\usepackage{amsthm}%,showidx}

\newcommand{\bm}[1]{{\mbox {\boldmath $#1$}}}

\newcommand{\bms}[1]{{\mbox {\boldmath $\scriptstyle #1$}}}

\usepackage{amsmath}
\usepackage{amscd}
\usepackage{amsmath}
\usepackage{amssymb}
\usepackage{latexsym}
\makeatletter

\@addtoreset{equation}{section}
\makeatother

\newtheorem{thm}{Theorem}[section]
\newtheorem{pr}[thm]{Proposition}
\newtheorem{df}[thm]{Definition}
\newtheorem{lm}[thm]{Lemma}
\newtheorem{cor}[thm]{Corollary}
\newtheorem{cn}[thm]{Conjecture}

\newtheorem{ex}[thm]{Example}
\newcommand{\sm}{\raisebox{2.33pt}{~\rule{6.4pt}{1.3pt}~}}

\input xy \xyoption{all} \CompileMatrices

\begin{document}

\title{The characteristic cycle
and the singular support\\
of a constructible sheaf}
\author{Takeshi Saito}

\maketitle

\begin{abstract}
We define 
the characteristic cycle of an 
\'etale sheaf
as a cycle on the cotangent
bundle of a smooth variety
in positive characteristic
using the singular support
recently defined by Beilinson.
We prove a formula \`a la Milnor
for the total dimension of the space of
vanishing cycles and an index formula
computing the Euler-Poincar\'e 
characteristic, generalizing
the Grothendieck-Ogg-Shafarevich
formula to higher dimension.

An essential ingredient of the construction
and the proof is a partial generalization
to higher dimension of the
semi-continuity of the Swan conductor
due to Deligne-Laumon.
We prove the index formula
by establishing certain functorial
properties of characteristic cycles.
\end{abstract}

As is observed by Deligne
in \cite{bp},
a strong analogy
between the wild ramification
of $\ell$-adic sheaf in
positive characteristic
and the irregular singularity
of partial differential equation
on a complex manifold
suggests to define
the characteristic cycle of an $\bar{\mathbf F}_\ell$-sheaf
as a cycle on the cotangent bundle
of a smooth variety in positive characteristic 
$p\neq \ell$.
Kashiwara and Schapira show that
the characteristic cycles
in a transcendental setting
are in fact defined directly
for constructible sheaves
without using ${\cal D}$-modules
and develop a theory in their book \cite{KSc} 
.

Recently, Beilinson \cite{Be}
defined the singular support 
as a closed conical subset of the cotangent bundle 
that controls the local acyclicity
of morphisms.
We define the characteristic cycle
as a ${\mathbf Z}$-linear 
combination of its irreducible components,
characterized by a Milnor formula 
$$-\dim{\rm tot}\ \phi_u
(j^*{\cal F},f)
=
(CC {\cal F},df)_{T^*W,u}
\leqno{\rm (\ref{eqMil})}$$
for the
total dimension of the space of
vanishing cycles.
Roughly speaking,
we realize a program sketched in
\cite{bp}.
The Milnor formula proved by
Deligne \cite{Milnor} is the case
where the sheaf is constant.
In the first version of this article,
the characteristic cycle $CC{\cal F}$
was defined as a 
${\mathbf Z}[\frac1p]$-linear 
combination.
The proof of the integrality
is due to Beilinson,
based on a suggestion by Deligne.

We also prove an index formula 
$$
\chi(X,{\cal F})
=
(CC {\cal F},
T^{*}_{X}X)_{T^{*}X}.
\leqno{\rm (\ref{eqEP})}$$
computing the Euler-Poincar\'e 
characteristic.
The Grothendieck-Ogg-Shafarevich
formula \cite{SGA5} is the case where
the variety is a curve.
Imitating the construction in
\cite[Appendix]{Gi}
in an analytic context,
we define the characteristic 
class of constructible sheaf
in Definition \ref{dfccX},
which gives an analogue of 
the MacPherson-Chern class \cite{Mac}
in positive characteristic.

To define the characteristic cycle,
it suffices to determine the coefficient
of each irreducible component
of the singular support.
We do this by imposing
the Milnor formula  (\ref{eqMil})
for morphisms defined by pencils
by choosing an embedding to
a projective space.
To prove that the coefficients are
independent of the choice
and that the characteristic cycle 
satisfies the Milnor formula
(\ref{eqMil}) in general,
we use
the continuity Proposition \ref{prMsc} of the total
dimension of the space of
vanishing cycles.
This is a partial generalization to higher
dimension of the semi-continuity
of Swan conductor \cite{DL}.

The proof of the integrality
of characteristic cycles due to
Beilinson given in Section \ref{sZ}
is done by constructing 
a function $f$ such that the section $df$
meets transversely to 
each irreducible component
of the singular support.
To complete the proof 
in an exceptional case
that occurs only in characteristic $p=2$,
we consider the pull-back
to the product with ${\mathbf A}^1$
and use the compatibility 
Proposition \ref{prsm*}
with pull-back by smooth morphisms 
of the construction
of characteristic cycle.

The index formula
Theorem \ref{thmEP}
computing the Euler-Poincar\'e 
characteristic
is deduced from the compatibility 
Theorem \ref{thmi*} 
\cite[2e Conjecture, p.~10]{bp} of
the construction of characteristic
cycles with the pull-back
by properly $C$-transversal morphisms
(Definition \ref{dfCt})
for the singular support $C$.
The compatibility 
Theorem \ref{thmi*} also implies
a description 
Theorem \ref{thmram} 
of the characteristic cycle in terms
of ramification theory \cite{AS}, \cite{nonlog}.
In the tamely ramified case,
the description has been proved by
a different method in \cite{Yang}.

We present slight modifications
of the proofs of Theorems \ref{thmi*} 
and \ref{thmEP}
due to Beilinson
in Sections \ref{sspb} and \ref{ssREP}.
The author thanks him for
allowing to include them in 
this article.
They are based on functorial properties
of characteristic cycles 
under the Radon transform.
Contrary,
the original proof which is not produced
in this article was based on ramification theory.

The only result from Section \ref{svan}
necessary for the definition
of characteristic cycle 
Theorem \ref{thmM} in Section
\ref{scc} is the continuity
Proposition \ref{prMsc}.
The contents of Section \ref{sFt}
where we establish a characterization
Proposition \ref{prfh}
of the singular support
in terms of the ${\cal F}$-transversality
introduced in Definition \ref{dfprpg}
depend only on the
first two subsections in
Section \ref{sss}
where we recall basic properties
of the singular support from \cite{Be}.

We describe briefly the content
of each section.
In Section \ref{ssvan},
we introduce and study
flat functions on a scheme
quasi-finite over a base scheme,
used to formulate the partial generalization 
of the semi-continuity
of Swan conductor to higher dimension.
After briefly recalling the generalization of
the formalism of vanishing cycles
with general base scheme
and its relation with local acyclicity,
we recall from \cite{Or}
basic properties Proposition \ref{prM}
in the case
where the locus of
non local acyclicity is quasi-finite,
in Section \ref{ssla}.
We recall and reformulate
the semi-continuity
of Swan conductor from \cite{DL}
using 
the formalism of vanishing cycles
with general base scheme
and give a partial generalization 
to higher dimension in Section \ref{sssc}.

We briefly recall definitions
and results on closed conical subsets
of the cotangent bundle and on
singular supports
from \cite{Be}
in Sections \ref{sCc} and \ref{sss}.
In Section \ref{ssram},
we give a description
Proposition \ref{prnd} of
the singular support 
in terms of ramification theory
using a characterization
Corollary \ref{corRn} of
the singular support.

We define the characteristic cycle
as characterized by the Milnor
formula (\ref{eqMil})
in Section \ref{ssCC}.
After some preliminary
on morphisms defined by pencils
and their universal family in Section \ref{sspl},
we fix some terminology
and notation to formulate
the Milnor formula
in Section \ref{ssic}.
In Proposition \ref{prfla}
in Section \ref{ssic},
we state a certain flatness criterion for
a function on isolated characteristic points
to be defined by intersections
numbers with a cycle
on the cotangent bundle.
We prove the existence of
characteristic cycle in
Theorem \ref{thmM} in Section \ref{ssCC}
by showing that the total dimension
of the space of vanishing
cycles satisfies the flatness criterion
using the continuity Proposition 
\ref{prMsc}.
We also establish some elementary properties
of characteristic cycles in Section \ref{ssCC}.
In Section \ref{sZ},
we give a proof of the integrality 
of characteristic cycle by Beilinson.

After some preliminaries
on the Chow groups of projective space
bundles in Section \ref{ssCh},
we define the characteristic class
in Section \ref{ssccc}.

We state and prove the compatibility 
Theorem \ref{thmi*}
of the construction of characteristic cycles
with properly $C$-transversal morphisms
for the singular support $C$
in Section \ref{sspb}.
In Section \ref{ssREP},
we prove the index formula
Theorem \ref{thmEP}
computing the Euler number.
It is deduced
from Theorem \ref{thmi*}
and the compatibility of 
characteristic classes
with Radon transform.
We deduce
a description of the characteristic
cycle Theorem \ref{thmram} in terms
of ramification theory
from Theorem \ref{thmi*}.

We introduce the notion 
of ${\cal F}$-transversality
in Definition \ref{dfprpg}
using a canonical morphism
(\ref{eqprpg})
and establish a characterization
Proposition \ref{prfh}
of the singular support
in terms of the ${\cal F}$-transversality
in Section \ref{sFt}.

The author thanks Pierre Deligne
for sending him an unpublished notes
\cite{bp}.
This article is the result of
an attempt to understand its contents.
The author thanks Luc Illusie
for sending him a preprint
\cite{TS}
and for the introduction to
the formalism of generalized vanishing cycles.
The author thanks Ofer Gabber for
suggesting an improvement on
the statement of Proposition \ref{prfh}.

The author thanks 
Alexander Beilinson greatly 
for his generous help on
various stages of the study
of the subjects of this article.
First for the preprint \cite{Be}
and an earlier suggestion of use of Radon transform \cite{KT}.
The proof of integrality
of characteristic cycle
which he kindly suggested the author
to include in Section \ref{sZ} is due to him.
He also kindly suggested to
include his conceptual proofs of
Theorems \ref{thmi*} and \ref{thmEP},
that allow to replace the original
more technical proofs.

The author greatly acknowledges
Ahmed Abbes
for pointing out
the similarity of the theory
of \'etale sheaves and the
content of the book
\cite{KSc},
which was the starting point
of the whole project,
and also for suggesting
a link between the vanishing topos
and the semi-continuity of Swan conductor.
The author would also like to thank
Luc Illusie and Ahmed Abbes
for inspiring discussion.
The research was partially supported
by JSPS Grants-in-Aid 
for Scientific Research
(A) 26247002.

\tableofcontents

%\newpage
\section{Vanishing topos and
the semi-continuity of the Swan conductor}
\label{svan}
\subsection{Calculus on vanishing topos}\label{ssvan}

Let $f\colon X\to S$
be a morphism of schemes.
By abuse of notation,
let $X, S$ also denote
the associated \'etale toposes.
For the definition of
the {\em vanishing topos}
$X\overset \gets\times_SS$
and the morphisms 
$$\begin{CD}
X@>{\Psi_f}>> X\overset \gets\times_SS
@>{p_2}>>S\\
@.@V{p_1}VV\\
@.X\end{CD}$$
of toposes,
we refer to \cite[1.1, 4.1, 4.3]{Il} and
\cite[1.1]{TS}.
For a geometric point $x$ of 
a scheme $X$,
we assume in this article
that the residue field
of $x$ is a separable closure
of the residue field at the
image of $x$ in $X$,
if we do not say otherwise explicitly.

For a geometric point $x$ of $X$,
the fiber 
$x\times_X(X\overset \gets\times_SS)$
of $p_1\colon
X\overset \gets\times_SS\to X$
at $x$ is the vanishing topos
$x\overset \gets\times_SS$
and is canonically
identified with the strict localization
$S_{(s)}$
at the geometric point $s=f(x)$
of $S$ defined by the image of $x$
(cf.\  \cite[(1.8.2)]{TS}).

A point on the topos $X\overset \gets\times_SS$
is defined by a triple denoted $x\gets t$
consisting of a geometric point $x$ of $X$,
a geometric point $t$ of $S$
and a specialization $s=f(x)\gets t$
namely a geometric point $S_{(s)}\gets t$
of the strict localization
lifting $S\gets t$.
The fiber 
$(X\overset \gets\times_SS)
\times_{
S\overset \gets\times_SS}
(s\gets t)$
of the canonical morphism
$X\overset \gets\times_SS\to 
S\overset \gets\times_SS$
at a point $s\gets t$
is canonically identified with
the geometric fiber $X_{s}$.
The fiber products
$X_{(x)}\times_{S_{(s)}}S_{(t)}$
and
$X_{(x)}\times_{S_{(s)}}t$
are called the {\em Milnor tube}
and the {\em Milnor fiber}
respectively.

For a commutative diagram
$$\xymatrix{
X\ar[rr]^f\ar[rd]_p&&
Y\ar[ld]^g\\
&S}$$
of morphisms of schemes,
the morphism
$\overset{\gets} g\colon X\overset \gets\times_YY\to
X\overset \gets\times_SS$
is defined by functoriality and
we have a canonical isomorphism
$\Psi_p\to \overset{\gets} g
\circ \Psi_f$.
On the fibers of a geometric point $x$
of $X$, the morphism
$\overset{\gets} g$
induces a morphism
\begin{equation}
g_{(x)}\colon Y_{(y)}
=x\overset\gets\times_YY\to S_{(s)}
=x\overset\gets\times_SS
\label{eqpx}
\end{equation}
on the strict localizations
at $y=f(x)$ and $s=p(x)$ \cite[(1.7.3)]{TS}.
In particular for $Y=X$,
we have a canonical isomorphism
$\Psi_p\to \overset{\gets} g
\circ \Psi_{\rm id}$.

Let $\Lambda$ be a finite local
ring with residue
field $\Lambda_0$ of characteristic 
$\ell$ invertible on $S$.
Let $D^+(-)$ denote
the derived category of 
complexes of $\Lambda$-modules
bounded below
and let $D^b(-)$ denote
the subcategory consisting
of complexes with bounded cohomology.
In the following,
we assume that $S$ and $X$ are quasi-compact
and quasi-separated.
We say that an object of
$D^b(X\overset \gets\times_SS)$
is constructible
if there exist finite partitions
$X=\coprod_\alpha X_\alpha$
and
$S=\coprod_\beta S_\beta$ 
by locally closed
constructible subschemes such that
the restrictions to $X_\alpha
\overset \gets\times_SS_\beta$
of cohomology sheaves are locally constant
and constructible \cite[1.3]{TS}.

Let $D^b_c(-)$ denote
the subcategory of $D^b(-)$ consisting
of constructible objects
and
let $D_{\rm ctf}(-)
\subset D^b_c(-)$ denote
its subcategory consisting
of objects of finite tor-dimension.
If $\Lambda$ is a field,
we have
$D_{\rm ctf}(-)
=D^b_c(-)$.

We canonically identify
a function on the underlying set
of a scheme $X$ with 
the function on the set
of isomorphism classes of
geometric points $x$ of $X$.
Similarly, we call a function on the set
of isomorphism classes of
points $x\gets t$ of $X\overset \gets\times_SS$
a function on $X\overset \gets\times_SS$.
We say that a function on 
$X\overset \gets\times_SS$
is a {\em constructible function}
if there exist finite partitions
$X=\coprod_\alpha X_\alpha$
and
$S=\coprod_\beta S_\beta$ as above
such that the restrictions
to $X_\alpha
\overset \gets\times_SS_\beta$
are locally constant.

For an object ${\cal K}$ of
$D_{\rm ctf}(X\overset \gets\times_SS)$,
the rank function 
$\dim {\cal K}_{x\gets t}$
is defined as a constructible function
on $X\overset \gets\times_SS$.
If $\Lambda=\Lambda_0$ is a field, we have
$\dim {\cal K}_{x\gets t}=
\sum_q(-1)^q\dim {\cal H}^q{\cal K}_{x\gets t}$.
In general, 
we have
$\dim_\Lambda {\cal K}
=
\dim_{\Lambda_0} 
{\cal K}\otimes_\Lambda^L\Lambda_0$.

\begin{df}\label{df11}
Let $Z$ be a quasi-finite scheme 
%of finite type 
over $S$ 
and let $\varphi\colon Z\to {\mathbf Q}$
be a function.
We define the {\em derivative} $\delta(\varphi)$
of $\varphi$ as a function
on $Z\overset \gets\times_SS$
by 
\begin{equation}
\delta(\varphi)(x\gets t)
=\varphi(x)
-\sum_{z\in Z_{(x)}\times
_{S_{(s)}}t}\varphi(z)
\label{eqdel}
\end{equation}
where $s=f(x)$.
If the derivative $\delta(\varphi)$
is $0$ (resp.\ $\delta(\varphi)\geqq 0$),
we say that the function $\varphi$
is {\em flat} (resp.\ {\em increasing})
over $S$.
If the morphism $f\colon Z\to S$ is finite,
we define a function $f_*\varphi$ on $S$
by 
\begin{equation}
f_*\varphi(s)=
\sum_{x\in Z_s}\varphi(x).
\label{eqf*}
\end{equation}
\end{df}

\begin{lm}\label{lmscZ}
Let $S$ be a noetherian scheme,
$Z$ be a quasi-finite scheme 
%of finite type 
over $S$ and 
$\varphi\colon Z\to {\mathbf Q}$ be a function.

{\rm 1.}
Assume that
$Z$ is \'etale over $S$.
Then $\varphi$ is locally constant
(resp.\ upper semi-continuous)
if and only if
it is flat over $S$
(resp.\ constructible
and increasing over $S$).

{\rm 2.} 
The function $\varphi$ is
constructible if and only if its
{\em derivative} $\delta(\varphi)
\colon Z\overset \gets\times_SS
\to {\mathbf Q}$
defined in {\rm (\ref{eqdel})}
is constructible.
Consequently, if $\varphi$ is flat over $S$,
then $\varphi$ is constructible.

{\rm 3.} 
Assume that
$\varphi$ is {\em flat}
over $S$.
Then, $\varphi=0$
if and only if
$\varphi(x)=0$ for
the generic point $x$ of
every irreducible component of $Z$.

{\rm 4.}
Assume that 
the morphism $f\colon Z\to S$ is finite.
If $\varphi$ is constructible,
the function $f_*\varphi$
is also constructible.
Assume that $\varphi$ is constructible
and is increasing over $S$.
Then the function $f_*\varphi$ on $S$ is 
upper semi-continuous.
Further, the function $\varphi$ is flat over $S$ if 
and only if
$f_*\varphi$ is locally constant.
\end{lm}

\proof{
1. 
Since the question is \'etale local on $Z$,
we may assume that $Z\to S$ 
is an isomorphism.
Then the assertion is clear.

2.
Assume $\delta(\varphi)$ is constructible.
By noetherian induction,
it suffices to show the following:
For every geometric point $t$ of $S$
and the closure $T\subset S$ of its image,
there exists a dense open subset
$V\subset T$ such that
$\varphi$ is locally constant on
$Z\times_SV$.
Replacing $S$ by $T$,
it suffices to consider the case where 
$t$ dominates
the generic point of an irreducible scheme
$S$.
For a geometric point $x$ of $Z$
above $t$,
we have $\delta(\varphi)(x\gets t)=0$.
By further replacing $X$ by
an \'etale neighborhood of $x$,
it suffices to consider the case
where $Z$ is \'etale over $S$
and $\delta(\varphi)=0$.
Then, by 1,
$\varphi$ is locally constant and hence
constructible.

Assume $\varphi$ is constructible.
For a closed subset $T$ of $S$,
the subtopos
$(Z\sm Z\times_ST)\overset\gets\times_ST$
is empty.
Hence, 
by noetherian induction,
it suffices to show the following:
For every geometric point $t$ of $S$
and the closure $T\subset S$ of its image,
there exists a dense open subset
$V\subset T$ such that
$\delta(\varphi)$ is locally constant on
$(Z\times_SV)\overset\gets\times_SV$.
Similarly as above,
it suffices to consider the case where 
$Z$ is \'etale over $S$
and $\varphi$ is locally constant.
Then, by 1.,
$\delta(\varphi)=0$ and is constructible.

3.
By (\ref{eqdel}),
a function flat over $S$ is uniquely determined
by the values at the generic points
of irreducible components.

4.
Assume that $f\colon Z\to S$ is finite.
If $\varphi$ is constructible,
there exists a dense open subscheme
$U$ of $S$ such that
$\varphi$ is locally constant
on $Z\times_SU$.
Hence 
$f_*\varphi$ is constructible by
noetherian induction.
For a specialization $s\gets t$,
we have
$\sum_{x\in Z_s}\delta(\varphi)(x\gets t)
=f_*\varphi(s)-
f_*\varphi(t)
=\delta(f_*\varphi)(s\gets t)$.
Hence we may assume
$Z=S$ and then the assertion is clear.
\qed}
\medskip

We give an example
of flat function.
Let $S$ be a noetherian
scheme, $X$ be a
scheme of finite type over $S$
and $Z\subset X$ be a closed
subscheme quasi-finite over $S$.
Let $A$ be a complex of 
${\cal O}_X$-modules
such that the cohomology sheaves
${\cal H}^q(A)$ are coherent
${\cal O}_X$-modules supported on $Z$
for all $q$ and
that $A$ is of finite tor-dimension 
as a complex of ${\cal O}_S$-modules.
For a geometric point $z$ of $Z$ 
and its image $s$ in $S$,
let ${\cal O}_{S,s}$ denote the {\em strict}
localization
and $k(s)$ the separably closed
residue field of ${\cal O}_{S,s}$.
Then, the $k(s)$-vector spaces
$Tor^{{\cal O}_{S,s}}_q
(A_z,k(s))$ are of finite dimension
and are $0$ except for finitely many $q$.
We define a function
$\varphi_A\colon Z\to {\mathbf Z}$
by
\begin{equation}
\varphi_A(z)
=
\sum_q
(-1)^q
\dim_{k(s)}
Tor^{{\cal O}_{S,s}}_q
(A_z,k(s)).
\label{eqA}
\end{equation}

\begin{lm}\label{lmA}
Let noetherian schemes
$Z\subset X\to S$ 
and a complex $A$ be as above.

{\rm 1.}
The function 
$\varphi_A\colon Z\to {\mathbf Z}$
defined by {\rm (\ref{eqA})} is
constructible and flat over $S$.

{\rm 2.}
Suppose that $S$ and $Z$ are integral
and that the image of the generic point
$\xi$ of $Z$ is the generic point
$\eta$ of $S$.
If $A={\cal O}_Z$,
the value of the function $\varphi_A$
at a geometric point of $Z$ above $\xi$ is
the inseparable degree $[k(\xi):k(\eta)]
_{\rm insep}$.
\end{lm}

The condition that
$A$ is of finite tor-dimension 
as a complex of ${\cal O}_S$-modules
is satisfied if $S$ is regular.

\proof{1.
Since the assertion is
\'etale local on $Z$,
we may assume that $Z$ is
finite over $S$, 
that $X$ and $S$ are affine and that
$z$ is the unique point in the geometric fiber $Z
\times_S{\rm Spec}\ k(s)$.
Then, the complex $Rf_*A$
is a perfect complex of
${\cal O}_S$-modules
and $\varphi_A(z)$ equal the rank of
$Rf_*A$.
Hence, the assertion follows.

2.
We may assume that
$S={\rm Spec}\ k(\eta)$
and $Z={\rm Spec}\ k(\xi)$
and the assertion follows.
\qed}
\medskip

We generalize the definition
of derivative
to functions on vanishing topos.

\begin{df}\label{df14}
Let 
\begin{equation}
\xymatrix{
Z\ar[rr]^f\ar[rd]_p&&
Y\ar[ld]^g\\
&S}
\label{eqfgp}
\end{equation} 
be a commutative diagram of
morphisms of schemes
such that $Z$ is quasi-finite over $S$.
Let $\psi\colon 
Z\overset\gets\times_YY\to {\mathbf Q}$
be a function such that
$\psi(x\gets w)=0$
unless $w$ is not supported
on the image of $f_{(x)}
\colon Z_{(x)}\to Y_{(y)}$
where $y=f(x)$.
We define the derivative
$\delta(\psi)$ as 
a function on 
$Z\overset\gets\times_SS\to {\mathbf Z}$
by
\begin{equation}
\delta(\psi)(x\gets t)
=\psi(x\gets y)
-\sum_{w\in Y_{(y)}\times_{S_{(s)}}t}
\psi(x\gets w)
\label{eqdelY}
\end{equation}
where $y=f(x)$ and $s=p(x)$.
We say that $\psi$ is {\em flat} over $S$
if $\delta(\psi)=0$.
\end{df}

The sum in the right hand side
of (\ref{eqdelY})
is a finite sum by the assumption
that $Z$ is quasi-finite over $S$
and the assumption on the support of $\psi$.
If $Z=Y$, we recover the
definition (\ref{eqdel})
by applying (\ref{eqdelY})
to the pull-back $p_2^*\varphi\colon
Z\overset\gets\times_ZZ\to {\mathbf Z}$
by $p_2\colon
Z\overset\gets\times_ZZ\to Z$.

The following elementary Lemma will be used
in the proof of a generalization
of the continuity of the Swan conductor.

\begin{lm}\label{lmdelp}
Let the assumption on the diagram
{\rm (\ref{eqfgp})} be as in
Definition {\rm \ref{df14}} and let
$\varphi$ be a function on $Z$.
We define a function
$\psi$ on $Z\overset\gets\times_YY$
by
\begin{equation}
\psi(x\gets w)=
\sum_{z\in Z_{(x)}\times_{Y_{(y)}}w}
\varphi(z)
\label{eqpsi}
\end{equation}
where $y=f(x)$.
Then the derivative $\delta(\varphi)$
on $Z\overset\gets\times_SS$
defined by {\rm (\ref{eqdel})}
equals 
$\delta(\psi)$
defined by {\rm (\ref{eqdelY})}.
\end{lm}

\proof{
It follows from
$\psi(x\gets y)=\varphi(x)$
and 
$Z_{(x)}\times_{S_{(s)}}t
=
\coprod_{w\in Y_{(y)}\times_{S_{(s)}}t}
(Z_{(x)}\times_{Y_{(y)}}w)$.
Note that except for finitely many
geometric points $w$ of
$Y_{(y)}\times_{S_{(s)}}t$
those supported on the image of
$Z_{(x)}$, the fiber
$Z_{(x)}\times_{Y_{(y)}}w$
is empty.
\qed}

\subsection{Nearby cycles and the local acyclicity}\label{ssla}

For a morphism $f\colon X\to S$,
the morphism $\Psi_f\colon X\to 
X\overset \gets\times_SS$
defines the nearby cycles functor
$R\Psi_f:D^+(X)\to D^+(X\overset \gets\times_SS)$.
The canonical morphism
$p_1^*\to R\Psi_f$
of functors is defined by adjunction
and by the isomorphism ${\rm id}\to p_1\circ \Psi_f.$
The cone of the morphism
$p_1^*\to R\Psi_f$
defines the vanishing cycles functor
$R\Phi_f:D^+(X)\to D^+(X\overset \gets\times_SS)$.
If $S$ is the spectrum of 
a henselian discrete valuation ring
and if $s,\eta$ denote
its closed and generic points,
we recover the classical construction
of complexes $\psi, \phi$ of nearby cycles
and vanishing cycles as 
the restrictions to
$X_s\overset \gets\times_S\eta$
of $R\Psi_f$ and $R\Phi_f$
respectively.

Recall that $f\colon X\to S$
is said to be locally acyclic relatively
to a complex ${\cal F}$ 
of $\Lambda$-modules
on $X$
{\rm \cite[Definition 2.12]{TF}}
if the canonical morphism
\begin{equation}
\begin{CD}
{\cal F}_x
@>>>
R\Gamma(X_{(x)}\times_{S_{(s)}}t,
{\cal F}|_{X_{(x)}\times_{S_{(s)}}t})
\end{CD}
\label{eqla}
\end{equation}
is an isomorphism
for every $x\gets t$.
Recall that $f\colon X\to S$
is said to be universally locally acyclic relatively
to ${\cal F}$,
if for every morphism $S'\to S$,
its base change is locally acyclic relatively
to the pull-back of ${\cal F}$.

\begin{lm}[{\rm cf.\ \cite[Proposition 7.6.2]{Fu}}]\label{lmctf}
Let $f\colon X\to S$
be a morphism of finite type and
let ${\cal F}\in D_{\rm ctf}(X)$
be a complex of finite tor-dimension.

{\rm 1.}
Suppose that ${\cal F}$ is of tor-amplitude $[a,b]$ and that $f\colon X\to S$ is of relative dimension $d$. Then, 
for points $x\gets t$
of $X\overset{\gets}\times_SS$,
the complex
$R\Gamma(X_{(x)}\times_{S_{(s)}}t,
{\cal F}|_{X_{(x)}\times_{S_{(s)}}t})$
of $\Lambda$-modules
is of %finite tor-dimension
tor-amplitude $[a,b+d]$ 
and,
for a $\Lambda$-module $M$,
the canonical morphism
\begin{equation}
R\Gamma(X_{(x)}\times_{S_{(s)}}t,
{\cal F}|_{X_{(x)}\times_{S_{(s)}}t})
\otimes_\Lambda^LM
\to  
R\Gamma(X_{(x)}\times_{S_{(s)}}t,
{\cal F}|_{X_{(x)}\times_{S_{(s)}}t}
\otimes_\Lambda^LM)
\label{eqFu}
\end{equation}
is an isomorphism.

{\rm 2.}
Let $\Lambda_0$
be the residue field of $\Lambda$.
Then,
$f\colon X\to S$ is locally acyclic
(resp.\ universally locally acyclic)
relatively to 
${\cal F}$ if and only if
it is so relatively to
${\cal F}_0=
{\cal F}\otimes_\Lambda^L\Lambda_0$.
\end{lm}

\proof{
1.
By the assumption that
$f\colon X\to S$
is of finite type and of
relative dimension $d$,
the functor 
$R\Gamma(X_{(x)}\times_{S_{(s)}}t,-)$
is of cohomological dimension $\leqq d$
by \cite[Corollaire 3.2]{XIV}.
Hence, similarly as \cite[(4.9.1)]{Rapport},
the canonical morphism (\ref{eqFu})
is an isomorphism.
Since the complex
${\cal F}|_{X_{(x)}\times_{S_{(s)}}t}
\otimes_\Lambda^LM$
is acyclic outside $[a,b]$,
the complex
$R\Gamma(X_{(x)}\times_{S_{(s)}}t,$
${\cal F}|_{X_{(x)}\times_{S_{(s)}}t}
\otimes_\Lambda^LM)$
is acyclic outside $[a,b+d]$.
Thus, the complex
$R\Gamma(X_{(x)}\times_{S_{(s)}}t,
{\cal F}|_{X_{(x)}\times_{S_{(s)}}t})$ is 
of tor-amplitude $[a,b+d]$.

2.
It suffices to show the 
assertion for local acyclicity.
If the canonical morphism 
(\ref{eqla})
is an isomorphism for $x\gets t$,
then (\ref{eqla})
for ${\cal F}_0$ is an isomorphism
by the isomorphism (\ref{eqFu}) 
for $M=\Lambda_0$.

To show the converse,
let $I^\bullet$ be a filtration
by ideals of $\Lambda$
such that ${\rm Gr}\Lambda$
is a $\Lambda_0$-vector space.
Then, $I^\bullet$ defines a filtration
on ${\cal F}={\cal F}\otimes^L\Lambda$ and
a canonical isomorphism
${\rm Gr}{\cal F}\to
{\cal F}_0\otimes_{\Lambda_0}
{\rm Gr}\Lambda$.
Hence if (\ref{eqla})
is an isomorphism for ${\cal F}_0$ then
(\ref{eqla})
for ${\cal F}$ is an isomorphism.
\qed}
\medskip

We consider a commutative diagram
\begin{equation}
\xymatrix{
X\ar[rr]^f\ar[rd]_p&&
Y\ar[ld]^g\\
&S}
\label{eqXYS}
\end{equation} 
of schemes.
The canonical isomorphism
$\overset \gets g\circ \Psi_f
\to \Psi_p$
induces an isomorphism of functors
\begin{equation}
R\overset \gets g_{*}\circ R\Psi_f
\to R\Psi_p
\label{eqPsPs}
\end{equation}

For an object ${\cal K}$
of $D^+(X\overset \gets\times_YY)$
and a geometric point $x$ of $X$,
the restriction of $R\overset{\gets} g_*{\cal K}$
on $x\overset \gets\times_SS=S_{(s)}$
for $s=f(x)$
is canonically identified
with $Rg_{(x)*}({\cal K}|_{Y_{(y)}})$
for $y=f(x)$
in the notation of (\ref{eqpx})
by \cite[(1.9.2)]{TS}.
For the stalk at a point $x\gets t$ of
$X\overset \gets\times_SS$,
this identification gives a canonical isomorphism
\begin{equation}
R\overset{\gets} g_*{\cal K}_{x\gets t}
\to
R\Gamma(Y_{(y)}\times_{S_{(s)}}
S_{(t)},{\cal K}|_{Y_{(y)}\times_{S_{(s)}}
S_{(t)}}).
\label{eqMt}
\end{equation}
For an object ${\cal F}$ of $D^+(X)$,
(\ref{eqMt}) applied to $Y=X$ 
gives a canonical identification
\begin{equation}
R\Psi_p{\cal F}_{x\gets t}
\to
R\Gamma(X_{(x)}\times_{S_{(s)}}
S_{(t)},{\cal F}|_{X_{(x)}\times_{S_{(s)}}
S_{(t)}})\label{eqMt2}
\end{equation}
with the cohomology
of the Milnor tube \cite[(1.1.15)]{TS}.

A cartesian diagram
$$\begin{CD}
X@<<< X_T\\
@VfVV @VV{f_T}V\\
S@<i<< T
\end{CD}$$
of schemes defines
a 2-commutative diagram
$$
\begin{CD}
X_T@<{p_1}<<
X_T\overset \gets\times_TT
@<{\Psi_{f_T}}<< X_T\\
@ViVV@V{\overset \gets i}VV@VViV\\
X@<{p_1}<<
X\overset \gets\times_SS
@<{\Psi_f}<< X
\end{CD}$$
and
the base change morphisms define
a morphism
\begin{equation}
\begin{CD}
@>>> \overset{\gets*} ip_1^*
@>>> \overset{\gets*} iR\Psi_f
@>>> \overset{\gets*} iR\Phi_f@>>>\\
@.@V{\simeq}VV @VVV @VVV\\
@>>> p_1^*i^*
@>>> R\Psi_{f_T}i^*
@>>> R\Phi_{f_T}i^*@>>>\\
\end{CD}
\label{eqbc}
\end{equation}
of distinguished triangles of functors.
For an object ${\cal F}$ of $D^+(X)$,
we say that the formation of
$R\Psi_f{\cal F}$ commutes with
the base change $T\to S$ if
the middle vertical arrow defines
an isomorphism
$\overset{\gets*} iR\Psi_f{\cal F}
\to R\Psi_{f_T}i^*{\cal F}$.

For a point $x\gets t$
of $X\overset\gets\times_SS$,
if $T\subset S$ denotes
the closure of the image of $t$ in $S$,
the left square of (\ref{eqbc})
induces a commutative diagram
\begin{equation}
\xymatrix{
(p_1^*{\cal F})_{x\gets t}
={\cal F}_x\ar[rr]\ar[rrd]
&&
R\Psi_f{\cal F}_{x\gets t}
=
R\Gamma(X_{(x)}\times_{S_{(s)}}
S_{(t)},{\cal F}|_{X_{(x)}\times_{S_{(s)}}
S_{(t)}})\ar[d]\\
&&
R\Psi_{f_T}({\cal F}|_{X_T})_{x\gets t}
=
R\Gamma(X_{(x)}\times_{S_{(s)}}
t,{\cal F}|_{X_{(x)}\times_{S_{(s)}}
t}).}
\label{eqMf}
\end{equation}
The vertical arrow
is the canonical morphism from the cohomology 
of the Milnor tube
to that of the Milnor fiber
and the slant arrow
is the canonical morphism (\ref{eqla}).
Recall that we assume
that the residue field of $t$
is a separable closure of
the residue field at
the image in $S_{(s)}$.

We interpret the local acyclicity
in terms of vanishing topos.

\begin{pr}\label{lmapp}
Let $f\colon X\to S$
be a morphism of schemes.
Then, for an object ${\cal F}$ of $D^+(X)$, 
the conditions %{\rm (1)} and {\rm (2)}
in {\rm 1.{}\!} and {\rm 2.{}\!} 
below are equivalent to each other respectively.

{\rm 1.}
{\rm (1)}
For every finite morphism $g\colon T\to S$,
for every geometric point $x$ of $X$
and 
for every specialization
$t\gets u$ of
geometric points of $T$
such that 
$s=f(x)=g(t)$
as geometric points of $S$,
the canonical morphism
\begin{equation}
R\Gamma(X_{(x)}\times_{S_{(s)}}
T_{(u)},{\cal F}|_{X_{(x)}\times_{S_{(s)}}
T_{(u)}})
\to
R\Gamma(X_{(x)}\times_{S_{(s)}}
u,{\cal F}|_{X_{(x)}\times_{S_{(s)}}
u})
\label{eqSTt}
\end{equation}
that is a vertical arrow in {\rm (\ref{eqMf})}
for the base change $f_T\colon X_T
=X\times_ST\to T$
is an isomorphism.

{\rm (2)}
The formation of $R\Psi_f{\cal F}$
commutes with finite base change $T\to S$.

{\rm 2. ({\cite[Corollaire 2.6]{app}})}
{\rm (1)}
The morphism $f\colon X\to S$ is 
(resp.\ universally) locally 
acyclic relatively to ${\cal F}$.

{\rm (2)}
The canonical morphism
$p_1^*{\cal F}\to
R\Psi_f{\cal F}$
is an isomorphism
and the formation of $R\Psi_f{\cal F}$
commutes with finite (resp.\ arbitrary) base change $T\to S$.

{\rm (3)}
The canonical morphism
$p_1^*{\cal F}_T\to
R\Psi_{f_T}{\cal F}_T$
is an isomorphism
for every finite 
(resp.\ every) morphism $T\to S$
and the pull-back ${\cal F}_T$
of ${\cal F}$ on $X_T=X\times_ST$.
\end{pr}

\proof{
1.
Let $T\to S$ be a finite morphism
and $x\mapsto s$ and $t\gets u$ be as in
the condition (1).
Then, 
for the geometric point $x'$ of
$X_T$ defined by a unique point of
$x\times_st$,
the Milnor tube
$X_{T (x')}\times_{T_{(t)}}
T_{(u)}$
is canonically isomorphic to
$X_{(x)}\times_{S_{(s)}}
T_{(u)}$ since $T\to S$ is finite.
Thus, the morphism
(\ref{eqSTt}) is identified with
the morphism
$R\Gamma(X_{(x)}\times_{S_{(s)}}
T_{(u)},{\cal F}|_{X_{(x)}\times_{S_{(s)}}
T_{(u)}})
\to
R\Gamma(X_{(x)}\times_{S_{(s)}}
u,{\cal F}|_{X_{(x)}\times_{S_{(s)}}
u})$
that is the vertical arrow in {\rm (\ref{eqMf})}
for $f_T\colon X_T\to T$
at $x'\gets u$.

Let $T'\subset T$
be the closure of the image of $u$.
We consider the
base change morphisms
\begin{equation}
\begin{CD}
R\Psi_f{\cal F}_{x\gets u}
=&
R\Gamma(X_{(x)}\times_{S_{(s)}}
S_{(u)},{\cal F}|_{X_{(x)}\times_{S_{(s)}}
S_{(u)}})\\
&@VVV\\
R\Psi_{f_T}({\cal F}|_{X_T})_{x\gets u}
=&
R\Gamma(X_{(x)}\times_{S_{(s)}}
T_{(u)},{\cal F}|_{X_{(x)}\times_{S_{(s)}}
T_{(u)}})\\
&@VVV\\
R\Psi_{f_{T'}}({\cal F}|_{X_{T'}})_{x\gets u}
=&
R\Gamma(X_{(x)}\times_{S_{(s)}}
u,{\cal F}|_{X_{(x)}\times_{S_{(s)}}
u}).
\end{CD}
\label{eqSTft}
\end{equation}
The condition (1) implies
that the lower arrow and
the composition are isomorphisms.
Hence, the upper arrow
is an isomorphism
and we have
(1)$\Rightarrow$(2).

Conversely, 
the condition (2) implies
that the upper arrow and
the composition are isomorphisms.
Hence, the lower arrow
that is the same as
(\ref{eqSTt}) is an isomorphism.
Thus, we have
(2)$\Rightarrow$(1).

2. 
First, we show the
equivalence of (1) and (2)
in the cases without resp.
The condition (1) is equivalent to
that the slant arrow in (\ref{eqMf})
is an isomorphism
for every point $x\gets t$
of $X\overset \gets\times_SS$.
Hence the condition (2) implies the condition (1)
by (2)$\Rightarrow$(1) in 1.\
and the commutativity of the diagram (\ref{eqMf}).

Conversely, by \cite[Corollaire 2.6]{app},
if the condition (1) is satisfied,
the formation of
$Rf_{(x)*}({\cal F}|_{X_{(x)}})$
commutes with finite base change
for every geometric point $x$ of $X$
where $f_{(x)}\colon X_{(x)}\to S_{(s)}$
is the morphism on the strict localizations
induced by $f$.
Thus, for every finite morphism
$T\to S$ and every point $x\gets u$
of $X_T\overset\gets\times_TT$,
the upper arrow in (\ref{eqSTft})
is an isomorphism.
Hence %by (1)$\Rightarrow$(2) in 1.\
the formation of $R\Psi_f{\cal F}$
commutes with finite base change $T\to S$.
Further the vertical arrow in (\ref{eqMf})
is an isomorphism. Thus
the canonical morphism
$p_1^*{\cal F}\to
R\Psi_f{\cal F}$
is an isomorphism
further by 
the commutativity of the diagram (\ref{eqMf}).

If $p_1^*{\cal F}\to
R\Psi_f{\cal F}$
is an isomorphism,
the formation of
$R\Psi_f{\cal F}$ commutes
with base change $T\to S$
if and only if
$p_1^*{\cal F}_{X_T}\to
R\Psi_{f_T}{\cal F}_{X_T}$
is an isomorphism
by the left square of (\ref{eqbc}).
Hence we have
an equivalence (2)$\Leftrightarrow$(3).
The equivalence (1)$\Leftrightarrow$(3)
in the cases with resp.\ follows
immediately from that without resp.
\qed}

\begin{pr}\label{prM}
Let $f\colon X\to S$
be a morphism of finite type of noetherian schemes
and let $Z\subset X$ be a closed
subscheme quasi-finite over $S$.
Let ${\cal F}$
be an object of
$D^b_c(X)$
such that the restriction of $f\colon X\to S$
to the complement $X\sm Z$ is 
(resp.\ universally) locally acyclic
relatively to the restriction of ${\cal F}$.

{\rm 1. (cf.\ {\cite[Proposition 6.1]{Or}})}
$R\Psi_f{\cal F}$ and $R\Phi_f{\cal F}$
are constructible
and their formations commute
with finite (resp.\ arbitrary)
base change.
The constructible object 
$R\Phi_f{\cal F}$ is supported
on $Z\overset \gets\times_SS$.

{\rm 2. }
Let $x$ be a geometric point of $X$
and $s=f(x)$ be the geometric point
of $S$ defined by the image of $x$ by $f$.
Let $t$ and $u$ be geometric points
of $S_{(s)}$
and  $t\gets u$ be a specialization. 
Then, there exists a distinguished triangle
\begin{equation}
\begin{CD}
@>>>
R\Psi_f{\cal F}_{x\gets t}
@>>>
R\Psi_f{\cal F}_{x\gets u}
@>>>
{\displaystyle
\bigoplus_{z\in (Z\times_XX_{(x)})
\times_{S_{(s)}}t}}
R\Phi_f{\cal F}_{z\gets u}
@>>>
\end{CD}
\label{eqDel}
\end{equation}
where 
$R\Psi_f{\cal F}_{x\gets t}
\to
R\Psi_f{\cal F}_{x\gets u}$
is the cospecialization.
\end{pr}

The commutativity of the formation of $R\Psi_f{\cal F}$ with any base change implies
its constructibility by \cite[8.1, 10.5]{Or}
as noted after \cite[Theorem 1.3.1]{TS}.

\proof{
1.
The constructibility 
is proved by
taking a compactification
in \cite[Proposition 6.1]{Or}.
The commutativity with base change
is proved similarly by
taking a compactification
and applying the proper base change
theorem.

The assertion on the support of $R\Phi_f{\cal F}$
follows from Proposition \ref{lmapp}.2 (1)$\Rightarrow$(2).

2.
Let $t$ and $u$
be geometric points of
$S_{(s)}$ and $t\gets u$ be a specialization.
By replacing $S$ by the strict localization
$S_{(s)}$
and shrinking $X$, 
we may assume that $S=S_{(s)}$, that $X$ is affine
and that $Z=Z\times_XX_{(x)}$ is finite over $S$.

We consider the diagram
$$\begin{CD}
s@. t\\
@V{i_{s}}VV @VV{i_t}V\\
S@<j<< S_{(t)}@<k<< u
\end{CD}$$
and let the morphisms obtained by
the base change $X\to S$
denoted by the same letters,
by abuse of notation.
Similarly as the sliced vanishing cycles
in the proof of \cite[Proposition 6.1]{Or},
we consider an object $\Phi_{t\gets u}{\cal F}$ on 
$X\times_SS_{(t)}$ fitting in the distinguished
triangle
$\to j^*{\cal F}\to Rk_*(j\circ k)^*{\cal F}\to 
\Phi_{t\gets u}{\cal F}\to $.
Since the formation of $R\Psi_f{\cal F}$ commutes
with finite base change by 1.,
we have a distinguished
triangle (\ref{eqDel})
with the third term replaced by
$\Delta_x=(Rj_*\Phi_{t\gets u}{\cal F})_x$.
Further,
the third term itself is canonically isomorphic
to the direct sum of
$\Delta_z=(\Phi_{t\gets u}{\cal F})_z$
for $z\in Z_t$.

Since  $R\Phi_f{\cal F}$ is acyclic
outside $Z\overset\gets\times_SS$ by 1.,
the canonical morphisms
$i_s^*{\cal F}\to 
i_s^*Rj_*j^*{\cal F}$
and
$i_s^*{\cal F}\to 
i_s^*R(j\circ k)_*(j\circ k)^*{\cal F}$
are isomorphisms on $X_s\sm Z_s$.
Hence, 
the restriction $i_s^*Rj_*\Phi_{t\gets u}{\cal F}$
is acyclic on  $X_s\sm Z_s$.
Similalry, 
the restriction $i_t^*\Phi_{t\gets u}{\cal F}$
is acyclic on  $X_t\sm Z_t$.

We take a compactification $\bar X$ of $X$
and an extension $\bar {\cal F}$ of ${\cal F}$ to $\bar X$.
Define $\Phi_{t\gets u}\bar {\cal F}$ on 
$\bar X\times_SS_{(t)}$ similarly as
$\Phi_{t\gets u}{\cal F}$
and set $Y=\bar X\sm X$.
By the proper base change theorem,
the canonical morphisms
$R\Gamma(\bar X_s,
i_s^*Rj_*\Phi_{t\gets u}\bar {\cal F})
\gets
R\Gamma(\bar X,
Rj_*\Phi_{t\gets u}\bar {\cal F})
\to
R\Gamma(\bar X\times_SS_{(t)},
\Phi_{t\gets u}\bar {\cal F})
\to
R\Gamma(\bar X_t,
i_t^*\Phi_{t\gets u}\bar {\cal F})$
are isomorphisms
and similarly for the restrictions to $Y$.
Hence, we obtain a commutative diagram
$$\begin{CD}
\Delta_{x}=&R\Gamma(Z_s,
i_s^*Rj_*\Phi_{t\gets u}{\cal F})
@>>>
R\Gamma_{c}(X_s,
i_s^*Rj_*\Phi_{t\gets u}{\cal F})
\\
&@VVV@VVV\\
\bigoplus_{z\in Z_t}
\Delta_z=
&R\Gamma(Z_t,
i_t^*\Phi_{t\gets u} {\cal F})
@>>>
R\Gamma_c(X_t,
i_t^*\Phi_{t\gets u} {\cal F})
\end{CD}$$
of isomorphisms
and the assertion follows.
\qed}

\begin{cor}\label{corsl}
We keep the assumptions in Proposition {\rm \ref{prM}}
and 
let $x$ be a geometric point of $X$
and $s=f(x)$ be the geometric point
of $S$ defined by the image of $x$ by $f$
as in Proposition {\rm \ref{prM}.2}.
Then, the restriction of
the constructible sheaf
$R^q\Psi_f{\cal F}$
on $x\overset\gets\times_SS
=S_{(s)}$
is locally constant 
outside the image of
the finite scheme $Z\times_XX_{(x)}$
for every $q$.
\end{cor}

\proof{
Let $t$ and $u$
be geometric points of $S_{(s)}$
not in the image of $Z\times_XX_{(x)}$
and $t\gets u$ be a specialization.
Since $R\Psi_f{\cal F}$
is constructible,
it suffices to show that
the cospecialization morphism
$R\Psi_f{\cal F}_{x\gets t}
\to
R\Psi_f{\cal F}_{x\gets u}$ is an isomorphism.
Then by the assumption on
the local acyclicity,
the complex $\Phi_{t\gets u}{\cal F}$ 
in the proof of Proposition \ref{prM}.2 is acyclic.
Hence the assertion follows from (\ref{eqDel})
and the isomorphism
$R\Phi_f{\cal F}_{z\gets u}\to 
(\Phi_{t\gets u}{\cal F})_z$.
\qed}

\begin{cor}\label{corsc}
We keep the assumptions in Proposition {\rm \ref{prM}}.
We further assume that 
${\cal F}$ is of finite tor-dimension.

{\rm 1.}
The complexes
$R\Psi_f{\cal F}$
and 
$R\Phi_f{\cal F}$
are of finite tor-dimension.
Consequently,
the functions
$\dim R\Psi_f{\cal F}$
and 
$\dim R\Phi_f{\cal F}$
are defined and constructible.

{\rm 2.}
Define a constructible
function $\delta_{\cal F}$ 
on $X\overset\gets\times_SS$
supported on $Z\overset \gets\times_SS$ 
by $\delta_{\cal F}(x\gets t)=
\dim R\Phi_f{\cal F}_{x\gets t}$.
Let $\Lambda_0$
denote the residue field
of $\Lambda$ and assume that
$R\Phi_f{\cal F}
\otimes^L_{\Lambda_0}\Lambda$ 
is acyclic except at degree $0$.
Then, we have $\delta_{\cal F}\geqq 0$
and the equality $\delta_{\cal F}=0$
is equivalent to the condition that
the morphism $f$ is (resp.\ universally) locally acyclic
relatively to ${\cal F}$.
\end{cor}

\proof{
1. By Proposition \ref{prM}.1,
Proposition \ref{lmapp}.1
and Lemma \ref{lmctf}.1,
the complex
$R\Psi_f{\cal F}$
is of finite tor-dimension
and hence
$R\Phi_f{\cal F}$
is also of finite tor-dimension.
Since they are constructible 
by Proposition \ref{prM}.1,
the functions
$\dim R\Psi_f{\cal F}$
and 
$\dim R\Phi_f{\cal F}$
are defined and constructible.

2.
The positivity $\delta_{\cal F}\geqq 0$
follows from 
the assumption that
$R\Phi_f{\cal F}$ is acyclic except at degree $0$.
Further the equality $\delta_{\cal F}=0$
is equivalent to 
$R\Phi_f{\cal F}=0$.
Since the formation of 
$R\Psi_f{\cal F}$
commutes with finite (resp.\ arbitrary) 
base change by Proposition \ref{prM}.1,
it is further
equivalent to the condition that
the morphism $f$ is (resp.\ universally) locally acyclic
relatively to ${\cal F}$
by Proposition \ref{lmapp}.2.
\qed}

\begin{lm}\label{lmP}
Assume that
$\Lambda=\Lambda_0$
is a field.
Then, the assumption 
that $R\Phi_f{\cal F}$ 
is acyclic except at degree $0$
in Corollary {\rm \ref{corsc}.2}
is satisfied if the following conditions
are satisfied:
The scheme $S$ is noetherian,
the restriction of $f\colon X\to S$
to $X\sm Z$ is 
universally locally acyclic
relatively to the restriction of ${\cal F}$
and the following condition {\rm (P)}
is satisfied.

{\rm (P)}
For every morphism
$T\to S$ from the spectrum $T$ of
a discrete valuation ring,
the pull-back of ${\cal F}[1]$
to $X_T$ is perverse.
\end{lm}

\proof{
Let $x\gets t$ be a point of
$X\overset \gets\times_SS$ 
and let $T\to S$ be a morphism
from the spectrum $T$ of
a discrete valuation ring
such that the image of $T\to S$
is the same as that of $\{f(x),t\}$.
Since the formation of $R\Phi_f{\cal F}$ commutes
with arbitrary base change
by Proposition \ref{prM}.1,
the base change morphism
$R\Phi_f{\cal F}_{x\gets t}
\to R\Phi_{f_T}({\cal F}|_{X_T})_{x\gets t}$
is an isomorphism.
The complex
$R\Phi_{f_T}({\cal F}|_{X_T})$ 
is a perverse sheaf 
by the assumption (P)
and by the theorem of Gabber
\cite[Corollaire 4.6]{au}.
Since 
$R\Phi_{f_T}({\cal F}|_{X_T})$ 
vanishes outside the closed fiber $Z_s$,
this implies that
the complex $R\Phi_f{\cal F}$
is acyclic except at degree $0$.
\qed}
\medskip

The condition {\rm (P)}
is satisfied
if $f\colon X\to S$ is smooth of relative dimension $d$
and
${\cal F}=j_!{\cal G}[d]$
for the open immersion $j\colon U\to X$
of the complement $U=X\sm D$
of a Cartier divisor $D$
and a locally constant constructible 
sheaf ${\cal G}$
of free $\Lambda$-modules on $U$.

\subsection{Semi-continuity of the Swan conductor}\label{sssc}

In this subsection,
we assume that $\Lambda=\Lambda_0$ is a field
for simplicity.
If $\Lambda$ is not a field,
the same results hold without modifications
for constructible complexes of
$\Lambda$-modules
of finite tor-dimension,
by taking 
$\otimes^L_{\Lambda}\Lambda_0$.

We reformulate the main result of
Deligne-Laumon in \cite{DL}
in Proposition \ref{prDL} below.
Let $f\colon X\to S$
be a flat morphism of relative dimension $1$
and let $Z\subset X$ be a closed subscheme.
Assume that 
$X\sm Z$ is smooth over $S$
and that $Z$ is quasi-finite over $S$.
Let ${\cal F}$ be a constructible
complex of $\Lambda$-modules on $X$ such that
the restrictions of the cohomology sheaves
on $X\sm Z$ are locally constant.

Let $s\to S$ be a geometric point
such that the residue
field is an {\em algebraic closure}
of the residue field of the image of $s$
in $S$.
For a geometric point 
$x$ of $Z$ above $s$,
the normalization of
the strict localization
$X_{s,(x)}$ is the finite disjoint
union $\amalg_i{\rm Spec}\ {\cal O}_{K_i}$ 
where ${\cal O}_{K_i}$ are strictly local discrete valuation
rings with algebraically
closed residue field $k(s)$.
Let $\bar\eta_i$ denote
the geometric point
defined by a separable closure $\bar K_i$ of
the fraction field $K_i$.
For a $\Lambda$-representation $V$
of the absolute Galois group $G_{K_i}
={\rm Gal}(\bar K_i/K_i)$,
the Swan conductor
${\rm Sw}_{K_i}V\in {\mathbf N}$
is defined \cite{DL}
and the total dimension is
defined as the sum
$\dim{\rm tot}_{K_i}V
=\dim V
+{\rm Sw}_{K_i}V$.
%We will recall the definition of the total dimension more precisely in a more general situation in Section \ref{ssram}.

The stalk
${\cal H}^q({\cal F})_{{\bar \eta}_i}$
for each integer $q$
defines a $\Lambda$-representation
of the absolute Galois group $G_{K_i}$
and hence the total dimension
$\dim{\rm tot}_{K_i}{\cal F}_{{\bar \eta}_i}$
is defined as the alternating sum
$\sum_q(-1)^q
\dim{\rm tot}_{K_i}{\cal H}^q({\cal F})_{{\bar \eta}_i}$.
We define the Artin conductor by
\begin{equation}
a_x({\cal F}|_{X_s})
=
\sum_i
\dim{\rm tot}_{K_i}{\cal F}_{{\bar \eta}_i}
-
\dim{\cal F}_x.
\label{eqax}
\end{equation}
We define 
a function $\varphi_{\cal F}$
on $X$ supported on $Z$ by
\begin{equation}
\varphi_{\cal F}(x)=
a_x({\cal F}|_{X_s})
\label{eqphis}
\end{equation}
for $s=f(x)$.
The derivative $\delta(\varphi_{\cal F})$
on $X\overset\gets\times_SS$
is defined by (\ref{eqdel}).

\begin{pr}[{\cite[Th\'eor\`eme 2.1.1]{DL}}]\label{prDL}
Let $S$ be a noetherian scheme
and 
$f\colon X\to S$
be a flat morphism of relative dimension $1$.
Let $Z\subset X$
be a closed subscheme quasi-finite
over $S$
such that $U=X\sm Z$ is smooth over $S$.
Let ${\cal F}$ be a constructible
complex of $\Lambda$-modules 
%of finite tor-dimension
on $X$ such that
the restrictions of cohomology sheaves
on $X\sm Z$ are locally constant.

{\rm 1.}
The objects $R\Psi_f{\cal F}$
and $R\Phi_f{\cal F}$
are constructible
and their formations commutes
with any base change.
The function
$\varphi_{\cal F}$ {\rm (\ref{eqphis})}
satisfies
\begin{equation}
\dim R\Phi_f{\cal F}_{x\gets t}=
\delta(\varphi_{\cal F})(x\gets t)
\label{eqDL}
\end{equation}
and is constructible.

{\rm 2.}
Assume ${\cal F}=j_!{\cal G}[1]$
for the open immersion 
$j\colon U=X\sm Z\to X$
and a locally constant constructible 
sheaf ${\cal G}$ %of free $\Lambda$-modules
on $U$
and that $Z$ is flat over $S$.
Then, we have $\delta(\varphi_{\cal F})\geqq 0$.
The function $\varphi_{\cal F}$
is {\em flat} over $S$
if and only if
$f\colon X\to S$ is universally locally acyclic
relatively to ${\cal F}=j_!{\cal G}[1]$.
\end{pr}

\proof{We sketch and/or recall an outline of
proof with some modifications.

{\rm 1.}
The constructibility of $R\Psi_f{\cal F}$
and $R\Phi_f{\cal F}$
and the commutativity with base
change follow
from Proposition \ref{prM}.1
and the local acyclicity of smooth morphism.

By devissage, the proof of 
(\ref{eqDL}) is reduced to 
the case where ${\cal F}=j_!{\cal G}[1]$
for the open immersion 
$j\colon U=X\sm Z\to X$
and a locally constant sheaf ${\cal G}$ on $U$.
By the commutativity with base
change,
the equality (\ref{eqDL})
is reduced to the case where $S$
is the spectrum of a complete
discrete valuation ring
with algebraically closed residue field.
Further by base change
and the normalization, we may assume 
that $X$ is normal
and that its generic fiber is smooth.
By devissage, we may assume
that $Z$ is flat over $S$.

In this case, (\ref{eqDL})
was first proved in \cite{DL}, under an
extra assumption that $X$ is smooth,
by constructing a
good compactification using a deformation argument.
Later it was reproved together with a generalization
 in \cite[Remark (4.6)]{Kato}
using the semi-stable reduction theorem
of curves
without using the deformation argument.

Since $R\Phi_f{\cal F}$ is constructible,
the equality (\ref{eqDL}) implies that
the function $\delta(\varphi_{\cal F})$
on $Z\overset\gets\times_SS$
is constructible.
Hence $\varphi_{\cal F}$
on $Z$ is constructible
by Lemma \ref{lmscZ}.2.

2.
The complex $R\Phi_f{\cal F}$
is acyclic except at degree $0$
by Lemma \ref{lmP}.
Hence the assertions follow
from the equality (\ref{eqDL}) 
and Corollary \ref{corsc}.2.
\qed}

\begin{cor}\label{corDL}
Assume further that $Z$ is finite 
and flat over $S$
and that ${\cal F}=j_!{\cal G}$
for a locally constant sheaf ${\cal G}$
on $U$.
Then, 
the function $f_*\varphi_{\cal F}$ 
{\rm (\ref{eqf*})} on $S$
is lower semi-continuous.
The function $f_*\varphi_{\cal F}$
is locally constant if and only if
$f\colon X\to S$ is universally locally acyclic
relatively to ${\cal F}=j_!{\cal G}$.
\end{cor}

\proof{
It follows from Proposition \ref{prDL},
Lemma \ref{lmscZ}.4 and Corollary \ref{corsc}.2.
The lower semi-continuity
replaces the upper semi-continuity
because of the shift $[1]$
in Proposition \ref{prDL}.2.
\qed}
\medskip

We give a slight generalization
of Proposition \ref{prDL}
using vanishing topos.
Let 
\begin{equation}\xymatrix{
Z\ar[r]^{\subset}&
X\ar[rr]^f\ar[rd]_p&&
Y\ar[ld]^g\\
&&S}
\label{eqfpg}
\end{equation}
be a commutative diagram of morphisms
of finite type of noetherian schemes
such that
$g\colon Y\to S$
is flat of relative dimension $1$
and that $Z\subset X$ 
is a closed subscheme
quasi-finite over $S$.
For geometric point
$x\to X$,
we set $y=f(x)$ and $s=p(x)$
and define 
$T_{(x)}\subset Y_{(y)}$ to be the image
of the finite scheme
$Z\times_XX_{(x)}$ over $S_{(s)}$
by $f_{(x)}\colon X_{(x)}\to Y_{(y)}$.
Assume that, for every geometric point
$x\to X$,
the complement
$Y_{(y)}\sm T_{(x)}$
is essentially smooth over $S_{(s)}$.

Let ${\cal K}$
be an object of
$D^b_c(X\overset \gets\times_YY)$
such that, 
for every geometric point
$x\to X$, the restrictions of
cohomology sheaves on
$Y_{(y)}\sm T_{(x)}
\subset Y_{(y)}=x\overset \gets\times_YY$ 
are locally constant.
Then, similarly as (\ref{eqphis}),
we define 
a function $\psi_{\cal K}$
on $X\overset\gets\times_YY$ by
\begin{equation}
\psi_{\cal K}(x\gets w)=
a_w({\cal K}|_{Y_{(y)}\times_{S_{(s)}t}})
\label{eqphi2}
\end{equation}
where $y=f(x),s=p(x)$ and $t=g(w)$
with {\em algebraically closed} residue field $k(t)$.
We also define a function $\delta(\psi_{\cal K})$
on $X\overset\gets\times_SS$ by
(\ref{eqdelY}).

\begin{pr}\label{prc}
Let the notation be as above.
Let ${\cal K}$
be an object of
$D_c^b(X\overset \gets\times_YY)$
and $x\gets t$ be a point of
$X\overset \gets\times_SS$.
Set $y=f(x)$ and $s=p(x)$
and assume that the restriction of
cohomology sheaf ${\cal H}^q{\cal K}$
on $Y_{(y)}\sm T_{(x)}$ 
is locally constant for every $q$.
Then, we have
\begin{equation}
\dim R\overset\gets g_*{\cal K}_{x\gets t}-
\dim R\overset\gets g_*{\cal K}_{x\gets s}=
\delta(\psi_{\cal K})(x\gets t).
\label{eqDL2}
\end{equation}
\end{pr}

\proof{
By the  canonical isomorphisms
$R\overset\gets g_*{\cal K}_{x\gets t}
\to
R\Gamma(Y_{(y)}\times_{S_{(s)}}S_{(t)},
{\cal K}|_{Y_{(y)}\times_{S_{(s)}}S_{(t)}})$
and
$R\overset\gets g_*{\cal K}_{x\gets s}
\to
R\Gamma(Y_{(y)},
{\cal K}|_{Y_{(y)}})={\cal K}_y$
(\ref{eqMt}),
we obtain a distinguished triangle
$\to
R\overset\gets g_*{\cal K}_{x\gets s}
\to
R\overset\gets g_*{\cal K}_{x\gets t}
\to
R\Phi_{g_{(y)}}({\cal K}|_{Y_{(y)}})_{y\gets t}
\to $.
Hence it follows
from Proposition \ref{prDL}.1.
\qed}
\medskip

In fact, (\ref{eqDL})
is a special case of (\ref{eqDL3}) below
where $X=Y$.

\begin{cor}\label{corc}
We keep the notation in Proposition {\rm \ref{prc}}.
Let ${\cal F}$ be an object of
$D^b_c(X)$ such that ${\cal K}=
R\Psi_f{\cal F}$ is an object of
$D^b_c(X\overset\gets\times_YY)$.
Assume that ${\cal K}$
and a point $x\gets t$ of
$X\overset \gets\times_SS$
satisfies the condition in Proposition 
{\rm \ref{prc}}.
Then, we have
\begin{equation}
\dim R\Phi_p{\cal F}_{x\gets t}
=
\delta(\psi_{\cal K})(x\gets t).
\label{eqDL3}
\end{equation}
\end{cor}

\proof{
By the isomorphisms
$R\Psi_p{\cal F}
\to R\overset\gets g_*{\cal K}$
and $R\overset\gets g_*{\cal K}_{x\gets s}
\to {\cal K}_y\to {\cal F}_x$,
we obtain a distinguished triangle
$\to R\overset\gets g_*{\cal K}_{x\gets s}
\to R\overset\gets g_*{\cal K}_{x\gets t}
\to R\Phi_p{\cal F}_{x\gets t}\to$.
Hence it follows from
(\ref{eqDL2}).
\qed}
\medskip

We consider the diagram (\ref{eqfpg})
satisfying the condition there
and assume further that $g\colon Y\to S$ 
is smooth.
Let ${\cal F}$
be an object of
$D^b_c(X)$ and 
assume that $p\colon X\to S$
is locally acyclic
relatively to ${\cal F}$
and that the restriction of $f\colon X\to Y$
to the complement $X\sm Z$ is
locally acyclic
relatively to the restriction of ${\cal F}$.

We define a function 
$\varphi_{{\cal F},f}$
on $Z$ as follows.
For a geometric point
$x$ of $Z$,
set $y=f(x)$ and $s=p(x)$.
We regard $s$ as a geometric point
of $S$ such that the residue
field is an {\em algebraic closure}
of the residue field of the image of $s$
in $S$.
The base change 
$f_s\colon X_s\to Y_s$
of $f\colon X\to Y$
is a morphism to a smooth curve
over the algebraically closed field $k(s)$.
The strict localization
$Y_{s,(y)}$ is ${\rm Spec}\ {\cal O}_{K_y}$
for a strictly local
discrete valuation ring ${\cal O}_{K_y}$
with an algebraically closed
residue field $k(y)=k(s)$
since $Y_s$ is a smooth curve over $s$.
The cohomology of the 
stalk of the vanishing cycles complex
$\phi_x({\cal F}|_{X_s},f_s)$
define $\Lambda$-representations
of the absolute Galois group
$G_{K_y}$ and hence the total dimension
$\dim{\rm tot}_y
\phi_x({\cal F}|_{X_s},f_s)$
is defined as the alternating sum.
Similarly as (\ref{eqphis}),
we define a function 
$\varphi_{{\cal F},f}$
on $Z$ by
\begin{equation}
\varphi_{{\cal F},f}(x)
=
\dim{\rm tot}_y
\phi_x({\cal F}|_{X_s},f_s)
\label{eqphKf}
\end{equation}

\begin{pr}\label{prMsc}
Let 
\begin{equation*}\xymatrix{
Z\ar[r]^{\subset}&
X\ar[rr]^f\ar[rd]_p&&
Y\ar[ld]^g\\
&&S}
\leqno{\rm (\ref{eqfpg})}
\end{equation*}
be a commutative diagram of morphisms
of finite type of noetherian schemes
such that
$g\colon Y\to S$
is {\em smooth} of relative dimension $1$
and that $Z\subset X$ 
is a closed subscheme
quasi-finite over $S$.

Let ${\cal F}$
be an object of
$D^b_c(X)$ and 
assume that $p\colon X\to S$ is
locally acyclic
relatively to  ${\cal F}$
and that the restriction of $f\colon X\to Y$
to the complement $X\sm Z$ is
locally acyclic
relatively to the restriction of ${\cal F}$.
Then, the function 
$\varphi_{{\cal F},f}$ 
{\rm (\ref{eqphKf})} on $Z$
is constructible and {\em flat} over $S$.
If $Z$ is \'etale over $S$,
it is locally constant.
\end{pr}

\proof{
Let $x$ be a geometric point of
$Z$ and let $y=f(x)$ and
$s=p(x)$ be its images.
We regard $s$ as a geometric point
of $S$ such that the residue
field is an {\em algebraic closure}
of the residue field of the image of $s$
in $S$.
Let $Y_{s,{(y)}}={\rm Spec}\ {\cal O}_{K_y}$
be the strict localization
of the geometric fiber and let
$x\gets u$ be the point
of $X\overset\gets \times_YY$
defined by a separable closure
of $K_y$.
The complex $R\Phi_f{\cal F}$
is constructible
and its construction commutes
with base change
by Proposition \ref{prM}.1.
Hence, we have a canonical
isomorphism
$\phi_x({\cal F}|_{X_s},f_s)
\to
R\Phi_f{\cal F}_{x\gets u}$
and
\begin{equation}
\varphi_{{\cal F},f}(x)
=
\dim{\rm tot}_y
\phi_x({\cal F}|_{X_s},f_s)
=
\dim{\rm tot}_yR\Phi_f{\cal F}_{x\gets u}.
\label{eqphw}
\end{equation}

We apply Proposition \ref{prc} to
${\cal K}=R\Psi_f{\cal F}$.
The assumption in Proposition \ref{prc} that
${\cal H}^q{\cal K}
=R^q\Psi_f{\cal F}$
on
$Y_{(y)}\sm T_{(x)}\subset
Y_{(y)}=x\overset \gets\times_YY$
is locally constant for every $q$
is satisfied
for every geometric point $x$ of $X$
by Corollary \ref{corsl}.
Hence the function 
$\psi_{\cal K}$  (\ref{eqphi2})
for ${\cal K}=R\Psi_f{\cal F}$
is defined as a function
on $X\overset\gets\times_YY$.

In order to apply Lemma \ref{lmdelp}, 
we show
$$\psi_{\cal K}(x\gets w)
=\sum_{z\in Z_{(x)}\times_{Y_{(y)}}w}
\varphi_{{\cal F},f}(z)$$
for a point $x\gets w$
of $Z\overset\gets\times_YY$
such that $w$ is supported on 
the image $T_{(x)}\subset Y_{(y)}$
of $Z_{(x)}$.
By the assumption that $Y\to S$ is
smooth, the fiber
$Y_{(w)}\times_{S_{(t)}}t$
for $t=p(w)$
is the spectrum of a discrete valuation ring.
Let $v$ be its geometric generic point
regarded as a geometric point of
$Y_{(w)}$.

By (\ref{eqphi2}) and (\ref{eqax}), we have
$$\psi_{\cal K}(x\gets w)
=
\dim{\rm tot}_w(R\Psi_f{\cal F}_{x\gets v})
-
\dim (R\Psi_f{\cal F}_{x\gets w}).$$
We apply Proposition \ref{prM}.2
to $f\colon X\to Y$
and specializations $y\gets w\gets v$
to compute the right hand side.
Then, the distinguished triangle (\ref{eqDel})
implies that the right hand side equals
$\sum_{z\in Z_{(x)}\times_{Y_{(y)}}w}
\dim{\rm tot}_w(R\Phi_f{\cal F}_{z\gets v})$.
By (\ref{eqphw})
applied for $z$,
it is further equal to $
\sum_{z\in Z_{(x)}\times_{Y_{(y)}}w}
\varphi_{{\cal F},f}(z)
$ as required.

Therefore, by applying Lemma \ref{lmdelp}, 
we obtain
$\delta(\psi_{\cal K})=
\delta(\varphi_{{\cal F},f})$
as functions on $Z\overset\gets\times_SS$.
Since $R\Phi_p{\cal F}=0$,
the function $\psi_{\cal K}$ is flat over $S$
by (\ref{eqDL3}).
Hence, 
the function $\varphi_{{\cal F},f}$ is also flat
over $S$.
Since it is flat over $S$, 
the function $\varphi_{{\cal F},f}$ is 
constructible by Lemma \ref{lmscZ}.2.

If $Z$ is \'etale over $S$, the function
$\varphi_{{\cal F},f}$ is locally constant by
Lemma \ref{lmscZ}.1.
\qed}

\section{Closed conical subsets
on the cotangent bundle}\label{sCc}

In this preliminary section,
we recall and study some notions
introduced in \cite{Be}
related to closed conical subsets
of the cotangent bundles.

\subsection{$C$-transversality}\label{ssCt}

We say that a closed subset 
$C\subset E$ of a vector bundle
$E$ over $X$ is {\em conical} if it is stable
under the action of 
the multiplicative group ${\mathbf G}_m$.
Equivalently,
it is defined by a graded ideal ${\cal I}$
of the graded algebra $S^\bullet_{{\cal O}_X} 
{\cal E}^\vee$,
if the vector bundle $E={\mathbf V}({\cal E})$ 
is associated to a locally free 
${\cal O}_X$-module ${\cal E}$.
For a closed conical subset
$C\subset E$,
we call its intersection
$B$ with the $0$-section
regarded as a closed
subset of $X$ the {\em base} of $C$.

Let 
${\mathbf P}(C)={\rm Proj}_X
(S^\bullet {\cal E}^\vee/{\cal I})
\subset {\mathbf P}(E)={\rm Proj}_X
S^\bullet {\cal E}^\vee$
denote the projectivization.
The projectivization
${\mathbf P}(C)$ is empty if and
only if $C$ is a subset of the
$0$-section.
The projectivization 
${\mathbf P}(C)$ itself does not
determine $C$ 
but the pair with the base
$B$ determines $C$ uniquely.

We study the intersection of 
a closed conical subset 
with the inverse image of the 0-section
by a morphism of vector bundles.
For a morphism $X\to Y$ 
of finite type of locally noetherian schemes
and a closed subset $Z\subset X$,
we say that the restriction of $X\to Y$ on $Z$
is {\em finite} (resp.\ {\em proper}) 
if for every closed subscheme
structure or equivalently
for the reduced closed subscheme
structure on $Z$,
the induced morphism
$Z\to Y$ is finite (resp.\ proper).

\begin{lm}[{\rm \cite[Lemma 1.2 (ii)]{Be}}]\label{lmnc}
Let $E\to F$ be 
a morphism of vector bundles
over a locally noetherian scheme $X$.
For a closed conical subset $C
\subset E$,
the following
conditions are equivalent:

{\rm (1)}
The intersection of $C$
with the inverse image $K$ of
the $0$-section of $F$ by
$E\to F$ is a subset of
the $0$-section of $E$.

{\rm (2)}
The restriction of
$E\to F$ on $C$ is finite.
\end{lm}

\proof{
(1)$\Rightarrow$(2):
Since the question is local on 
$X$, we may assume $X$ is
affine.
Then the assertion
follows from the elementary Lemma below.

(2)$\Rightarrow$(1):
The intersection $C\cap K\subset E$ is 
a closed conical subset
finite over $X$.
Hence, it is a subset of the
0-section.
\qed}

\begin{lm}\label{lmRI}
Let $R=\bigoplus_{n\geqq 0}
R_n$ be a graded ring
and $I=\bigoplus_{n> 0}
R_n\subset R$ be the graded ideal.
Let $M=\bigoplus_{n\geqq 0}
M_n$ be a graded $R$-module.
If the $R/I$-module
$M/IM$ is finitely generated,
then the $R$-module
$M$ is finitely generated.
\end{lm}

\proof{
Let $x_1,\ldots,x_m\in M$
be a lifting of a system of generators 
of $M/IM$ consisting
of homogeneous elements
and let $N\subset M$ be 
the sub $R$-module
generated by $x_1,\ldots,x_m$.
By induction on $n$,
the morphism
$N/I^nN\to M/I^nM$ is 
surjective for every $n\geqq 0$.
Thus, we have $M_n\subset N$
for every $n\geqq 0$
and the assertion follows.
\qed}
\medskip

In the rest of this article,
$k$ denotes a field of characteristic 
$p\geqq 0$
and $X$ denotes a smooth scheme over $k$,
unless otherwise stated.
The cotangent bundle
$T^*X$ is the covariant
vector bundle over $X$
associated to the locally free ${\cal O}_X$-module
$\Omega^1_{X/k}$.
The $0$-section of $T^*X$ is 
identified with the conormal bundle
$T^*_XX$ of $X\subset X$.

For a closed conical subset
$C\subset T^*X$, 
we define the condition for
a morphism coming into $X$
to be $C$-transversal.

\begin{df}[{\rm \cite[1.2]{Be}}]\label{dfCh}
Let $X$ be a smooth scheme 
over a field $k$
and let $C\subset T^*X$ be a closed conical subset of 
the cotangent bundle.
Let $h\colon W\to X$ be a 
morphism of smooth schemes over $k$.
Define  
\begin{equation}
h^*C=W\times_XC\subset W\times_XT^*X
\label{eqhC}
\end{equation}
to be the pull-back of $C$ and
let $K\subset W\times_XT^*X$ 
be the inverse image of the $0$-section 
by the canonical morphism
$dh\colon W\times_XT^*X\to T^*W$.

{\rm 1.}
For a point $w\in W$,
we say that $h\colon W\to X$
is {\em $C$-transversal} at $w$
if the fiber $(h^*C\cap K)\times_Ww$
of the intersection 
is a subset of the $0$-section
$W\times_XT^*_XX\subset W\times_XT^*X$.

We say that $h\colon W\to X$
is {\em $C$-transversal}
if the intersection $h^*C\cap K$
is a subset of the $0$-section
$W\times_XT^*_XX\subset W\times_XT^*X$.

{\rm 2}.
If $h\colon W\to X$
is $C$-transversal,
we define a closed conical subset
\begin{equation}
h^{\circ}C\subset T^*W
\label{eqhoC}
\end{equation}
to be the image of
$h^*C$ by $W\times_XT^*X\to T^*W$.
\end{df}

By Lemma {\rm \ref{lmnc}},
if $h\colon W\to X$ is $C$-transversal,
then $h^{\circ}C$ is a closed conical subset
of $T^*W$.

\begin{lm}\label{lmCh}
Let $X$ be a smooth scheme 
over a field $k$
and let $C\subset T^*X$ be a closed conical subset of 
the cotangent bundle.
Let $h\colon W\to X$ be a morphism 
of smooth schemes over $k$.

{\rm 1.}
If $h$ is smooth,
then $h$ is $C$-transversal
and the canonical morphism
$h^{*}C\to h^{o}C$ is an isomorphism.

{\rm 2.}
If $C\subset T^{*}_{X}X$ is 
a subset of the $0$-section,
then $h$ is $C$-transversal.

{\rm 3.}
{\rm (cf.\ \cite[Lemma 2.2 (i)]{Be})}
For a morphism $g\colon V\to W$
of smooth schemes over $k$,
the following conditions are equivalent:

{\rm (1)}
$h$ is $C$-transversal
on a neighborhood of
$g(V)\subset W$ and 
$g\colon V\to W$ is
$h^{\circ}C$-transversal.

{\rm (2)}
The composition $h\circ g\colon V\to X$
is $C$-transversal.

{\rm 4.} {\rm (\cite[Lemma 1.2 (i)]{Be})}
The subset of $W$ consisting of
points $w\in W$ where
$h\colon W\to X$ is $C$-transversal
is an open subset of $W$.

{\rm 5.}
Let $D=\bigcup_{i=1}^mD_i$
be a divisor with simple normal crossings
of $X$ relatively to $X\to {\rm Spec}\ k$
and let 
\begin{equation}
C_D=\bigcup_{I\subset\{1,\ldots\,m\}}
T^*_{D_I}X\subset T^*X
\label{eqCD}
\end{equation}
be the union of the conormal bundles
of the intersections
$D_I=\bigcap_{i\in I}D_i$
of irreducible components
for all subsets $I\subset\{1,\ldots,m\}$
of indices, including $D_\varnothing=X$.
Then, $h\colon W\to X$ is $C$-transversal
if and only if
$h^*D=D\times_XW$
is a divisor with simple normal crossings
relatively to $W\to {\rm Spec}\ k$
and 
$h^*D_i\subset W$ are smooth divisors
for $i=1,\ldots,m$.
\end{lm}

The assertion 3 implies that
$C$-transversal morphisms
have similar properties as
\'etale morphisms.
For ${\cal F}$-transversal morphisms 
introduced in Definition \ref{dfprpg},
properties corresponding to 1-3
will be proved 
in Lemma \ref{lmprpg}.

\proof{
1. 
Since
the canonical morphism
$dh\colon W\times_XT^*X\to T^*W$
is a closed immersion,
the intersection $h^*C\cap K$
is a subset of the $0$-section 
$K\subset W\times_XT^*X$ and
the morphism $h^*C\to h^{\circ}C$ is
an isomorphism.

2. 
Since $h^*C\subset
W\times_XT^*X$ is a subset of 
the $0$-section,
the assertion follows.

3.
Since 
$V\times_XT^*X\to T^*V$
is the composition
$V\times_XT^*X\to V\times_WT^*W\to T^*V$,
the condition (2) is equivalent to
that both 
the intersection of
$(hg)^*C=g^*(h^*C)
\subset V\times_XT^*X$ with the
inverse image of the $0$-section by
$V\times_XT^*X\to V\times_WT^*W$
and the intersection of the image of
$(hg)^*C$ in
$V\times_WT^*W$
with the
inverse image of the $0$-section by
$V\times_WT^*W\to T^*V$
are subsets of the $0$-sections.
By 4 below,
the first condition means that
$h$ is $C$-transversal
on a neighborhood of
$g(V)$ and then the second means that
$g\colon V\to W$ is
$h^{\circ}C$-transversal.

4.
The complement of the subset is
the image of the closed subset
${\mathbf P}(h^*C\cap K)\subset
{\mathbf P}(W\times_XT^*X)$
of a projective space bundle over $W$.

5. The $C$-transversality is equivalent to
the injectivity of
$W\times_XT^*_{D_I}X\to
T^*W$ for all $I\subset \{1,\ldots, m\}$.
Hence the assertion follows.
\qed}
\medskip

For a closed conical subset
$C\subset T^*X$, 
we define the condition for
a morphism going out of $X$
to be $C$-transversal.

\begin{df}[{\rm \cite[1.2]{Be}}]\label{dfCf}
Let $X$ be a smooth scheme 
over a field $k$
and let $C\subset T^*X$ be a closed conical subset 
of the cotangent bundle.
Let $f\colon X\to Y$ 
be a morphism of
smooth schemes over $k$.

{\rm 1.}
For a point $x\in X$,
we say that $f\colon X\to Y$
is {\em $C$-transversal} at $x$
if the fiber $df^{-1}(C)\times_Xx$
of the inverse image
of $C$ by the canonical morphism
$df\colon X\times_YT^*Y\to T^*X$
is a subset of the $0$-section
$X\times_YT^*_YY\subset X\times_YT^*Y$.

We say that $f\colon X\to Y$
is {\em $C$-transversal}
if the inverse image $df^{-1}(C)$
of $C$ by the canonical morphism
$df\colon X\times_YT^*Y\to T^*X$
is a subset of the $0$-section
$X\times_YT^*_YY\subset X\times_YT^*Y$.

{\rm 2.}
We say that a pair
of morphisms $h\colon W\to X$ and
$f\colon W\to Y$ of smooth schemes over $k$
is {\em $C$-transversal}
if $h\colon W\to X$ 
is $C$-transversal
and if
$f\colon W\to Y$ 
is $h^{\circ}C$-transversal.
\end{df}

\begin{lm}\label{lmCf}
Let $X$ be a smooth scheme 
over a field $k$
and let $C\subset T^*X$ be a closed conical subset 
of the cotangent bundle.
Let $f\colon X\to Y$ 
be a morphism of
smooth schemes over $k$.

{\rm 1.}
For $Y={\rm Spec}\ k$,
the canonical
morphism $f\colon X\to {\rm Spec}\ k$ 
is $C$-transversal.

{\rm 2.}
Assume that $C$ is the $0$-section
$T^*_XX\subset T^*X$.
Then, $f$ is $C$-transversal
if and only if $f$ is smooth.

{\rm 3.}
Assume that $f$ is \'etale.
Then, $f$ is $C$-transversal
if and only if 
$C$ is a subset of the $0$-section
$T^*_XX\subset T^*X$.

{\rm 4.}
Assume that $f\colon X\to Y$ 
is $C$-transversal and let
$g\colon Y\to Z$ be a smooth morphism.
Then, the composition
$g\circ f\colon X\to Z$ 
is $C$-transversal.

{\rm 5.} {\rm (\cite[Lemma 1.2 (i)]{Be})}
The subset of $X$ consisting of
points $x\in X$ where
$f\colon X\to Y$ is $C$-transversal
is an open subset of $Y$.

{\rm 6.}
Assume that
$f\colon X\to Y$ 
is $C$-transversal.
Then, the morphism 
$f\colon X\to Y$ is smooth
on a neighborhood of the
base $B$ of $C$.

{\rm 7.}
Assume that $Y$ is a curve
and let $x\in X$ be a point.
Then, $f\colon X\to Y$ is 
{\em not} $C$-transversal at
$x$ if and only if
the image of the fiber
$(X\times_YT^*Y)\times_Xx$
by the canonical morphism
$df\colon X\times_YT^*Y\to T^*X$
is a subset of the fiber $C\times_Xx$.

{\rm 8.}
Let $D=\bigcup_{i=1}^mD_i$
be a divisor with simple normal crossings
of $X$ relatively to $X\to {\rm Spec}\ k$
and let 
$C_D=\bigcup_{I\subset\{1,\ldots\,m\}}
T^*_{D_I}X\subset T^*X$
{\rm (\ref{eqCD})}
be the union of the conormal bundles
of the intersections
$D_I=\bigcap_{i\in I}D_i$
of irreducible components
for all subsets $I\subset\{1,\ldots,m\}$
of indices.
Then, $f\colon X\to Y$ is $C$-transversal
if and only if
$D\subset X$
has simple normal crossings
relatively to $f\colon X\to Y$.

{\rm 9}.
Let $h\colon W\to X$ and $f\colon W\to Y$
be morphisms of smooth schemes
over $k$. Then, the following conditions
are equivalent:

{\rm (1)}
The pair $(h,f)$ is $C$-transversal.

{\rm (2)}
The morphism $(h,f)\colon W\to X\times Y$ 
is $C\times T^*Y$-transversal.
\end{lm}

\medskip
In next subsection,
we will see that
the property 1 is related to
the generic local acyclicity
\cite[Corollaire 2.16]{TF}.
The property 2 is related to
the local acyclicity of smooth morphism
(see also Lemma \ref{lmlcst}.1).
The property 3 is related to
the characterization of
locally constant sheaves
\cite[Proposition 2.11]{cst}
(see also Lemma \ref{lmlcst}.3).
The property 4 is related to
\cite[Corollaire 2.7]{app}.

\proof{
1.
Since $T^*Y=0$ for $Y={\rm Spec}\ k$,
the assertion follows.

2.
Assume that
$C$ is the $0$-section
$T^*_XX\subset T^*X$.
Then, $f\colon X\to Y$
is $C$-transversal
if and only if the canonical morphism
$X\times_YT^*Y\to T^*X$
is an injection.
Hence, this is equivalent to that
$f$ is smooth.

3.
Assume that $f$ is \'etale.
Then, since
$X\times_YT^*Y\to T^*X$
is an isomorphism,
the morphism
$f$ is $C$-transversal if and only if
$C$ is a subset of the $0$-section
$T^*_XX\subset T^*X$.

4.
Assume that 
$g\colon Y\to Z$ is smooth.
Then, since $X\times_ZT^*Z\to 
X\times_YT^*Y$ is an injection,
the $C$-transversality for $f$ 
implies that for $g\circ f$.

5.
The complement of the subset is
the image of the closed subset
${\mathbf P}(df^{-1}C)\subset
{\mathbf P}(X\times_YT^*Y)$
of a projective space bundle over $X$.

6.
If $f\colon X\to Y$ 
is $C$-transversal,
then $df\colon X\times_YT^*Y\to T^*X$
is an injection on a neighborhood of $B$
and hence
$f\colon X\to Y$ is smooth
on a neighborhood of $B$.

7.
Since the fiber
$(X\times_YT^*Y)\times_Xx$
is a line and the inverse image
of the fiber $C\times_Xx$
is its conical subset,
the inverse image of $C\times_Xx$
is not the subset of the $0$-section
if and only if it is equal to
the line
$(X\times_YT^*Y)\times_Xx$
itself.

8.
The $C$-transversality 
is equivalent to the injectivity of
$D_I\times_YT^*Y\to T^*D_I$
for all subsets $I\subset \{1,\ldots,m\}$.
Hence, it is 
equivalent to the smoothness of
$D_I\to Y$
for $I\subset \{1,\ldots,m\}$.
Thus the assertion follows.

9.
The conditions (1) and (2)
are rephrased as conditions
on elements
$(\alpha,\beta)
\in {\rm Ker}(W\times_XT^*X\oplus
W\times_YT^*Y
\to T^*W)_w$
in the fiber at $w\in W$ 
of the kernel as follows:

(1) If $\alpha$ is contained in the inverse
image of $C$ and $\beta=0$,
then $\alpha=0$.
Further, if
$\alpha$ is contained in the inverse
image of $C$ then $\beta=0$.

(2) 
If $\alpha$ is contained in the inverse
image of $C$ then $(\alpha,\beta)=(0,0)$.

\noindent
Thus, the conditions are equivalent.
\qed}

\begin{df}[{\rm \cite[1.2]{Be}}]\label{dff!C}
Let $f\colon X\to Y$ be a
morphism of smooth schemes over $k$
and $C\subset T^*X$ be a 
closed conical subset.
Assume that $f$ is {\em proper}
on the base $B$ of $C$.
We define a closed conical subset 
\begin{equation}
f_{\circ}C \subset T^*Y
\label{eqf!C}
\end{equation}
to be the closure of the image
by the first arrow
$$\begin{CD}T^*Y
@<<< 
X\times_YT^*Y@>>> T^*X
\end{CD}$$
of the inverse image of $C$
by the second arrow.
\end{df}

\begin{lm}\label{lmncf}
Let $f\colon X\to Y$
be a morphism
of smooth schemes over $k$.
Let $C\subset T^*X$ be a closed conical 
subset.
Assume that $f$ is {\em proper}
on the base $B$ of $C$.
Then, for a morphism
$g\colon Y\to Z$ of smooth schemes over $k$,
the following conditions
are equivalent:

{\rm (1)}
The morphism $g$ is $f_\circ C$-transversal.

{\rm (2)}
The composition $gf\colon X\to Z$ 
is $C$-transversal.
\end{lm}

\proof{
We consider the commutative diagram
$$\begin{CD}
@.
X\times_ZT^*Z@>>>
X\times_YT^*Y@>>> T^*X\\
@.@VVV @VVV\\
T^*Z@<<<
Y\times_ZT^*Z@>>> T^*Y.
\end{CD}$$
The condition (1) 
(resp.\ (2)) means that
the subset in 
$T^*Z$ obtained by taking
inverse images and images
in the diagram
starting from
$C\subset T^*X$
via $T^*Y$ 
(resp.\ via $X\times_ZT^*Z$) is a subset of
the $0$-section.
They are equivalent
since the square is cartesian.
\qed}%\medskip

\begin{lm}\label{lmChf}
Let $C\subset T^*X$ be a closed conical subset
and let 
\begin{equation}
\begin{CD}
X@<h<<W\\
@VfVV @VVgV\\
Y@<i<<Z
\end{CD}
\label{eqfghi}
\end{equation}
be a cartesian diagram of
smooth schemes over $k$.

{\rm 1.}
Assume that the horizontal arrows
are regular immersions of the same codimension.
Then, 
the following conditions are equivalent:

{\rm (1)}
The morphism 
$f\colon X\to Y$ 
is $C$-transversal on a neighborhood of $W$.

{\rm (2)}
The pair of the immersion
$h\colon W\to X$
and the morphism 
$g\colon W\to Z$ 
is $C$-transversal.

{\rm 2.}
Assume that $f\colon X\to Y$
is smooth and
$C$-transversal.
Then, 
the pair of $h\colon W\to X$ and 
$g\colon W\to Z$ is $C$-transversal.

{\rm 3.}
Assume that $f$ is proper on
the base $B$ of $C$
and is flat on a neighborhood
of $B$.
Assume that the horizontal arrows
are immersions.
Then, the following conditions
are equivalent:

{\rm (1)}
$i$ is $f_{\circ}C$-transversal,

{\rm (2)}
$h$ is $C$-transversal.

Further, if these equivalent conditions
are satisfied, 
we have $i^{\circ}f_{\circ}C
=g_{\circ}h^{\circ}C$.
\end{lm}

\proof{
1.
We have a commutative diagram
\begin{equation}
\begin{CD}
0@>>> T^*_WX@>>>
W\times_XT^*X@>>> T^*W@>>>0\\
@.@AAA @AAA @AAA@.\\
0@>>> W\times_ZT^*_ZY@>>>
W\times_YT^*Y@>>> W\times_ZT^*Z@>>>0
\end{CD}
\label{eqTWZ}
\end{equation}
of exact sequences
of vector bundles on $W$.
The condition (1) is equivalent to
that the inverse image in 
$W\times_YT^*Y$ of 
$h^*C\subset W\times_XT^*X$
by the middle vertical arrow is
a subset of the $0$-section,
by Lemma \ref{lmCh}.4.
The condition (2) is equivalent to
that the inverse image in 
$T^*_WX$ of 
$h^*C\subset W\times_XT^*X$
by the upper left horizontal arrow 
and the inverse image in 
$W\times_ZT^*Z$ of 
$h^{\circ}C\subset T^*W$
by the right vertical arrow 
are subsets of the $0$-sections.
Since the left vertical arrow 
is an isomorphism,
the conditions (1) and (2) 
are equivalent.

2.
First, we assume $Z\to Y$ is smooth.
Then $h$ is smooth and $C$-transversal.
Further 
the arrows in the right square
of (\ref{eqTWZ})
are injections
and the square is cartesian.
Thus the inverse image of
$h^{\circ}C\subset 
W\times_XT^*X
\subset T^*W$
in $W\times_ZT^*Z$
is the same as the inverse image
in $W\times_YT^*Y$
and hence $g$ is $h^{\circ}C$-transversal.

If $Z\to Y$ is an immersion,
the assertion follows from 1.
In general,
it suffices to decompose 
$Z\to Y$ as the composition
$Z\to Z\times Y\to Y$
of the graph and the projection
by Lemma \ref{lmCh}.3.

3. 
By replacing $X$ by a neighborhood
of $B$, we may assume
that $f$ is flat.
The condition (1) (reap.\ (2))
means that the inverse image of
$C\subset T^*X$
in $T^*_WX$ 
(resp.\ its image in $T^*_ZY$)
is a subset of the $0$-section.
Since the left vertical arrow
in (\ref{eqTWZ}) is an isomorphism,
these conditions are equivalent.

The subset
$i^{\circ}f_{\circ}C$
(resp.\ $g_{\circ}h^{\circ}C$)
of $T^*Z$
is the image of
the subset of
$W\times_ZT^*Z$
obtained from
$h^*C\subset
W\times_XT^*X$
by taking inverse images
and images
in the right square of (\ref{eqTWZ})
via $W\times_YT^*Y$
(resp.\ via $T^*W$).
Since the square is cartesian,
they are equal.
\qed}
%\medskip

\subsection{The universal family of
hyperplane sections}\label{ssRn}

Assume that $X$ 
smooth over a field $k$
is quasi-projective
and let ${\cal L}$ be an ample invertible
${\cal O}_X$-module.
Let $E$ be a  $k$-vector space of finite
dimension and
$E\to \Gamma(X,{\cal L})$ 
be a $k$-linear mapping inducing
a surjection
$E\otimes_k{\cal O}_X\to {\cal L}$
and an immersion 
\begin{equation}
i\colon X\to {\mathbf P}
={\mathbf P}(E^\vee)=
{\rm Proj}\ S^\bullet E.
\label{eqiXP}
\end{equation}
We use a contra-Grothendieck notation
for a projective space
${\mathbf P}(E)(k)=(E\sm \{0\})/k^\times$.

The dual ${\mathbf P}^\vee=
{\mathbf P}(E)$
of ${\mathbf P}$ parametrizes
hyperplanes $H\subset {\mathbf P}$.
The universal hyperplane
$Q=\{(x,H)\mid x\in H\}
\subset 
{\mathbf P}\times
{\mathbf P}^\vee$
is defined by
the identity ${\rm id}\in
{\rm End}(E)$
regarded as a section 
$F\in \Gamma(
{\mathbf P}\times
{\mathbf P}^\vee,
{\cal O}(1,1))
=E\otimes E^\vee$.
By the canonical injection
$\Omega^1_{{\mathbf P}/k}(1)
\to E\otimes {\cal O}_{\mathbf P}$,
the universal hyperplane
$Q$ is identified 
with the covariant projective
space bundle
${\mathbf P}(T^*{\mathbf P})$
associated to the cotangent bundle
$T^*{\mathbf P}$.
The image of the conormal bundle
$T^*_Q(
{\mathbf P}\times {\mathbf P}^\vee)
\to
Q\times_
{{\mathbf P}\times {\mathbf P}^\vee}
T^*({\mathbf P}\times {\mathbf P}^\vee)
\to
Q\times_
{\mathbf P}
T^*{\mathbf P}$
by the projection
is identified
with the universal sub line bundle
of the pull-back 
$Q\times_
{\mathbf P}T^*{\mathbf P}$
on 
$Q={\mathbf P}(T^*{\mathbf P})$.

The fibered product
$X\times_{\mathbf P}Q
={\mathbf P}(X\times_{\mathbf P}T^*{\mathbf P})$
is the intersection
of  $X\times {\mathbf P}^\vee$
with $Q$
in ${\mathbf P}\times {\mathbf P}^\vee$
and is the universal family
of hyperplane sections.
We consider
the universal family of hyperplane sections
\begin{equation}
\begin{CD}
X@<p<<X\times_{\mathbf P}Q
@>{p^\vee}>> {\mathbf P}^\vee={\mathbf P}(E).
\end{CD}
\label{eqhsfb}
\end{equation}

Let $C\subset T^*X$ be a
closed conical subset.
Define a closed conical subset
$\widetilde C\subset
X\times_{\mathbf P}T^*{\mathbf P}$
to be the pull-back of $C$ by
the surjection
$X\times_{\mathbf P}T^*{\mathbf P}
\to T^*X$
and let
\begin{equation}
{\mathbf P}(\widetilde C)\subset
{\mathbf P}(X\times_{\mathbf P}T^*{\mathbf P})
=
X\times_{\mathbf P}Q
\label{eqPS}
\end{equation}
be the projectivization.
The subset
${\mathbf P}({\widetilde C})
\subset
X\times_{\mathbf P}Q
\subset X\times {\mathbf P}^\vee$
consists of the points $(x,H)$
such that the fiber
$T^*_{X\times_{\mathbf P}Q}
(X\times{\mathbf P}^\vee)
\times_{X\times_{\mathbf P}Q}(x,H)
\subset
(X\times_{\mathbf P}T^*{\mathbf P})\times_Xx$ 
is a subset of $\widetilde C$
since the conormal
bundle $T^*_{X\times_{\mathbf P}Q}
(X\times{\mathbf P}^\vee)
\subset 
X\times_{\mathbf P} T^*{\mathbf P}$
is the universal sub line bundle
on the projective bundle
$X\times_{\mathbf P}Q=
{\mathbf P}(X\times_{\mathbf P} T^*{\mathbf P})$.
If $C=T^*_XX$ is the $0$-section,
the lifting $\widetilde C$
is the conormal bundle
$T^*_X{\mathbf P}$.
Further if $i\colon X\to {\mathbf P}$
is a closed immersion,
the image $p^\vee(
{\mathbf P}(\widetilde C))
\subset {\mathbf P}^\vee$
is the dual variety of $X$.

If $V\subset {\mathbf P}$ is
an open subscheme such that
$i\colon X\to {\mathbf P}$
induces a closed immersion
$i^{\circ}\colon X\to V$,
then
$\widetilde C\subset
X\times_{\mathbf P}T^*{\mathbf P}
=X\times_VT^*V
\subset T^*V$
is identified with $i^{\circ}_{\circ}C$
(\ref{eqf!C}).

\begin{lm}\label{lmPC}
Let $C\subset T^*X$ be a
closed conical subset.
The complement 
$X\times_{\mathbf P}Q
\sm {\mathbf P}(\widetilde C)$
is the largest open subset
where the pair
$X\gets X\times_{\mathbf P}Q
\to {\mathbf P}^\vee$
is $C$-transversal.
\end{lm}

\proof{
By Lemma \ref{lmCf}.9,
the largest open subset
$U\subset X\times _{\mathbf P}Q$
where the pair
$(p,p^\vee)$
is $C$-transversal equals
the largest open subset
where $(p,p^\vee)\colon 
X\times_{\mathbf P}Q\to X\times
{\mathbf P}^\vee$
is $C\times T^*{\mathbf P}^\vee$-transversal.
The kernel 
$L={\rm Ker}\bigl(
(X\times_{\mathbf P}Q)\times_XT^*X
\oplus
(X\times_{\mathbf P}Q)\times_
{{\mathbf P}^\vee}T^*{\mathbf P}^\vee
\to 
T^*(X\times_{\mathbf P}Q)\bigr)$
is canonically identified with
the restriction of
the universal sub line bundle of
$T^*{\mathbf P}$ on
$Q={\mathbf P}(T^*{\mathbf P})$.
Hence,  $U$ is the complement
of ${\mathbf P}(\widetilde C)$.
\qed}
\medskip

We consider a similar 
but slightly different situation.
Let $X$ be a smooth
scheme over a field $k$
and $h\colon X\to {\mathbf P}
={\mathbf P}^n$
be an immersion over $k$
to a projective space.
Let $Q={\mathbf P}(T^*{\mathbf P})
\subset {\mathbf P}
\times {\mathbf P}^\vee$
be the universal family
of hyperplanes
and consider the
commutative diagram 
\begin{equation}
\xymatrix{
X\times_{\mathbf P}Q
\ar[r]\ar[d]_p
\ar[rrrd]^{\!\! p^{\vee}}
& 
Q\ar[rrd]^{{\bms p}^{\vee}}
\ar[d]_{\!\!\!\!\!\! {\bms p}}&&
\\
X\ar[r]^h&
{\mathbf P}
&&{\mathbf P}^\vee}
\label{eqXPh}
\end{equation}
with cartesian square.
We have $Q
={\mathbf P}(X\times_{\mathbf P}T^*{\mathbf P})$.

\begin{lm}\label{lmh!}
Let ${\mathbf P}={\mathbf P}^n$
be a projective space 
and let ${\mathbf P}^\vee$
be the dual projective space.
Let $C^\vee
\subset T^*{\mathbf P}^\vee$
be a closed conical subset 
such that every irreducible
component is of dimension $n$.
Define a closed conical subset $C
\subset T^*{\mathbf P}$ by
$C={\bm p}_{\circ}
{\bm p}^{\vee \circ}C^\vee$.
Then, every irreducible
component of $C$ is of dimension $n$.
\end{lm}

\proof{
It suffices to consider
the case where $C^\vee$ is irreducible.
We have ${\mathbf P}(C)=
{\mathbf P}(C^\vee)$
in $Q={\mathbf P}(T^*{\mathbf P})
={\mathbf P}(T^*{\mathbf P}^\vee)$.
Unless the base of $C^\vee$
is finite,
$C$ contains the $0$-section
$T^*_{\mathbf P}{\mathbf P}$.
If the base of $C^\vee$
consists of a closed point
corresponding to a hyperplane
$H\subset {\mathbf P}$,
we have $C=T^*_H{\mathbf P}$.
Thus the assertion follows.
\qed}

\begin{pr}\label{prhCt}
Let the notation be as in
{\rm (\ref{eqXPh})} and
let $C^\vee
\subset T^*{\mathbf P}^\vee$
be a closed conical subset.
Define a closed conical subset 
$C\subset T^*{\mathbf P}$
by $C={\bm p}_{\circ}{\bm p}^{\vee \circ}C^\vee$
and 
its projectivization
${\mathbf P}(C)
\subset
{\mathbf P}(T^*{\mathbf P})$.
Then, the complement
$X\times_{\mathbf P}Q
\sm (X\times_{\mathbf P}Q)
\cap {\mathbf P}(C)$
is the largest open subset $U
\subset X\times_{\mathbf P}Q$
where the pair
$(p^\vee, p)$
is $C^\vee$-transversal.
\end{pr}

\proof{
The proof is similar to that of
Lemma \ref{lmPC}.1.
By Lemma \ref{lmCf}.9,
the largest open subset
$U\subset X\times _{\mathbf P}Q$
where the pair
$(p^\vee, p)$
is $C^\vee$-transversal equals
the largest open subset
where $(p,p^\vee)\colon 
X\times_{\mathbf P}Q\to X\times
{\mathbf P}^\vee$
is $T^*X\times C^\vee$-transversal.
The kernel 
$L={\rm Ker}\bigl(
(X\times_{\mathbf P}Q)\times_XT^*X
\oplus
(X\times_{\mathbf P}Q)\times_
{{\mathbf P}^\vee}T^*{\mathbf P}^\vee
\to 
T^*(X\times_{\mathbf P}Q)\bigr)$
is canonically identified with
the restriction of
the universal sub line bundle of
$T^*{\mathbf P}^\vee$ on
$Q={\mathbf P}(T^*{\mathbf P}^\vee)$.
Hence,  $U$ is the complement
of the intersection
$(X\times_{\mathbf P}Q)
\cap {\mathbf P}(C^\vee)$.
Since $L$ is also canonically identified with
the restriction of
the universal sub line bundle
of $T^*{\mathbf P}$ on
$Q={\mathbf P}(T^*{\mathbf P})$,
we have
${\mathbf P}(C^\vee)
={\mathbf P}(C)$
and the assertion follows.
\qed}
\medskip

Let $\Delta_X\subset X\times_{\mathbf P}Q$
be the sub projective space
bundle
\begin{equation}
\Delta_X={\mathbf P}(T^*_X{\mathbf P})
\subset
{\mathbf P}(X\times_{\mathbf P}T^*{\mathbf P})
=X\times_{\mathbf P}Q.
\label{eqRX}
\end{equation}

\begin{cor}\label{corhCt}
Let the notation be as in
Proposition {\rm \ref{prhCt}}.

{\rm 1.}
The following conditions are
equivalent:

{\rm (1)}
The immersion $h\colon X\to
{\mathbf P}$ is $C$-transversal.

{\rm (2)}
The pair $(p^\vee, p)$
is $C^\vee$-transversal
on a neighborhood
of $\Delta_X\subset X\times_{\mathbf P}Q$.

{\rm 2.}
Assume that
the immersion $h\colon X\to
{\mathbf P}$ is $C$-transversal.
Then, $p^\vee\colon X\times_{\mathbf P}Q
\to {\mathbf P}^\vee$ is
$C^\vee$-transversal and we have
$h^{\circ}C=p_{\circ}p^{\vee \circ}C^\vee$.
\end{cor}

\proof{
1.
By Proposition \ref{prhCt},
the condition (2) is equivalent to
${\mathbf P}(C)\cap
\Delta_X=\varnothing$.
This is equivalent to
that $T^*_X{\mathbf P}
\cap C$
is a subset of the $0$-section
and hence to (1).

2.
By Lemma \ref{lmChf}.3
applied to $C={\bm p}_{\circ}
{\bm p}^{\vee \circ}C^\vee$,
the assumption that
$h\colon X\to {\mathbf P}$
is $C$-transversal
implies that the immersion
$i\colon X\times_{\mathbf P}Q
\to Q$ is ${\bm p}^{\vee \circ}C^\vee$-transversal
and $h^{\circ}C=
h^{\circ}{\bm p}_{\circ}
{\bm p}^{\vee \circ}C^\vee=
p_{\circ}i^{\circ}{\bm p}^{\vee \circ}C^\vee$.
Hence by Lemma \ref{lmCh}.3,
$p^\vee\colon X\times_{\mathbf P}Q
\to {\mathbf P}^\vee$ is
$C^\vee$-transversal and
we have
$i^{\circ}{\bm p}^{\vee \circ}C^\vee=
p^{\vee \circ}C^\vee$.
Thus the assertion is proved.
\qed}

%\medskip

\begin{pr}\label{prfp}
Let the notation be as in
Proposition {\rm \ref{prhCt}} and
let $f\colon X\to {\mathbf A}^1$
be a smooth morphism over $k$.
Define sub vector bundles
of $X\times_{\mathbf P}T^*{\mathbf P}$
by the cartesian diagram
$$\begin{CD}
T^*_X{\mathbf P}
@>{\subset}>> V@>{\subset}>> 
X\times_{\mathbf P}T^*{\mathbf P}\\
@VVV@VVV@VVV\\
T^*_XX
@>{\subset}>> 
X\times_{{\mathbf A}^1}T^*{\mathbf A}^1
@>{\subset}>> T^*X
\end{CD}$$
and closed subsets
$\Delta_X={\mathbf P}(T^*_X{\mathbf P})
\subset \Delta_f={\mathbf P}(V)
\subset 
X\times_{\mathbf P}Q
={\mathbf P}(X\times_{\mathbf P}T^*{\mathbf P})$
to be the associated projective subspace bundles.
Then, the complement 
$X\times_{\mathbf P}Q
\sm (\Delta_f\cap {\mathbf P}(C))$
is the largest open subset
$U\subset X\times_{\mathbf P}Q$
where the pair
$(p^\vee,f\circ p)$
is $C^\vee$-transversal.
\end{pr}

\proof{
By Lemma \ref{lmCf}.9,
the largest open subset
$U\subset X\times_{\mathbf P}Q$
where the pair
$(p^\vee,f\circ p)$
is $C^\vee$-transversal
equals
the largest open subset
where $(f\circ p,p^\vee)\colon
X\times_{\mathbf P}Q
\to {\mathbf A}^1\times{\mathbf P}^\vee$
is $T^*{\mathbf A}^1
\times C^\vee$-transversal.
Hence, similarly as in the proof
of Proposition \ref{prhCt},
$U\subset X\times_{\mathbf P}Q$ 
is the complement of
$\Delta_f\cap 
((X\times_{\mathbf P}Q)
\cap {\mathbf P}(C^\vee))=
\Delta_f\cap 
{\mathbf P}(C)$.
\qed}

\begin{cor}\label{corfp}
Let the notation be as in
Propositions {\rm \ref{prhCt}}
and {\rm \ref{prfp}} and assume that
$h\colon X\to {\mathbf P}$
is $C$-transversal.
Then,
the composition $f\circ p\colon 
X\times_{\mathbf P}Q\to X\to {\mathbf A}^1$
is $p^{\vee \circ}C^\vee$-transversal
on the complement
of the intersection 
with $\Delta_f\cap {\mathbf P}(C)$.
Further, the intersection
$\Delta_f\cap {\mathbf P}(C)
\subset 
X\times_{\mathbf P}Q$
is finite over the complement
of the largest open subset
of $X$ where
$f\colon X\to {\mathbf A}^1$
is $h^{\circ}C$-transversal.
\end{cor}

\proof{
The first assertion is clear from
Proposition \ref{prfp}.

We show the second assertion.
By Corollary \ref{corhCt}.2,
$p^\vee\colon X\times_{\mathbf P}Q
\to {\mathbf P}^\vee$ is
$C^\vee$-transversal and
we have $h^{\circ}C=
p_{\circ}p^{\vee \circ}C^\vee$.
Hence by Lemma \ref{lmncf},
if $f\colon X\to {\mathbf A}^1$
is $h^{\circ}C$-transversal,
then $fp\colon X\times_{\mathbf P}Q
\to {\mathbf A}^1$ is
$p^{\vee \circ} C^\vee$-transversal.
Thus, by Proposition \ref{prfp},
the intersection
$\Delta_f\cap {\mathbf P}(C)$
is a subset of the inverse
image of
the complement
of the largest open subset
of $X$ where
$f\colon X\to {\mathbf A}^1$
is $h^{\circ}C$-transversal.
Therefore, it suffices to show that
$\Delta_f\cap {\mathbf P}(C)$
is finite over $X$.

By the assumption that
$h\colon X\to {\mathbf P}$
is $C$-transversal,
the morphism
$p\colon X\times_{\mathbf P}Q
\to X$ is $p^{\vee \circ}C^\vee$-transversal
on a neighborhood of $\Delta_X$
by Corollary \ref{corhCt}.1.
Since $f\colon X\to {\mathbf A}^1$
is smooth, the composition
$fp\colon X\to {\mathbf A}^1$
is $p^{\vee \circ}C^\vee$-transversal
on a neighborhood of $\Delta_X$
by Lemma \ref{lmCf}.4.
Hence, by Proposition \ref{prfp},
the intersection
$\Delta_f\cap {\mathbf P}(C)$
does not meet $\Delta_X$.
Since $\Delta_X$ is a hyperplane bundle
of a projective space bundle
$\Delta_f$ over $X$,
the intersection
$\Delta_f\cap {\mathbf P}(C)$
is finite over $X$.
\qed}

\subsection{Image of the projectivization}

We further recall from \cite{Be}
a definition and properties.

\begin{df}[{\rm cf.\ \cite[4.2]{Be}}]\label{dfwi}
Let $f\colon X\to S$ be a morphism
of separated schemes of finite type over a field $k$
and $Y,Z\subset X$ be closed subsets.
We say that $Y$ and $Z$ {\em well intersect}
with respect to $f$ if we have
\begin{equation}
\dim (Y\times_SZ\sm Y\times_XZ)
\leqq \dim Y+\dim Z-\dim S.
\label{eqdfwi}
\end{equation}
\end{df}

We slightly modified the
original definition by replacing
the equality by an inequality.

\begin{lm}\label{lmJ}
Let $f\colon X\to S$ be a morphism of
separated schemes of finite type
over a field $k$
and $Y,Z\subset X$ be closed subsets.

{\rm 1.}
Let $g\colon P\to S$ be
a morphism of
separated schemes of finite type
over a field $k$
and $X\to P$ be an immersion
over $k$.
Then, $Y$ and $Z$ well intersect
with respect to $f$
if and only if they
well intersect
with respect to $g$.

{\rm 2.}
Let $S'\to S$ be a faithfully flat morphism of
relative constant dimension of
separated schemes of finite type
over a field $k$
and let $f'$ and $Y',Z'\subset X'$
denote the base changes of $f$ and $Y,Z\subset X$
by $S'\to S$.
Then, $Y$ and $Z$ well-intersects
with respect to $f$ if and only if
$Y'$ and $Z'$ well-intersects
with respect to $f'$.

{\rm 3.}
Let $X\to U\to T$ be a morphism of
separated schemes of finite type
over a field $k$
and $V,W\subset U$ be closed subsets.
Assume that $X\to S\times T$
is an immersion and that
$Y,Z\subset X$ are inverse images of
$V,W\subset U$.
Assume further that
the morphisms
$X\to S$ and $X\to U$
and the base change
\begin{equation}
(X\times_SX)
\times_{(T\times T)}
(T\times T)^{\circ}\to
(U\times U)
\times_{(T\times T)}
(T\times T)^{\circ}
=U\times U
\sm U\times_T U
\label{eqXXUU}
\end{equation}
of the morphism
$X\times_SX\to U\times U$
by the open immersion
$(T\times T)^{\circ} =
T\times T\sm T
\to T\times T$ of the complement
of the diagonal are faithfully flat of relative
constant dimensions.
Then, $Y$ and $Z$ well-intersects
with respect to $f$.
\end{lm}

\proof{1.
Since $X\to P$ is an immersion,
we have $Y\times_XZ=
Y\times_PZ$
and the assertion follows.

2.
By the assumption, every term in (\ref{eqdfwi})
with $'$ equals the corresponding term
without $'$ added the relative dimension
and the assertion follows.

3.
By the assumption that
$X\to U$ is faithfully 
flat of relative constant dimension,
we have
\begin{equation}
\dim Y=\dim V+(\dim X-\dim U),\quad
\dim Z=\dim W+(\dim X-\dim U).
\label{eqYV}
\end{equation}
By the assumption that $X\to S\times T$
is an immersion,
the complement of the diagonal
$X\times_SX\sm X$ 
equals the base change 
$(X\times_SX)
\times_{(T\times T)}
(T\times T)^{\circ}$
and 
$Y\times_SZ\sm Y\times_XZ$
is the inverse image of
$V\times W\sm
V\times_T W$
by the morphism
(\ref{eqXXUU}).
Since the morphism
(\ref{eqXXUU}) is faithfully flat
of relative constant dimension,
we have
\begin{align}
&\dim (Y\times_SZ\sm Y\times_XZ)
\label{eqVW}
\\
&\leqq
\dim V+\dim W
+(\dim (X\times_SX\sm X)
-\dim (U\times U\sm U\times_TU)).
\nonumber
\end{align}
By the assumption that
$X\to S$ is faithfully flat
of relative constant dimension, 
we have 
$\dim X\times_SX=2\dim X-\dim S$
and the right hand side of
(\ref{eqVW}) equals
$\dim V+\dim W
+2(\dim X-\dim U)-\dim S$
and the assertion follows
by (\ref{eqYV}).
\qed}

\begin{lm}[{\rm cf. \cite[Lemma 4.2]{Be}}]\label{lmwi}
Let $Y,Z\subset X$ be irreducible
closed subsets.
Assume that 
$Y$ and $Z$ well intersects
with respect to $f$
and we have
$\dim Y,\dim Z<\dim S$.

{\rm 1.}
If $Y=Z$,
then $Y\to \overline{f(Y)}$ 
is generically radicial.

{\rm 2.}
If $Y\not \subset Z$,
then we have
$f(Y)\not \subset
\overline{f(Y)\cap f(Z)}$
and 
$f(Y)\not \subset f(Z)$.
\end{lm}

\proof{
1.
We have
$\dim (Y\times_SY\sm Y)
\leqq 2\dim Y-\dim S< \dim Y.$

2.
By 1, we have
$\dim Y=\dim \overline{f(Y)}$.
We have $\dim \overline{f(Y)\cap f(Z)}
\leqq
\dim Y\times_SZ
\leqq
\max(
\dim (Y\times_SZ\sm Y\times_XZ),
\dim Y\times_XZ)$.
Since $Y$ and $Z$ assumed to
well intersect,
we have
$\dim (Y\times_SZ\sm Y\times_XZ)
\leqq \dim Y+\dim Z-\dim S
< \dim Y$.
If $Y\not\subset Z$,
we have
$\dim Y\times_XZ
< \dim Y$.
Thus, we have
$\dim \overline{f(Y)\cap f(Z)}
<\dim Y=\dim \overline{f(Y)}$.

If we had $f(Y)\subset f(Z)$,
we would have
$f(Y)\subset
\overline{f(Y)\cap f(Z)}$
and a contradition.
\qed}
\medskip

For a $k$-vector space $E$
of finite dimension
and a $k$-linear mapping
$E\to \Gamma(X,{\cal L})$
defining an immersion
$X\to {\mathbf P}={\mathbf P}(E^\vee)$,
we consider the following condition:
\begin{itemize}
\item[{\rm (E)}]
For every pair of distinct closed
points $u\neq v$ of the base change
$X_{\bar k}$ to an algebraic closure
$\bar k$ of $k$,
the composition
\begin{equation}
E\to \Gamma(X,{\cal L})
\otimes_k\bar k
\to 
{\cal L}_u/{\mathfrak m}_u^2{\cal L}_u
\oplus
{\cal L}_v/{\mathfrak m}_v^2{\cal L}_v
\label{eqE}
\end{equation}
is a surjection.
\end{itemize}
For an integer $d\geqq 1$,
the linear mapping
$S^dE
\to\Gamma(X,{\cal L}^{\otimes d})$ 
of the symmetric power
defines an immersion
$X\to {\mathbf P}(S^dE^{\vee})$.

\begin{lm}\label{lmloc}
Let $X$ be a quasi-projective
scheme over a field
$k$ and
${\cal L}$ be an ample
invertible ${\cal O}_X$-module.
Assume that $E\to \Gamma(X,{\cal L})$ 
defines an immersion $X\to 
{\mathbf P}={\mathbf P}(E^{\vee})$.
For $d\geqq 3$, 
the symmetric power
$S^dE\to \Gamma(X,{\cal L}^{\otimes d})$
satisfies the condition {\rm (E)} above.
\end{lm}

\proof{
We may assume $k=\bar k,\
X={\mathbf P}^n,\
{\cal L}={\mathcal O}(1),\
E=\Gamma(X,{\cal L})$
and $u=(0,\cdots,0,1),
v=(1,0,\cdots,0)$.
Then, the assertion is clear.
\qed}
\medskip

For an immersion 
$i\colon X\to {\mathbf P}$
to a projective space
and a closed conical
subset $C\subset T^*X$,
we consider the following condition:
\begin{itemize}
\item[(C)]
For every irreducible component
$C_a\subset T^*X$ of $C=\bigcup_aC_a$,
the inverse image
$\widetilde C_a
\subset X\times_{\mathbf P}T^*{\mathbf P}$
is not a subset of the
$0$-section.
\end{itemize}
If the condition (C) is satisfied,
there is a one-to-one correpspondence
between the irreducible components
of $C$ and those of ${\mathbf P}(\widetilde C)$ sending $C_a$ to ${\mathbf P}(\widetilde C_a)$.

\begin{pr}[{\rm \cite[Lemma 4.3]{Be}}]\label{prwi}
Let $X$ be a quasi-projective
smooth scheme of dimension $n$
over a field $k$ and
${\cal L}$ be an ample
invertible ${\cal O}_X$-module.
Let $E$ be a $k$-vector space
of finite dimension and 
$E\to \Gamma(X,{\cal L})$
be a $k$-linear mapping
defining an immersion
$X\to {\mathbf P}={\mathbf P}(E^\vee)$
and 
satisfying the condition {\rm (E)}
before Lemma {\rm \ref{lmloc}}.
Let $C\subset T^*X$
be a closed conical subset
satisfying the condition {\rm (C)} above.
Then, for irreducible
components $C_a,C_b$ of $C$,
the projectivizations
${\mathbf P}(\widetilde C_a),
{\mathbf P}(\widetilde C_b)
\subset {\mathbf P}(X\times_{\mathbf P}
T^*{\mathbf P})=X\times_{\mathbf P}Q$
well intersects with respect to
$p^\vee\colon 
X\times_{\mathbf P}Q\to {\mathbf P}^\vee$.
\end{pr}
\medskip

The proof is essentially the same as
that of \cite[Lemma 4.3]{Be}.
For the sake of completeness,
we rephrase the proof.

\proof{
Let ${\cal I}\subset {\cal O}_{X\times X}$
denote the ideal sheaf defining
the diagonal immersion $X\to X\times X$
and let $P\subset X\times X$ denote
the closed subscheme defined by 
${\cal I}^2$.
The projections $p_1,\ p_2\colon P\to X$
are finite flat of degree $n+1$.
Define a vector bundle $J$
of rank $n+1$ on $X$ to be that associated
to the invertible ${\cal O}_P$-module 
$p_2^*{\cal L}$ 
regarded as a locally free 
${\cal O}_X$-module by $p_1$.
The exact sequence
$0\to \Omega^1_X\to {\cal O}_P
\to {\cal O}_X\to 0$
tensored with $p_2^*{\cal L}$ 
defines an exact sequence
\begin{equation}
\begin{CD}
0@>>>
T^*X\otimes L
@>>>
J
@>>> L@>>>0
\end{CD}
\label{eqLJ}
\end{equation}
of vector bundles on $X$.

Define a morphism $X\times E\to J$
of vector bundles on $X$
to be that induced by the pull-back
by $p_2$ of the canonical morphism
${\cal O}_X\otimes E\to {\cal L}$.
For a closed geometric point $x$,
the fiber of
$X\times E\to J$
is identified with the morphism
$E\to {\cal L}_x/{\mathfrak m}_x^2{\cal L}_x$
induced by the canonical morphism
${\cal O}_X\otimes E\to {\cal L}$.

We regard
${\mathbf P}(\widetilde C_a),
{\mathbf P}(\widetilde C_b)
\subset X\times_{\mathbf P}Q$
as subsets of
$X\times {\mathbf P}^\vee$
as in Lemma \ref{lmJ}.1.
We consider the cartesian diagram
\begin{equation}
\begin{CD}
&X\times {\mathbf P}^\vee
@<<<
X\times E^{\circ}
@>>>
J
@>>> X
\\
&@VVV@VVV\\
{\mathbf P}(E)=&
{\mathbf P}^\vee
@<<<
E^{\circ}&=E\sm \{0\}.
\end{CD}
\label{eqEJ}
\end{equation}
The lower horizontal arrow is
the canonical surjection
and is faithfully flat of
relative dimension 1.
For an irreducible component $C_a$
of $C\subset T^*X$,
define a closed conical
subset $C_a\otimes L
\subset T^*X\otimes L
\subset J$
as a twist of $C_a$ by using
the exact sequence (\ref{eqLJ}).
The inverse image
$\widetilde C_a^{\circ}\subset
X\times E^{\circ}$
of $C_a\otimes L
\subset J$ by the second upper arrow
in (\ref{eqEJ})
equals the inverse image of
${\mathbf P}(\widetilde C_a)
\subset X\times_{\mathbf P}Q
\subset X\times {\mathbf P}^\vee$
by the first upper arrow
by the condition (C).
By Lemma \ref{lmJ}.2
applied to ${\mathbf P}(\widetilde C_a),
{\mathbf P}(\widetilde C_b)
\subset
X\times_{\mathbf P}Q\subset
X\times {\mathbf P}^\vee$
and the base change
$E^{\circ} \to{\mathbf P}^\vee$,
it suffices to show
that
$\widetilde C_a^{\circ}$
and
$\widetilde C_b^{\circ}
\subset
X\times E^{\circ}$
well-intersect with respect to
$X\times E^{\circ} \to E^{\circ}$.

The morphism
$X\times E\to J$
induces a surjection
of vector bundles
$(X\times X)^{\circ} \times E\to 
(J\times J)\times_{(X\times X)}(X\times X)^{\circ}$
on $(X\times X)^{\circ}$
by the condition (E).
Thus, by Lemma \ref{lmJ}.3 applied
to $E^{\circ}\gets X\times E^{\circ}\to J\to X$
taken as $S\gets X\to U\to T$,
the assertion follows.
\qed}

\begin{cor}\label{corrad}
Let $X$ be a quasi-projective
smooth scheme of dimension $n$
over a field $k$ and
${\cal L}$ be an ample
invertible ${\cal O}_X$-module.
Let $E$ be a $k$-vector space
of finite dimension and 
$E\to \Gamma(X,{\cal L})$
be a $k$-linear mapping
defining an immersion
$X\to {\mathbf P}={\mathbf P}(E^\vee)$
satisfying the condition {\rm (E)} 
before Lemma {\rm \ref{lmloc}}.
Let $C\subset T^*X$
be a closed conical subset
satisfying the condition {\rm (C)} before
Proposition {\rm \ref{prwi}}.
Assume that every irreducible component
$C_a$ of
$C=\bigcup_aC_a$ is of dimension $n$. 

{\rm 1.}
For every irreducible component
$C_a$ of $C$,
the restriction
${\mathbf P}(\widetilde C_a)\to {\mathbf P}^\vee$
of $p^\vee\colon X\times_{\mathbf P}Q
\to {\mathbf P}^\vee$
is generically radicial
and the closure 
$D_a=\overline{
p^\vee({\mathbf P}(\widetilde C_a))}
\subset {\mathbf P}^\vee$ is a divisor.

{\rm 2.}
For distinct irreducible
components $C_a\neq C_b$ of $C$,
we have $D_a\neq D_b$.
\end{cor}

\proof{
1.
Since the closed subset
${\mathbf P}(\widetilde C)\to X\times_{\mathbf P}Q$
is of codimension $n$
and $p^\vee\colon X\times_{\mathbf P}Q
\to {\mathbf P}^\vee$
is of relative dimension $n-1$,
we have
$\dim {\mathbf P}(\widetilde C)
=\dim {\mathbf P}^\vee-1$.
Hence it follows from
Proposition \ref{prwi}
and Lemma \ref{lmwi}.1.

2.
Similarly as the proof of 1,
the assertion follows from
Proposition \ref{prwi}
and Lemma \ref{lmwi}.2.
\qed}
%\medskip

\section{Singular support}\label{sss}

We recall definitions and results from
\cite{Be} in Section \ref{ssss}.
We give a description of
singular support of the $0$-extension
of a locally constant sheaf
on the complement of
a divisor with simple normal crossings
under a certain assumption
in Section \ref{ssram}.

\subsection{Singular support}\label{ssss}

We recall definitions and results from
\cite{Be}.
We assume that $X$
is a smooth scheme over a field $k$.

Let $\Lambda$ be a finite local ring
such that the characteristic $\ell$
of the residue field is
invertible in $k$
and let ${\cal F}$ be a constructible complex of
$\Lambda$-modules on 
the \'etale site of $X$.
Recall that a complex 
${\cal F}$ of $\Lambda$-modules
is {\em constructible} if every cohomology sheaf
${\cal H}^q{\cal F}$ is constructible
and if ${\cal H}^q{\cal F}=0$ except
for finitely many $q$.
We say that a constructible complex
${\cal F}$ of $\Lambda$-modules
is {\em locally constant} if every cohomology sheaf
is locally constant.

We say that a constructible complex
${\cal F}$ of $\Lambda$-modules
is of {\em finite tor-dimension} 
if there exists an integer $a$
such that ${\cal H}^q({\cal F}
\otimes^L_{\Lambda}M)=0$
for every $q<a$ 
and for every $\Lambda$-module $M$.
The dualizing complex ${\cal K}_X$
defined as $Rf^!\Lambda$
for the canonical morphism
$f\colon X\to {\rm Spec}\ k$
is canonically isomorphic to
$\Lambda(n)[2n]$ if
every irreducible component of $X$
is of dimension $n$.
For a constructible complex
${\cal F}$ of finite tor-dimension,
the dual $D_X{\cal F}
=R{\cal H}om({\cal F},{\cal K}_X)$
is also constructible of finite tor-dimension.

\begin{df}\label{dfms}
Let $X$ be a smooth scheme
over a field $k$ and
let $C\subset T^*X$ be a closed conical
subset. Let
${\cal F}$ be a constructible
complex of $\Lambda$-modules on $X$.

{\rm 1. \cite[1.3]{Be}}
We say that ${\cal F}$ is
{\em micro-supported} on $C$
if for every $C$-transversal pair
$(f,h)$ of morphisms
$h\colon W\to X$ and $f\colon W\to Y$
of smooth schemes,
the morphism $f\colon W\to Y$
is locally acyclic relatively to 
$h^*{\cal F}$.

{\rm 2. \cite[1.5]{Be}}
We say that ${\cal F}$ is
{\em weakly micro-supported} on $C$
if for every $C$-transversal pair
$(f,j)$ of an \'etale morphism
$j\colon W\to X$ satisfying the
condition {\rm (W)} below
and a morphism
$f\colon W\to Y={\mathbf A}^1_k$,
the morphism $f\colon W\to Y$
is locally acyclic relatively to 
$h^*{\cal F}$:

{\rm (W)}
If $k$ is infinite,
then $j\colon W\to X$
is an open immersion.
If $k$ is finite,
then $j\colon W\to X$
is the composition of
an open immersion $V\to X$
with $W=V_{k'}
=V\times_kk'\to V$
for a finite extension $k'$ of $k$.
\end{df}

\begin{lm}\label{lmmc}
Let $C$ and $C'$ be closed
conical subsets of $T^*X$
and ${\cal F}$ be a constructible
complex of $\Lambda$-modules on $X$.

{\rm 1.}
Assume that $C\subset C'$.
If ${\cal F}$ is micro-supported
(resp.\ weakly micro-supported) on $C$,
then
${\cal F}$ is micro-supported
(resp.\ weakly micro-supported) on $C'$.

{\rm 2.}
The complex ${\cal F}$ is weakly micro-supported
on $C$ if and only if
$C$ contains the images
of the fibers $(U\times_YT^*Y)
\times_Ux$ for all 
\'etale morphisms $j\colon U\to X$
and morphisms
$f\colon U\to Y$ to smooth curves
satisfying the condition {\rm (W)}
in Definition {\rm \ref{dfms}}
and for closed points $x$ of $U$
where $\phi_x(j^*{\cal F},f)\neq 0$.

{\rm 3.}
If ${\cal F}$ is weakly micro-supported on $C$,
then the base $B$ of $C$ contains
the support of ${\cal F}$
as a subset.

{\rm 4.}
We consider the following conditions
for a morphism $h\colon W\to X$
of smooth schemes over $k$:

{\rm (1)}
${\cal F}$ is micro-supported on $C$.

{\rm (2)}
$h^*{\cal F}$ is micro-supported
on $h^{\circ}C$.

If $h$ is $C$-transversal,
we have
{\rm (1)}$\Rightarrow${\rm (2)}.
If $h$ is \'etale and surjective,
conversely we have
{\rm (2)}$\Rightarrow${\rm (1)}.

{\rm 5. (\cite[Lemma 2.1 (ii)]{Be})}
Assume that 
${\cal F}$ is micro-supported on $C$.
Then, for $C$-transversal pair
$(f,h)$ of morphisms
$h\colon W\to X$ and $f\colon W\to Y$
of smooth schemes,
the morphism $f\colon W\to Y$
is {\em universally} locally acyclic relatively to 
$h^*{\cal F}$.

{\rm 6.} {\rm ({\cite[Lemma 2.2 (ii)]{Be}})}
Assume ${\cal F}$ is micro-supported
(resp.\ weakly micro-supported) on $C$
and let $f\colon X\to Y$ be a proper morphism
of smooth schemes over $k$.
Then $Rf_*{\cal F}$ is micro-supported
(resp.\ weakly micro-supported)
on $f_{\circ}C$.

{\rm 7.}
Let $\Lambda_0$
be the residue field of $\Lambda$.
Then,
${\cal F}$ is micro-supported on $C$
if and only if
${\cal F}\otimes_\Lambda^L
\Lambda_0$ is micro-supported on $C$.
\end{lm}

\proof{
1.
If a pair $(f,h)$
of morphisms $h\colon W\to X$
and $f\colon X\to Y$ is $C'$-transversal,
then $(f,h)$
is $C$-transversal.
Hence the assertion follows.

2.
Since the vanishing cycles
complex $\phi(j^*{\cal F},f)$
is constructible,
it follows from Lemma \ref{lmCf}.7.

3.
Define an open subset $U=X\sm B$
to be the complement of
the base $B$ of $C$.
Let $j\colon U\to X$
be the open immersion and
$f\colon U\to {\mathbf A}^1$
be the morphism collapsing to $0$.
Then the pair $(f,j)$
is $C$-transversal.
Hence the morphism
$f\colon U\to {\mathbf A}^1$ is locally acyclic
relatively to $j^*{\cal F}$
and we have ${\cal F}|_U=0$.
Thus the base $B$
contains ${\rm supp}\ {\cal F}$.

4.
{\rm (1)}$\Rightarrow${\rm (2)}:
Assume $h\colon W\to X$ is
$C$-transversal.
If the pair $(f,g)$ of morphisms $g\colon V\to W$ 
and $f\colon V\to Y$ 
of smooth schemes over $k$
is $h^{\circ}C$-transversal,
then the pair $(f,h\circ g)$ is
$C$-transversal by Lemma \ref{lmCh}.3.
Hence the assertion follows.

{\rm (2)}$\Rightarrow${\rm (1)}:
The local acyclicity is an \'etale local condition.

7.
It follows from Lemma \ref{lmctf}.2.
\qed}
\medskip

By Lemma \ref{lmmc}.2,
there exists a smallest closed conical
subset $C\subset T^*X$
such that ${\cal F}$ is weakly micro-supported
on $C$.
This smallest $C$ is called the
{\em weak singular support} and is denoted by
$SS^w{\cal F}$.
For a closed conical subset $C$ of $T^*X$,
a constructible complex
${\cal F}$ is weakly micro-supported on $C$
if and only if
$SS^w{\cal F}$ is a subset of $C$.
By Lemma \ref{lmmc}.2 and 3,
the base of
$SS^w{\cal F}$ equals the support of ${\cal F}$.

If there exists a smallest closed conical
subset $C\subset T^*X$
such that ${\cal F}$ is micro-supported
on $C$, then we call such $C$
the {\em singular support}
of ${\cal F}$ and let it denoted by
$SS{\cal F}$.
If the singular support exists,
we have
$SS^w{\cal F}
\subset
SS{\cal F}$ by definition.
In fact, Theorem \ref{thmss} below
includes the existence of
$SS{\cal F}$ and the equality
$SS^w{\cal F}
=SS{\cal F}$.

\begin{lm}\label{lmlcst}
Let $X$ be a smooth scheme over $k$
and ${\cal F}$ be a constructible
complex of $\Lambda$-modules on $X$.

{\rm 1.}
{\rm (\cite[Lemma 2.1 (iii)]{Be})}
The following conditions are equivalent:

{\rm (1)}
${\cal F}$ is micro-supported on 
the $0$-section $T^*_XX$.

{\rm (2)}
${\cal F}$ is locally constant.

If ${\cal F}$ is locally constant
and if the support of ${\cal F}$ is $X$,
then the singular support
$SS{\cal F}$ exists and 
both $SS{\cal F}$ and $SS^w{\cal F}$
are equal to the $0$-section
$T^*_XX$.

{\rm 2.}
Assume $\dim X=1$.
Let $U\subset X$ be the largest open
subset where the restriction ${\cal F}|_U$
is locally constant.
Then, ${\cal F}$ is micro-supported
on the union
\begin{equation}
T^*_XX\cup \bigcup_{x\in {X\!\sm \!U}}
(T^*X\times_Xx)
\label{eqScv}
\end{equation}
of the $0$-section and
the fibers of the complement.
Further if $X={\rm supp}\ {\cal F}$, 
then the singular support $SS{\cal F}$
and $SS^w{\cal F}$ 
equal the closed conical subset
{\rm (\ref{eqScv})}.

{\rm 3.}
Let $D=\bigcup_{i=1}^mD_i$
be a divisor with simple normal crossings
of $X$ relatively to $k$.
Let ${\cal G}$ be a locally constant 
constructible sheaf
of $\Lambda$-modules
on $U=X\sm D$ tamely ramified along $D$
and let $j\colon U\to X$ be the open immersion.
Then ${\cal F}=j_!{\cal G}$
is micro-supported on the union
\begin{equation}
C_D=\bigcup_{I\subset\{1,\ldots,m\}}
T^*_{D_I}X
\label{eqStame}
\end{equation}
of the conormal bundles
of the intersections
$D_I=\bigcap_{i\in I}D_i$
of irreducible components
for all subsets $I\subset\{1,\ldots,m\}$
of indices.

{\rm 4.}
Assume that ${\cal F}$ is micro-supported
on $C$ and 
let $f\colon X\to Y$ be a $C$-transversal
proper morphism.
Then $Rf_*{\cal F}$ is locally constant.

{\rm 5.}
{\rm ({\cite[Lemma 2.5 (i)]{Be}})}
For a closed immersion
$i\colon X\to Y$
of smooth schemes over $k$,
we have $SS^wi_*{\cal F}=i_{\circ}SS^w{\cal F}$.
\end{lm}

\proof{1.
(1)$\Rightarrow$(2)
Since the identity ${\rm id}_X\colon X\to X$
is $T^*_XX$-transversal by Lemma \ref{lmCf}.3,
the identity ${\rm id}_X\colon X\to X$
is locally acyclic relatively to ${\cal F}$.
This means that ${\cal F}$ is locally constant
by \cite[Proposition 2.11]{cst}.

(2)$\Rightarrow$(1)
Since a $T^*_XX$-transversal
morphism $f\colon X\to Y$
of smooth schemes 
is smooth by Lemma \ref{lmCf}.2,
it follows from 
the local acyclicity of smooth morphism
\cite[Th\'eor\`eme 2.1]{Artin}.

If ${\cal F}$ is locally constant
and if the support of ${\cal F}$ is $X$,
Lemma \ref{lmmc}.3 and
(2)$\Rightarrow$(1)
implies that the singular support
$SS{\cal F}$ exists and 
that the equality $SS{\cal F}=SS^w{\cal F}
=T^*_XX$ holds.

2.
It follows from 1.\ that
${\cal F}$ is micro-supported on the union
in (\ref{eqScv}).

Assume that ${\cal F}$ is 
weakly micro-supported on $C$.
If ${\cal F}$ is not locally constant at $x$,
then $\phi_x({\cal F},{\rm id}_X)\neq 0$.
Hence $C$ contains the fiber $x\times_XT^*X$
by Lemma \ref{lmmc}.2.
If ${\rm supp}\ {\cal F}=X$,
then $C$ contains the $0$-section
$T^*_XX$ by Lemma \ref{lmmc}.3.

3.
By Lemma \ref{lmCh}.5,
it suffices to show that
any $C$-transversal morphism
$f\colon X\to Y$
of smooth schemes is universally
locally acyclic relatively to ${\cal F}$.
By Lemma \ref{lmCf}.8,
the $C$-transversality of
$f\colon X\to Y$ is equivalent to
that the divisor $D$ has
simple normal crossing relatively to $Y$.
Hence, the assertion follows from
a variant of \cite[1.3.3 (i)]{app}
with $Rf_*$ replaced by $f_!$
proved similarly.

4.
Since $f$ is proper,
$Rf_*{\cal F}$ is
micro-supported on $f_{\circ}C$
by Lemma \ref{lmmc}.6.
Since $f$ is $C$-transversal,
$f_{\circ}C\subset T^*Y$ is a subset of the 
$0$-section.
Hence
$Rf_*{\cal F}$ is locally constant
by 1 and Lemma \ref{lmmc}.1.
\qed}

\begin{thm}\label{thmss}
Let $X$ be a smooth scheme of finite type
over a field $k$ and
let ${\cal F}$ be a constructible complex
of $\Lambda$-modules on $X$.

{\rm 1. (\cite[Theorem 1.3 (i)]{Be})}
The singular support
$SS{\cal F}$ exists.

{\rm 2. (\cite[Theorem 1.3 (ii)]{Be})} 
If every irreducible component of $X$
is of dimension $n$,
then
every irreducible component of $SS{\cal F}$
is of dimension $n=\dim X$.

{\rm 3. (\cite[Theorem 1.5]{Be})}
We have $SS{\cal F}=
SS^w{\cal F}$.

{\rm 4. (\cite[Theorem 1.4 (ii)]{Be})}
For an exact sequence
$0\to {\cal F}'\to {\cal F}\to {\cal F}''\to 0$
of perverse sheaves,
we have
\begin{equation}
SS{\cal F}
=
SS{\cal F}'\cup
SS{\cal F}''.
\label{eqsump}
\end{equation}
We have
\begin{equation}
SS{\cal F}
=
\bigcup_q
SS\ ^p\!{\cal H}^q{\cal F}.
\label{eqssu}
\end{equation}
\end{thm}

\begin{cor}\label{corss}

{\rm 1.}
For a closed conical subset $C$ of $T^*X$,
the following conditions are
equivalent:

{\rm (1)}
${\cal F}$ is micro-supported on $C$.

{\rm (2)}
${\cal F}$ is weakly micro-supported on $C$.

{\rm 2. (\cite[Theorem 1.4 (i)]{Be})}
For a smooth morphism $h\colon W\to X$,
we have
$SSh^*{\cal F}=h^{\circ}SS{\cal F}$.

{\rm 3. (\cite[Lemma 2.1 (i)]{Be})}
The base of $SS{\cal F}=SS^w{\cal F}$
equals the support of ${\cal F}$.

{\rm 4.}
Assume that ${\cal F}$
is of finite tor-dimension and
let $\Lambda_0$
be the residue field of $\Lambda$.
Then, for 
${\cal F}_0=
{\cal F}\otimes_\Lambda^L\Lambda_0$
we have
$SS{\cal F}=SS{\cal F}_0$.
\end{cor}

\proof{
1.
By Theorem \ref{thmss}.3,
both conditions are equivalent to that
$SS{\cal F}$ is a subset of $C$.

2.
By Lemma \ref{lmmc}.4 
{\rm (2)}$\Rightarrow${\rm (1)}
and by 1,
we have $SSh^*{\cal F}=h^*SS{\cal F}$
for an \'etale morphism $h\colon W\to X$.
Hence we may assume $W={\mathbf A}^n\times X$.
We have
$SSh^*{\cal F}\subset h^{\circ}SS{\cal F}$
by Lemma \ref{lmmc}.4
{\rm (1)}$\Rightarrow${\rm (2)}.
Hence the $0$-section $i\colon X\to 
W={\mathbf A}^n$ is $SSh^*{\cal F}$-transversal
and further 
we have
$SS{\cal F}\subset i^{\circ}SSh^*{\cal F}$
by Lemma \ref{lmmc}.4
{\rm (1)}$\Rightarrow${\rm (2)}.
Thus we have
$SSh^*{\cal F}=h^*SS{\cal F}$.

4.
It follows from Lemma \ref{lmctf}.2.
\qed}
\medskip

%%%%%%%%%%%%

Let $X$ be a smooth scheme over 
a field $k$ and let
$i\colon X\to {\mathbf P}$
be an immersion to a projective space
as in (\ref{eqiXP}).
Let ${\cal F}$ be
a constructible complex
of $\Lambda$-modules
micro-supported on
a closed conical subset $C
\subset T^*X$
of the cotangent bundle.
If $X$ is connected,
the condition (C) 
for $C=SS{\cal F}$ means
that the support of $i_*{\cal F}$
is not equal to ${\mathbf P}$
by Corollary \ref{corss}.3.

By Lemma \ref{lmPC} 
and Lemma \ref{lmmc}.5,
the morphism $p^\vee \colon
X\times_{\mathbf P}Q
={\mathbf P}(X\times_{\mathbf P}T^*{\mathbf P})
\to {\mathbf P}^\vee$
is universally locally acyclic relatively to
$p^*{\cal F}$ on the complement
of the projectivization 
${\mathbf P}(\widetilde C)\subset
{\mathbf P}(X\times_{\mathbf P}T^*{\mathbf P})$.
More precisely, the following holds.

\begin{thm}[{\rm cf.\
\cite[Theorem 3.2, Lemma 3.3]{Be}}]
\label{thmRn}
Let $i\colon X\to {\mathbf P}$
be an immersion and let
$\widetilde C\subset
X\times_{\mathbf P}T^*{\mathbf P}$
be the inverse image
of the singular support
$C=SS{\cal F}\subset T^*X$
by the surjection
$X\times_{\mathbf P}T^*{\mathbf P}\to
T^*X$.
Then, the complement
$U=X\times_{\mathbf P}Q\sm
{\mathbf P}(\widetilde C)$
of the projectivization 
${\mathbf P}(\widetilde C)\subset
{\mathbf P}(X\times_{\mathbf P}T^*{\mathbf P})
=X\times_{\mathbf P}Q$
is the largest open subset
where $p^\vee \colon
X\times_{\mathbf P}Q
={\mathbf P}(X\times_{\mathbf P}T^*{\mathbf P})
\to {\mathbf P}^\vee$
is universally locally acyclic relatively to
$p^*{\cal F}$.
\end{thm}

\proof{
First, we consider the case where
$X={\mathbf P}$ is a projective space.
Applying \cite[Theorem 3.2]{Be} to 
the Radon transform $R{\cal F}$ \cite{Rn},
we obtain
$E_{p^{\vee}}(p^{*}R^\vee R({\cal F}))=
{\mathbf P}(SSR({\cal F}))$
in the notation loc.\ cit.
Since, $R^\vee R({\cal F})$ is isomorphic to ${\cal F}$
except locally constant sheaf \cite{Rn},
we have
$E_{p^{\vee}}(p^{*}R^\vee R({\cal F}))=
E_{p^{\vee}}(p^{*}{\cal F})$.
By \cite[Lemma 3.3]{Be}, we have
${\mathbf P}(SS{\cal F})=
{\mathbf P}(SSR({\cal F}))
\subset Q$.
Hence the assertion follows
for $X={\mathbf P}$.

Let $V\subset {\mathbf P}$ be
an open subset including $X$ as a closed subset and
let $i^{\circ}\colon X\to V$ be the closed immersion.
Since the assertion is proved for ${\mathbf P}$,
it holds also for an open subscheme $V$.
Since $SS(i^{\circ}_*{\cal F})=\widetilde C$,
the complement
$V\times_{\mathbf P}Q
\sm {\mathbf P}(\widetilde C)$
is the largest open subset
of $V\times_{\mathbf P}Q$
where $p^{\vee}\colon V\times_{\mathbf P}Q
\to {\mathbf P}^{\vee}$ is universally locally 
acyclic relatively to $p^{*}i^{\circ}_*{\cal F}$
by the assertion already proved for $V$.
Since $p^{*}i^{\circ}_*{\cal F}=0$
outside $X\times_{\mathbf P}Q$,
the assertion follows.
\qed}
%\medskip

\begin{cor}\label{corRn}
Let $i\colon X\to {\mathbf P}$
be an immersion and
$C\subset T^*X$
be a closed conical subset
such that the base $B\subset X$
contains the support of ${\cal F}$.
Assume that
there exists a closed subset
$Z\subset X\times_{\mathbf P}Q
={\mathbf P}(X\times_{\mathbf P}T^*{\mathbf P})$
of codimension $> \dim X$
such that 
$p^\vee \colon
X\times_{\mathbf P}Q
={\mathbf P}(X\times_{\mathbf P}T^*{\mathbf P})
\to {\mathbf P}^\vee$
is universally locally acyclic on the complement
of ${\mathbf P}(\widetilde C)\cup Z\subset 
{\mathbf P}(X\times_{\mathbf P}T^*{\mathbf P})$.
Then ${\cal F}$ is micro-supported on $C$.
\end{cor}

\proof{
Let $C_0=SS{\cal F}$
denote the singular support.
By Corollary \ref{corss}.3,
the base $B_0$ of $C_0$
equals the support of ${\cal F}$.
By Theorem \ref{thmRn},
we have ${\mathbf P}(\widetilde C)\cup Z
\supset {\mathbf P}(\widetilde C_0)$.
Since $Z\subset 
X\times_{\mathbf P}Q$ is of codimension $>
\dim X$ by the assumption and
every irreducible component
of ${\mathbf P}(\widetilde C_0)\subset 
X\times_{\mathbf P}Q$ is of codimension $\dim X$
by Theorem \ref{thmss}.2,
we have ${\mathbf P}(\widetilde C)
\supset {\mathbf P}(\widetilde C_0)$.
Thus, we have
$C\sm B\supset C_0\sm B_0$.
By the assumption $B\supset B_0$,
we have $C\supset C_0=SS{\cal F}$
and ${\cal F}$ is micro-supported on $C$
by Lemma \ref{lmmc}.1.
\qed}
\medskip

\begin{thm}[{\rm \cite[Theorem 1.7]{Be}}]\label{thmssD}
Let $X$ be a {\em projective}
smooth scheme of dimension $n$
over a field $k$ and
${\cal L}$ be an ample
invertible ${\cal O}_X$-module.
Let $E$ be a $k$-vector space
of finite dimension and 
$E\to \Gamma(X,{\cal L})$
be a $k$-linear mapping
defining a {\em closed} immersion
$X\to {\mathbf P}={\mathbf P}(E^\vee)$
and 
satisfying the condition {\rm (E)} 
before Lemma {\rm \ref{lmloc}}.
Then, $D=p^\vee({\mathbf P}(SS{\cal F}))
\subset {\mathbf P}^\vee$
is a divisor and the complement
${\mathbf P}^\vee
\sm D$ is the largest open
subset where
$Rp^\vee_*p^*{\cal F}$
is locally constant.
\end{thm}

\proof{
It suffices to apply
\cite[Theorem 1.7]{Be}
to $i_*{\cal F}$.\qed}

\begin{cor}\label{corssD}
Let ${\cal F}$ be a constructible
complex of finite tor-dimension.
Then, we have
\begin{equation}
SS{\cal F}
=
SS D_X{\cal F}.
\label{eqssD}
\end{equation}
\end{cor}

One could deduce 
directly Corollary \ref{corssD}
and its consequence Lemma \ref{lmRf}.4
from the compatibility of vanishing cycles
with duality which seems missing
in the literature.

\proof{
We may assume that $\Lambda=\Lambda_0$
is a finite field by Corollary \ref{corss}.4.
Since the assertion is local on $X$,
we may take a closed immersion
$i\colon X\to U$
and an open immersion
$j\colon U\to {\mathbf P}$
to a projective space.
By Lemma \ref{lmlcst}.5 and
Corollary \ref{corss}.2,
we may assume $X$ is projective.

Let $C=SS{\cal F}$
and $C'=SS D_X{\cal F}$
be the singular supports.
We take a projective embedding
$i\colon X\to {\mathbf P}$
defined by $E$ satisfying
the conditions (E) 
before Lemma {\rm \ref{lmloc}}
and (C) before
Proposition {\rm \ref{prwi}}.
Let $D=p^\vee({\mathbf P}(\widetilde C)),
D'=p^\vee({\mathbf P}(\widetilde C'))
\subset {\mathbf P}^\vee$
be the images.
Then, by Theorem \ref{thmssD},
the complement
${\mathbf P}^\vee
\sm D$
is the largest open subset
where the (shifted and twisted) Radon transform
$Rp^\vee_*p^*{\cal F}$
is locally constant.
Similarly,
the complement
${\mathbf P}^\vee
\sm D'$
is the largest open subset
where 
$Rp^\vee_*p^*D_X{\cal F}$
is locally constant.
Since the Radon transform commutes
with duality up to shift and twist,
we have
$D=D'$. Hence we have $C=C'$
by Corollary \ref{corrad}.2.
\qed}
\medskip

In the proof of Proposition \ref{prperv},
we will use the following fact
proved in the course of 
the proof of \cite[Theorem 1.5]{Be}.
For a line $L\subset {\mathbf P}^\vee$,
the axis $A_L\subset {\mathbf P}$
is a linear subspace of codimension $2$
defined as the intersections
of hyperplanes parametrized by $L$.
On the complement
$X_L^{\circ}=X\sm (X\cap A_L)$,
a canonical morphism
\begin{equation}
\begin{CD}
p_L^{\circ}\colon
X_L^{\circ} @>>> L
\end{CD}
\label{eqpLo}
\end{equation}
is defined by sending a point
$x\in X_L^{\circ}$
to the unique hyperplane $H\in L$
containing $x$.

\begin{lm}[{\rm \cite[4.9 (ii), (iii)]{Be}}]\label{lmRn}
Let $i\colon X\to {\mathbf P}$
be an immersion and
let ${\cal F}$ be
a perverse sheaf of
$\Lambda$-modules on $X$.
After replacing 
the immersion $i\colon X\to {\mathbf P}$
by the composition with
$d$-th Veronese embedding
for $d\geqq 3$,
for an irreducible component 
$C_a$ of $C$,
let $D_a=\overline{p^\vee({\mathbf P}
(\widetilde C_a))}\subset {\mathbf P}^\vee$
denote the closure of the image.
Then there exist
dense open subsets
$D_a^{\circ} \subset D_a$
satisfying the following conditions:

For $C_a\neq C_b$,
we have $D_a^{\circ}\cap D_b=\varnothing$.
The inverse image
${\mathbf P}(\widetilde C_a)^{\circ}
=
{\mathbf P}(\widetilde C_a)
\times_{D_a}
D_a^{\circ}$
is finite and radicial over
$D_a^{\circ}$.
For every $(x,H)
\in 
{\mathbf P}(\widetilde C_a)^{\circ}$
and for every line $L\subset
{\mathbf P}^\vee$
such that $x\in X_L^{\circ}$,
that 
$L$ meets $D_a^{\circ}$
properly at $H=p_L^{\circ}(x)$
and that
the tangent line 
$TL\times_L \{ H\}$ of $L$
at $H=p_L^{\circ}(x)
\in {\mathbf P}^\vee$
is not perpendicular to
the fiber 
$T^*_Q({\mathbf P}\times{\mathbf P}^\vee)
\times_Q(x,H)
\subset
T^*{\mathbf P}^\vee
\times_{{\mathbf P}^\vee}\{H\}$,
we have 
$\varphi_x^{-1}({\cal F},p_L^{\circ})\neq 0$
and
$\varphi_x^q({\cal F},p_L^{\circ} )= 0$
for $q\neq -1$.
\end{lm}

\subsection{Singular support and ramification}\label{ssram}

We assume $k$ is perfect.
In this subsection, we describe the
singular support $SS{\cal F}$
of ${\cal F}=j_!{\cal G}$
for a locally constant sheaf
${\cal G}$ on the complement
$U=X\sm D$ of a divisor $D$
with simple normal crossings
of a smooth scheme $X$
over $k$
and the open immersion $j\colon U\to X$.
First, we study the tamely ramified case.

\begin{pr}\label{prtame}
Let ${\cal G}\neq 0$ be
a locally constant constructible sheaf
of $\Lambda$-modules
on the complement
$U=X\sm D$ of a divisor $D=\bigcup_{i=1}^mD_i$
with simple normal crossings
of a smooth scheme $X$
over a perfect field $k$.
Assume that ${\cal G}$
is {\em tamely ramified} along $D$.
For the open immersion $j\colon U\to X$
and ${\cal F}=j_!{\cal G}$,
we have 
\begin{equation}
SS{\cal F}
=
\bigcup_{I\subset\{1,\ldots,m\}}
T^*_{D_I}X
\label{eqFtame}
\end{equation}
\end{pr}

\proof{
By Lemma \ref{lmlcst}.3,
we have
$SS{\cal F}
\subset
\bigcup_{I\subset\{1,\ldots,m\}}
T^*_{D_I}X$.
By Corollary \ref{corss}.2,
it is reduced to
the case where $X={\mathbf A}^n$
and $U={\mathbf G}_m^n$.
By the induction on $n=\dim X$
and further by the
compatibility with smooth pull-back,
it suffices to show that the fiber
$T^*_0X$ at the origin $0\in X=
{\mathbf A}^n={\rm Spec}\ k[T_1,\ldots,T_n]$
is a subset of $SS^w{\cal F}$.
If $n=0$,
since ${\cal G}\neq 0$,
we have 
$SS{\cal F}=T^*X$ by Lemma \ref{lmmc}.3.

Assume $n>0$.
Let $D_i=(T_i=0)$
$\subset X$ and
$C=\bigcup_{I\subsetneqq\{1,\ldots,n\}}
T^*_{D_I}X$
be the union except the fiber
$T^*_0X$.
Since the morphism
$f\colon X\to Y=
{\rm Spec}\ k[T]$ defined by
$T\mapsto T_1+\cdots+T_n$
is $C$-transversal,
it suffices to show the following.
\qed}

\begin{lm}%[{\cite[Proposition 6]{epsilon}}]
\label{lmtame}
Let $S={\rm Spec}\ {\cal O}_{K}$ be the spectrum
of a henselian discrete valuation 
ring with algebraically closed residue field $k$,
$X$ be a smooth
scheme of finite type 
of relative dimension $n-1$ over $S$
and $D$ be a divisor of $X$
with simple normal crossings.
Let $x$ be a closed point
of the closed fiber of $X$
contained in $D$.
Let $t_1,\ldots,t_n\in {\mathfrak m}_x$
be elements of the maximal ideal such that
$\bar t_1,\ldots,\bar t_n\in 
{\mathfrak m}_x/{\mathfrak m}_x^2$
is a basis.
Assume that $D$
is defined by $t_1\cdots t_n$
and that the class of 
a uniformizer $\pi$ of $S$
in ${\mathfrak m}_x/
{\mathfrak m}_x^2$
is not contained in any subspaces
generated by $n-1$ elements
of the basis $\bar t_1,\ldots,\bar t_n$.

Let $\Lambda$ be a finite local ring
with residue characteristic $\ell$ invertible on $S$
and let ${\cal G}$ be a locally constant
constructible sheaf of $\Lambda$-modules
on the complement $U=X\sm D$
tamely ramified along $D$.
Let $j\colon U\to X$ denote
the open immersion.

{\rm 1.}
On a neighborhood of $x$,
the complex $\phi(j_!{\cal G})$
is acyclic except at $x$
and at degree $n-1$
and the action of the inertia group
$I_K={\rm Gal}(\bar K/K)$
on $\phi^{n-1}_x(j_!{\cal G})$
is tamely ramified.

{\rm 2.}
If ${\cal G}$ is a sheaf of
free $\Lambda$-modules,
then the $\Lambda$-module
$\phi^{n-1}_x(j_!{\cal G})$
is free of rank ${\rm rank}\ {\cal G}$.
\end{lm}

\proof{
1.
We may write $\pi=\sum_{i=1}^n
u_{i}t_{i}$ in ${\mathfrak m}_x$
and $u_{i}$ are invertible.
By replacing $t_{i}$ by 
$u_{i}t_{i}$,
we may assume that 
$X$ is \'etale over
${\rm Spec} \ {\cal O}_{K}
[t_{1},\ldots,t_{n}]/
(\pi-(t_{1}+\cdots+t_n))$.
By Abhyankar's lemma,
we may assume that ${\cal G}$
is trivialized by the abelian covering
$s_{i}^{m}=t_{i}$ for an integer
$m$ invertible on $S$.
Since the assertion is \'etale local on $X$,
we may assume
$X={\rm Spec} \ {\cal O}_{K}
[t_{1},\ldots,t_n]/
(\pi-(t_{1}+\cdots+t_n))$.
Hence by the variant of
\cite[1.3.3 (i)]{app} for $f_!$,
the complex $\phi(j_!{\cal G})$
is acyclic outside $x$.
Since the complex
$\phi(j_!{\cal G})[n-1]$ is a perverse 
sheaf by \cite[Corollaire 4.6]{au},
the complex
$\phi(j_!{\cal G})$ is acyclic except at 
degree $n-1$.

Let $p\colon X'\to X$ be the blow-up at $x$
and $j'\colon U\to X'$ be the open immersion.
Let $D'$ be the proper transform of $D$
and $E$ be the exceptional divisor.
Then, the union of $D'$ with the closed
fiber $X'_s$ has simple normal crossings.
Hence, the action of the inertia group 
$I_K$ on
$\phi(j'_!{\cal G})$ is tamely ramified
by \cite[Proposition 6]{epsilon}
and $\phi_x(j_!{\cal G})=R\Gamma(E,\phi(j'_!{\cal G}))$
is also tamely ramified.

2.
By 1,
we may assume $\Lambda$ is a field.
We consider the stratification of
$E$ defined by the intersections
with the intersections of
irreducible components of $D'$.
Then, on each stratum,
the restriction of the cohomology sheaves
$\phi^q(j'_!{\cal G})$ are locally constant
and are tamely ramified along the boundary
by \cite[Proposition 6]{epsilon}.
Further, the alternating sum of
the rank is $0$ except for
$E^{\circ} =E\sm (E\cap D')$
and equals ${\rm rank}\ {\cal G}$ on
$E^{\circ}$.
Hence, we have
$$\dim \phi_x(j_!{\cal G})
=\chi(E,\phi(j'_!{\cal G}))
=\chi_c(E^{\circ},\phi(j'_!{\cal G}))
={\rm rank}\ {\cal G}\cdot
\chi_c(E^{\circ}).$$
Since $E^{\circ}=E_{n-1}^{\circ}
\subset {\mathbf G}_m^{n-1}
={\rm Spec}\ k[t_1^{\pm 1},\ldots,
t_{n-1}^{\pm 1}]$
is the complement
of the intersection
with the hyperplane
$\sum_it_i+1=0$
and the intersection 
is isomorphic to $E_{n-2}^{\circ}$,
we have
$\chi_c(E^{\circ}_{n-1})=
\chi_c({\mathbf G}_m^{n-1})-
\chi_c(E^{\circ}_{n-2})=
(-1)^{n-1}$
by induction on $n$.
Thus the assertion follows.
\qed}
\medskip

To give a description of the singular
support in some wildly ramified case
using ramification theory,
we briefly recall ramification theory
\cite{AS}, \cite{nonlog}.
Let $K$ be a henselian discrete valuation field
with residue field of characteristic $p>0$
and $G_K={\rm Gal}(K_{\rm sep}/K)$
be the absolute Galois group.
Then, the (non-logarithmic)
filtration $(G_K^r)_{r\geqq 1}$
by ramification groups
is defined in \cite[Definition 3.4]{AS}.
It is a decreasing filtration by closed
normal subgroups
indexed by rational numbers $\geqq 1$.

For a real number $r\geqq 1$,
we define subgroups
$G_K^{r+}\subset G_K^{r-}$
by $G_K^{r+}=\overline{\bigcup_{s>r}
G_K^s}$
and $G_K^{r-}=\bigcap_{s<r}
G_K^s.$
It is proved in \cite[Proposition 3.7 (1)]{AS} that
$G_K^1$ is the inertia group
$I={\rm Ker}(G_K\to G_F)$ where
$G_F$ denotes the absolute Galois
group of the residue field $F$
and $G_K^{1+}$ is the wild inertia group $P$
that is the pro-$p$ Sylow subgroup of $I$.
It is also proved in \cite[Theorem 3.8]{AS} that
$G_K^{r-}=G_K^{r}$
for rational numbers $r>1$
and 
$G_K^{r-}=G_K^{r+}$
for irrational numbers $r>1$.

Let $\Lambda$ be a finite field
of characteristic $\neq p$
and let $V$ be a continuous representation
of $G_K$ on a $\Lambda$-vector space
of finite dimension.
Then, since $P=G_K^{1+}$ is a pro-$p$ group
and since 
$G_K^{r-}=G_K^{r}$
for rational $r$ and
$G_K^{r-}=G_K^{r+}$
for irrational $r$,
there exists a unique decomposition
\begin{equation}
V=\bigoplus_{r\geqq 1}V^{(r)}
\label{eqslp}
\end{equation}
called the slope decomposition
characterized by the condition that
the $G_K^{r+}$ fixed part
$V^{G_K^{r+}}$ is equal to the sum
$\bigoplus_{s\leqq r}V^{(s)}$.

We study a geometric case
where $X$ is a smooth scheme
over a perfect field $k$ of characteristic
$p>0$.
Let $D$ be a reduced and irreducible divisor
and $U=X\sm D$ be the complement.
Let ${\cal G}$ be a 
locally constant constructible sheaf of
$\Lambda$-vector spaces on $U$.
Let $\xi$ be the generic point of
an irreducible component of $D$.
Then, the local ring
${\cal O}_{X,\xi}$ 
is a discrete valuation
ring and the fraction field $K$
of its henselization is called
the local field at $\xi$.
The stalk of ${\cal H}^q{\cal G}$
at the geometric point of $U$
defined by a separable closure
$K_{\rm sep}$ defines
a $\Lambda$-vector space $V^q$
with an action
of the absolute Galois group $G_K$.
%and the total dimension $\dim{\rm tot}_{K_i} V_i^q$ is defined.
%We define the Artin conductor $a_{K_i}({\cal G})\in {\mathbf Z}$ at $\xi_i$ and the Artin divisor $a({\cal G})$ by 
%\begin{align}
%a_{K_i}({\cal G})&=
%\sum_q(-1)^q
%(\dim{\rm tot}_{K_i} V_i^q-
%\dim_{\Lambda}{\cal H}^q{\cal G}_{\bar \xi_i}),
%\label{eqAc}
%\\
%a({\cal G})
%&=\sum_i
%a_{K_i}({\cal G})
%\cdot [D_i]. 
%\label{eqAD}
%\end{align}
%If ${\cal H}^q{\cal G}$ is locally constant for every $q$ on a neighborhood of $\xi_i$, we have $a_{K_i}({\cal G})=0$. If $X$ is a curve, we recover the classical definition of the Artin conductor.

For a rational number $r>1$,
the graded quotient
${\rm Gr}^rG_K
=G_K^r/G_K^{r+}$ is 
a profinite abelian group annihilated by $p$
\cite[Corollary 2.28]{nonlog}
and its dual group is related
to differential forms as follows.
We define ideals
${\mathfrak m}_{K_{\rm sep}}^{(r)}$
and ${\mathfrak m}_{K_{\rm sep}}^{(r+)}$
of the valuation ring ${\cal O}_{K_{\rm sep}}$ by
${\mathfrak m}_{K_{\rm sep}}^{(r)}
=\{x\in K_{\rm sep}\mid {\rm ord}_Kx\geqq r\}$
and ${\mathfrak m}_{K_{\rm sep}}^{(r+)}
=\{x\in K_{\rm sep}\mid {\rm ord}_Kx>r\}$
where ${\rm ord}_K$ denotes the valuation 
normalized by ${\rm ord}_K(\pi)=1$
for a uniformizer $\pi$ of $K$.
The residue field $\bar F$ of ${\cal O}_{K_{\rm sep}}$
is an algebraic closure of $F$ and 
the quotient
${\mathfrak m}_{K_{\rm sep}}^{(r)}
/{\mathfrak m}_{K_{\rm sep}}^{(r+)}$
is an $\bar F$-vector space of dimension 1.
A canonical injection
\begin{equation}
ch\colon
Hom_{{\mathbf F}_p}(
{\rm Gr}^rG_K,
{\mathbf F}_p)
\to 
Hom_{\bar F}
({\mathfrak m}_{\bar K}^{(r)}/
{\mathfrak m}_{\bar K}^{(r+)},
\Omega^1_{X/k,\xi}
\otimes{\bar F})
\label{eqcf}
\end{equation}
is also defined
\cite[Corollary 2.28]{nonlog}.

We say that the ramification of
${\cal G}$ along $D$ is
isoclinic of slope $r\geqq 1$
if $V=V^{(r)}$ in the slope decomposition
(\ref{eqslp}).
The ramification of
${\cal G}$ along $D$ is
isoclinic of slope $1$
if and only if
the corresponding Galois representation $V$
is tamely ramified.
Assume that the ramification of
${\cal G}$ along $D$ is
isoclinic of slope $r>1$.
Assume also that 
$\Lambda$ contains a primitive
$p$-th root of 1 and
identify ${\mathbf F}_p$ with
a subgroup of 
$\Lambda^\times$.
Then, $V=V^{(r)}$ is further decomposed
by characters
\begin{equation}
V=\bigoplus_{\chi\colon
{\rm Gr}^rG_K\to 
{\mathbf F}_p}
\chi^{\oplus m(\chi)}.
\label{eqchd}
\end{equation}
For a character $\chi$ appearing
in the decomposition (\ref{eqchd}),
the twisted differential form
$ch(\chi)$ defined on a
finite covering of a dense open
scheme of $D$ is
called a characteristic form of ${\cal G}$.

Assume that 
$U=X\sm D$ is
the complement of a divisor
with simple normal crossings $D$
and let $D_1,\ldots,D_m$
be the irreducible components of $D$.
We say the ramification of
${\cal G}$ along $D$ is
isoclinic of slope $R=\sum_ir_iD_i$
if the ramification of
${\cal G}$ along $D_i$ is
isoclinic of slope $r_i$
for every irreducible component
$D_i$ of $D$.

In \cite[Definition 3.1]{nonlog},
we define the condition for
ramification of ${\cal G}$ along $D$
to be non-degenerate.
The condition implies that
the characteristic forms
are extended to differential forms
on the boundary without zero.
We say that 
the ramification of ${\cal G}$ is non-degenerate along $D$ 
if it admits \'etale locally
a direct sum decomposition
${\cal G}=\bigoplus_j{\cal G}_j$
such that each ${\cal G}_j$
is isoclinic of slope $R_j\geqq D$
for a ${\mathbf Q}$-linear combination $R_j$
of irreducible components of $D$
and that the ramification of ${\cal G}_j$
is non-degenerate along $D$ at multiplicity $R_j$. Note that 
there exists 
a closed subset of codimension at least 2
such that on its complement,
the ramification of ${\cal G}$ along $D$
is non-degenerate.

We introduce a slightly stronger
condition that implies local acyclicity.
We say that 
the ramification of ${\cal G}$ is 
{\em strongly} non-degenerate along $D$
if it satisfies the condition above
with $R_j\geqq D$
replaced by $R_j=D$
or $R_j>D$.
The inequality $R_j>D$
means that the coefficient in $R_j$
of every irreducible component of $D$
is $>1$.
Note that 
if the ramification of ${\cal G}$ along $D$
is non-degenerate,
on the complement of the singular locus
of $D$, 
the ramification of ${\cal G}$ along $D$
is strongly non-degenerate.

Let $j\colon U=X\sm D\to X$
denote the open immersion.
We define a closed conical subset
$S(j_!{\cal G})\subset T^*X$
following the definition 
of the characteristic cycle given
\cite[Definition 3.5]{nonlog}
in the strongly non-degenerate case.
We will show that $S(j_!{\cal G})$
is in fact equal to the singular support
in Proposition \ref{prnd} below.

Assume first that 
the ramification of ${\cal G}$
along $D$ is isoclinic of
slope $R=D$.
Then, the locally constant sheaf
${\cal G}$ on $U$ is tamely ramified along $D$.
In this case,
let the singular support $SSj_!{\cal G}$
(\ref{eqFtame}) denoted by
\begin{equation}
S(j_!{\cal G})
=
\bigcup_IT^*_{D_I}X
\label{eqtameS}
\end{equation}
where $T^*_{D_I}X$
denotes the conormal bundle
of the intersection $D_I=\bigcap_ID_i$
for a set of indices
$I\subset \{1,\ldots,m\}$.

Assume the ramification of ${\cal G}$
along $D$ is isoclinic of
slope $R=\sum_ir_iD_i>
D=\sum_iD_i$. 
For each irreducible component,
we have a decomposition by 
characters 
$V=\bigoplus_{\chi\colon
{\rm Gr}^{r_i}G_{K_i}\to 
{\mathbf F}_p}
\chi^{\oplus m(\chi)}$.
Further, the characteristic form
of each character $\chi$
appearing in the decomposition
defines a sub line bundle $L_\chi$
of the pull-back $D_\chi\times_XT^*X$
of the cotangent bundle 
to a finite covering $\pi_\chi
\colon D_\chi\to D_i$
by the non-degenerate assumption.
Then, define 
a closed conical subset
$C=S(j_!{\cal G})\subset T^*X$ in
the case $R>D$ by 
\begin{equation}
S(j_!{\cal G})
=
T^*_XX\cup
\bigcup_i
\bigcup_{\chi}
\pi_\chi(L_\chi).
\label{eqwildS}
\end{equation}
In the general strongly non-degenerate case,
we define a closed conical subset
$C=S(j_!{\cal G})\subset T^*X$
by the additivity
and \'etale descent.

\begin{pr}[{\rm \cite[Proposition 3.15]{nonlog}}]\label{prnd}
Assume that the ramification of
a locally constant constructible sheaf
${\cal G}$ of $\Lambda$-modules
on the complement $U=X\sm D$ 
of a divisor with simple normal crossings
is strongly non-degenerate
along $D$.
Then, 
$S(j_!{\cal G})$ defined by 
{\rm (\ref{eqtameS})},
{\rm (\ref{eqwildS})},
the additivity and by \'etale descent
satisfies
\begin{equation}
SS j_!{\cal G}
=
S(j_!{\cal G}).
\label{eqSSj}
\end{equation}
\end{pr}

\proof{
Since the assertion is \'etale local,
we may assume that ${\cal G}$
is isoclinic of slope $R=D$ or $R>D$
and that the ramification of ${\cal G}$
is non-degenerate along $D$ at multiplicity $R$.
If $R=D$, it is proved in Proposition \ref{prtame}.
To prove the case $R>D$,
we use the following Lemma.

\begin{lm}\label{lmDf2}
Let $E$ be a $k$-vector space
of finite dimension and 
$E\to \Gamma(X,{\cal L})$
be a $k$-linear mapping
defining an immersion
$X\to {\mathbf P}={\mathbf P}(E^\vee)$.
Let $T\subset X$ be
an integral closed subscheme
and define a subspace $E'={\rm Ker}(E\to 
\Gamma(T,{\cal L}\otimes_{{\cal O}_X}{\cal O}_T))$
and ${\mathbf P}^{\prime\vee}
={\mathbf P}(E')\subset
{\mathbf P}^{\vee}
={\mathbf P}(E)$.

{\rm 1.}
$T\times {\mathbf P}^{\prime\vee}
\subset
X\times {\mathbf P}^{\vee}$
is a subset of
$T\times_{\mathbf P}Q
\subset X\times_{\mathbf P}Q$
and the complement
$(T\times_{\mathbf P}Q)
\sm
(T\times {\mathbf P}^{\prime\vee})$
is the largest open subscheme where
$T\times_{\mathbf P}Q
\to {\mathbf P}^\vee$
is flat.

{\rm 2.}
The codimension of 
$E'\subset E$
is strictly larger than $\dim T$.
\end{lm}

\proof{1.
For a hyperplane $H\subset {\mathbf P}$,
we have $T\subset H$ 
if $H\in {\mathbf P}^{\prime\vee}$
and
$T\cap H$ is a divisor of $T$
if otherwise.
Hence we have
$T\times {\mathbf P}^{\prime\vee}
\subset
T\times_{\mathbf P}Q$
and on the complement of
$T\times {\mathbf P}^{\prime\vee}$,
the divisor 
$T\times_{\mathbf P}Q$
of $T\times{\mathbf P}^\vee$
is flat over ${\mathbf P}^\vee$.

2.
The immersion $X\to {\mathbf P}$
induces an immersion
$T\to {\mathbf P}(E/E')$.
Hence we have
$\dim T<\dim E/E'$.
\qed}
\medskip

We assume $R>D$.
First, we show the inclusion
$SSj_!{\cal G}
\subset S(j_!{\cal G})$.
Since the assertion is local on $X$,
we may assume that $X$ is quasi-projective.
Let $i\colon X\to {\mathbf P}$
be an immersion 
satisfying the condition {\rm (E)}
before Lemma {\rm \ref{lmloc}}
and the condition {\rm (C)} before
Proposition {\rm \ref{prwi}}.
For irreducible component $D_i$ of $D$,
let ${\mathbf P}^\vee_i={\mathbf P}(E_i)
\subset
{\mathbf P}^\vee={\mathbf P}(E)$
be the subspace defined by
$E_i={\rm Ker}(E\to \Gamma(D_i,{\cal L}))$.
Let $Z\subset X\times_{\mathbf P}Q$
be the union 
$Z=\bigcup_{i=1,\ldots,m}
D_i\times {\mathbf P}^\vee_i$.
Then, by Lemma \ref{lmDf2}.1,
the divisor 
$D\times_{\mathbf P}Q
\subset 
X\times_{\mathbf P}Q$ with
simple normal crossings
is flat over ${\mathbf P}^\vee$
outside $Z$.
Hence, by \cite[Proposition 3.15]{nonlog},
$p^\vee\colon 
X\times_{\mathbf P}Q
\to {\mathbf P}^\vee$
is universally locally acyclic relatively to
the pull-back of $j_!{\cal G}$
outside the union
$Z\cup {\mathbf P}(\widetilde C)$
for $C=S(j!{\cal F})$.
Further by Lemma \ref{lmDf2}.2,
the closed subset
$Z\subset X\times_{\mathbf P}Q$
is of codimension $>\dim X$.
Hence we have
$SSj_!{\cal G}
\subset S(j_!{\cal G})$
by
Corollary \ref{corRn}.

We show the equality
$SSj_!{\cal G}
= S(j_!{\cal G})$
keeping the assumption $R>D$.
By Theorem \ref{thmss}.2,
we may assume $D$ is irreducible
since the assertion is local.
Further we may assume 
there exists a unique $L_\chi$
by the additivity
since the assertion is \'etale local.
Then, the assertion follows from
the contraposition of
Lemma \ref{lmlcst}.1 (1)$\Rightarrow$(2).
\qed}

\section{Characteristic cycle and 
the Milnor formula}\label{sCC}

In this section,
$X$ denotes a smooth scheme
over a perfect field $k$.
Assume that every irreducible component
of $X$ is of dimension $n$,
unless otherwise stated.

\subsection{Morphisms defined by pencils}
\label{sspl}

Let $E$ be a $k$-vector space of
finite dimension and
${\mathbf P}={\mathbf P}(E^\vee)$
and its dual
${\mathbf P}^\vee={\mathbf P}(E)$
be as in Section \ref{ssRn}.
Let $X$ be a smooth scheme
over $k$
and let $i\colon X\to {\mathbf P}$
be an immersion.
Let $L\subset {\mathbf P}^\vee$
be a line.
The morphism $p_L\colon X_L\to L$
defined by pencil
is defined by the cartesian diagram
\begin{equation}
\begin{CD}
X_L@>>>
X\times_{\mathbf P}Q\\
@V{p_L}VV @VV{p^\vee}V\\
L@>>> {\mathbf P}^\vee.
\end{CD}
\label{eqpL}
\end{equation}
The axis $A_L\subset {\mathbf P}$
is the intersection of hyperplanes
parametrized by $L$.
If the axis $A_L$ meets $X$ properly,
the scheme $X_L$ is the blow-up of $X$
at the intersection
$X\cap A_L$.
The morphism
$p_L^{\circ}\colon
X_L^{\circ} \to L$
(\ref{eqpLo})
is the restriction of
$p_L\colon X_L\to L$
to the complement 
$X_L^{\circ} =X\sm (X\cap A_L)$.

\begin{lm}\label{lmAL}
Let $C\subset T^*X$ be a
closed conical subset
and 
${\mathbf P}(\widetilde C)\subset
X\times_{\mathbf P}Q=
{\mathbf P}(X\times_{\mathbf P}T^*{\mathbf P})$
be the projectivization.
Assume that every irreducible component 
of $X$ and every irreducible component
$C_a$ of $C$ are of dimension $n$.
Let $L\subset {\mathbf P}^\vee$
be a line and
$p_L^{\circ}\colon X_L^{\circ}\to L$
be the morphism defined
by the pencil.

The complement
$X_L^{\circ}\sm( X_L^{\circ} \cap
{\mathbf P}(\widetilde C))
\subset X_L^{\circ}$
is the largest open subset where
$p_L^{\circ}\colon X_L^{\circ}\to L$
is $C$-transversal.
\end{lm}

\proof{
Since
$X\times_{\mathbf P}Q\to X$
is smooth,
the immersion $X_L^{\circ}\to
X\times_{\mathbf P}Q$ is $C$-transversal
by Lemma \ref{lmCh}.3.
Hence, for $x\in X_L^{\circ}$, the morphism
$p_L^{\circ}\colon X_L^{\circ}\to L$
is $C$-transversal at $x$ if and only if
$(x,p_L(x))\in X\times_{\mathbf P}Q$
is not contained in ${\mathbf P}(\widetilde C)$
by Lemma \ref{lmPC}.1 and by
Lemma \ref{lmChf}.1
applied to the cartesian diagram (\ref{eqpL}).
Hence the assertion follows.
\qed}
\medskip

We construct the universal family of
(\ref{eqpL}).
Let ${\mathbf G}=
{\rm Gr}(1,{\mathbf P}^\vee)$ be
the Grassmannian variety
parametrizing lines in
${\mathbf P}^\vee$.
The universal line
${\mathbf D}\subset
{\mathbf P}^\vee\times
{\mathbf G}$ 
is canonically identified with the flag variety
parametrizing
pairs $(H,L)$ of points $H$ of
${\mathbf P}^\vee$ 
and lines $L$ passing through $H$.
It is the same as the flag variety
${\rm Fl}(1,2,E)$ 
parametrizing pairs
of a line and a plane including the line in $E$.

The projective space
${\mathbf P}^\vee$ and the
Grassmannian variety
${\mathbf G}$
are also equal to
the Grassmannian varieties
${\rm Gr}(1,E)$ and
${\rm Gr}(2,E)$ 
parametrizing lines and planes
in $E$ respectively.
The projections
${\mathbf P}^\vee
\gets {\mathbf D}
\to {\mathbf G}$
sending a pair $(H,L)$ to the line $L$
and to the hyperplane $H$
are the canonical morphisms
${\rm Gr}(1,E)\gets {\rm Fl}(1,2,E)\to {\rm Gr}(2,E)$.
By the projection ${\mathbf D}\to
{\mathbf P}^\vee$,
it is also canonically identified
with the projective space bundle
associated to the tangent bundle
${\mathbf D}={\mathbf P}(T{\mathbf P}^\vee)$
by sending a line passing through
a point to the tangent line
of the line at the point.
Let ${\mathbf A}\subset {\mathbf P}
\times {\mathbf G}$
be the universal family
of the intersections
of hyperplanes parametrized
by lines.
The scheme  ${\mathbf A}$
is also canonically identified
with the Grassmann bundle
${\rm Gr}(2,T^*{\mathbf P})$
over ${\mathbf P}$.

Let $X$ be a smooth scheme
over $k$
and let $i\colon X\to {\mathbf P}$
be an immersion.
We construct a commutative diagram
\begin{equation}
\begin{CD}
X\times_{\mathbf P}Q
@<<<
(X\times {\mathbf G})'
@>>> 
X\times {\mathbf G}
\\
@V{p^\vee}VV @VVV @VVV\\
 {\mathbf P}^\vee
 @<<<
{\mathbf D}
@>>>
{\mathbf G}.
\end{CD}
\label{eqGLP}
\end{equation}
where the left square is cartesian
as follows.
The left vertical arrow and the lower line
are the canonical morphisms
${\rm Gr}(1,X\times_{\mathbf P}T^*{\mathbf P})
\to
{\rm Gr}(1,E)
\gets
{\rm Fl}(1,2,E)
\to
{\rm Gr}(2,E)$.
The fiber product
${\mathbf D}\times_{{\mathbf P}^{\vee}} 
(X\times_{\mathbf P}Q)$
is canonically identified with 
the blow-up 
$(X\times {\mathbf G})'
\to X\times {\mathbf G}$
at the intersection ${\mathbf A}
\cap (X\times {\mathbf G})
=
X\times_{\mathbf P} {\mathbf A}
=
{\rm Gr}(2,X\times_{\mathbf P} {\mathbf A})$
by the elementary lemma below.

The canonical morphism 
$(X\times {\mathbf G})'
\to X\times {\mathbf G}$
is an isomorphism on 
the complement
\begin{equation}
(X\times {\mathbf G})^{\circ}
=
(X\times {\mathbf G})
\sm
(X\times_{\mathbf P} {\mathbf A}).
\label{eqXGo}
\end{equation}
The scheme
$(X\times {\mathbf G})^{\circ}$
regarded as an open subscheme of
$(X\times {\mathbf G})'$
is the complement of
${\rm Fl}(1,2,X\times_{\mathbf P}T^*{\mathbf P})$
regarded as a closed subscheme.
For a line
$L\subset {\mathbf P}^\vee$,
the diagram (\ref{eqpL}) is the fiber of 
the middle vertical arrow 
$(X\times {\mathbf G})'
\to {\mathbf D}$ of (\ref{eqGLP}) at the
point of ${\mathbf G}$ corresponding to $L$.
The complement
$X_L^{\circ}$ 
is identified with 
$X_L\cap (X\times {\mathbf G})^{\circ}$.

\begin{lm}\label{lmFEL}
Let $0\to F\to E\to L\to 0$
be an exact sequence of
vector bundles on a scheme $X$
such that $L$ is a line bundle
and let $1\leqq r\leqq {\rm rank}\ F-1$ 
be an integer.
Then, the Grassmannian bundles
parametrizing subbundles
of rank $r$ and of $r-1$
and the flag bundles
parametrizing pairs of subbundles
of rank $r$ and $r-1$
with inclusion form a cartesian diagram
\begin{equation}
\begin{CD}
{\rm Fl}(r-1,r,F)
@>>>
 {\rm Gr}(r,F)\\
@VVV@VVV\\
P={\rm Gr}(r-1,F)\times_
{{\rm Gr}(r-1,E)}
{\rm Fl}(r-1,r,E)
@>>>
 {\rm Gr}(r,E).
\end{CD}
\label{eqFlG}
\end{equation} 
The right vertical arrow
is a regular immersion
of codimension $r$
and the lower horizontal arrow
is the blow-up at the
closed subscheme
${\rm Gr}(r,F)\subset
{\rm Gr}(r,E)$.
\end{lm}

\proof{
Let $p\colon P={\rm Gr}(r-1,F)\times_
{{\rm Gr}(r-1,E)}
{\rm Fl}(r-1,r,E)\to X$
and
$$\xymatrix{
D={\rm Fl}(r-1,r,E)\ar[rd]_f
\ar[rr]^d
&&
G={\rm Gr}(r,E)\ar[ld]^g
\\&X
}$$
denote the canonical morphisms. 
Let $0\to {\cal E}\to
{\cal F}\to{\cal L}\to 0$
be the corresponding exact sequence
of locally free ${\cal O}_X$-modules.
Let ${\cal V}\subset g^*{\cal E}$
and ${\cal W}\subset d^*{\cal V}
\subset f^*{\cal E}$
denote the universal sub  
${\cal O}_G$-module of rank $r$
and the universal sub  
${\cal O}_D$-module of rank $r-1$
respectively.

On $P\times_{
{\rm Gr}(r,E)}
{\rm Gr}(r,F)$,
the pull-back of
${\cal W}\subset d^*{\cal V}$
defines a flag on the pull-back of
${\cal F}$.
Hence, the diagram (\ref{eqFlG})
is cartesian.
On the complement
${\rm Gr}(r,E)\sm
{\rm Gr}(r,F)$,
the restriction
of ${\cal V}$ and its
intersection with
the pull-back of
${\cal F}$ 
define a flag on the restriction of
$g^*{\cal E}$.
Hence the restriction
$P\sm {\rm Fl}(r-1,r,F)
\to {\rm Gr}(r,E)\sm
{\rm Gr}(r,F)$
is an isomorphism.

The ideal sheaf
${\cal I}\subset {\cal O}_G$ 
defining the closed subscheme ${\rm Gr}(r,F)
\subset{\rm Gr}(r,E)$
is characterized by the condition that
${\cal I}\cdot g^*{\cal L}$
equals the image of the
composition ${\cal V}\to g^*{\cal E}
\to g^*{\cal L}$.
Hence, the immersion
${\rm Gr}(r,F)
\to{\rm Gr}(r,E)$ is a regular
immersion of codimension $r={\rm rank}\ {\cal V}$.

Since the image of the
composition ${\rm pr}_2^*d^*{\cal V}
\to p^*{\cal E}
\to p^*{\cal L}$ is isomorphic to
the invertible ${\cal O}_P$-module
${\rm pr}_2^*(d^*{\cal V}/{\cal W})$,
the ideal 
${\cal I}\cdot {\cal O}_P$ 
defining the closed subscheme
${\rm Fl}(r-1,r,F)\subset P$
is locally generated by one element.
Since ${\rm Fl}(r-1,r,F)$ and $P$
are smooth over $X$,
${\rm Fl}(r-1,r,F)$ is a divisor of $P$
flat over $X$ and 
the ideal 
${\cal I}\cdot {\cal O}_P$ is
invertible.
Thus the morphism $P\to {\rm Gr}(r,E)$
is canonically lifted to the blow-up
$P\to {\rm Gr}(r,E)'$
by the universality of blow-up.

Since the pull-back $\pi^*{\cal I}$
on the blow-up $\pi\colon {\rm Gr}(r,E)'
\to {\rm Gr}(r,E)$ is an invertible
ideal, the image $\pi^*{\cal I}\cdot
\pi^*g^*{\cal L}$ of
the composition
$\pi^*{\cal V}\to
\pi^*g^*{\cal E}\to
\pi^*g^*{\cal L}$ is an invertible module.
Hence the kernel
$\pi^*{\cal V}\cap
\pi^*g^*{\cal F}$ is locally 
a direct summand of
$\pi^*g^*{\cal F}$ of rank $r-1$
and 
$\pi^*{\cal V}\cap
\pi^*g^*{\cal F}
\subset \pi^*{\cal V}$ 
defines a flag on
$\pi^*g^*{\cal E}$.
Thus 
the inverse morphism
${\rm Gr}(r,E)'\to P$
is defined
and the assertion follows.
\qed}
%\medskip

\subsection{Isolated characteristic points and
intersection numbers}
\label{ssic}

The characteristic cycle will
be defined as a cycle
characterized by the Milnor formula
at isolated characteristic points.
We will introduce the notion of
isolated characteristic point
and study the intersection number
appearing in the Milnor formula.

\begin{df}\label{dfisoc}
Let $X$ be a smooth scheme over a
field $k$.
Let $C\subset T^*X$ be a 
closed conical subset
of the cotangent bundle $T^*X$.
Let $h\colon W\to X$
be an \'etale morphism
and $f\colon W\to Y$ be a morphism
over $k$ to a smooth curve over $k$.

{\rm 1.}
We say that a closed point $u$ of $W$
is at most an {\em isolated 
$C$-characteristic point} of 
$f\colon W\to Y$% with respect to $C$
if there exists an open neighborhood 
$V\subset W$ of $u$ such that the pair 
$X\gets V\sm \{u\}\to Y$ is $C$-transversal.
We say that a closed point $u\in W$
is an {\em isolated $C$-characteristic point} of $f$
if it is at most an isolated $C$-characteristic
point but 
$X\gets W\to Y$ is not $C$-transversal
at $u$.

{\rm 2.}
Assume that every irreducible component 
of $X$ and every irreducible component
$C_a$ of $C$ are of dimension $n$.
Let $u$ be at most an isolated $C$-characteristic point of 
$f\colon W\to Y$ with respect to $C$
and $A=\sum_am_a[C_a]$ be
a linear combination of irreducible components of $C$.
Then, we define the intersection number
\begin{equation}
(A,df)_{T^*W,u}
\label{eqAdf}
\end{equation}
as the intersection number
$\sum_am_a(j^*C_a,f^*\omega)_{T^*W,u}$
supported on the fiber of $u$
for the section $f^*\omega$ of $T^*W$
defined by the pull-back of
a basis $\omega$ of $T^*Y$
on a neighborhood of $f(u)\in Y$.
\end{df}

The cotangent bundle $T^*W$
is of dimension $2n$ and
its closed subsets
$j^*C_a,f^*\omega$
are of dimension $n$.
Their intersections
$j^*C_a \cap f^*\omega$
consist of at most a unique isolated
point $f^*\omega(u)
\in T^*_uW$
on the fiber of $u$
and the intersection numbers
$(j^*C_a,f^*\omega)_{T^*W,u}$
are defined
if $u$ is at most an isolated 
$C$-characteristic point.
Further, since $C$ is conical,
the  intersection numbers
$(j^*C_a,f^*\omega)_{T^*W,u}$
are independent of the choice of
basis $\omega$
and the intersection number
$(A,df)_{T^*W,u}$
is well-defined.

We compute the intersection number
(\ref{eqAdf}) for
morphisms defined by pencils.

\begin{lm}\label{lmApL}
Let $C\subset T^*X$ be a
closed conical subset
and 
${\mathbf P}(\widetilde C)\subset
X\times_{\mathbf P}Q=
{\mathbf P}(X\times_{\mathbf P}T^*{\mathbf P})$
be the projectivization.
Assume that every irreducible component 
of $X$ and every irreducible component
$C_a$ of $C$ are of dimension $n$.
Let $L\subset {\mathbf P}^\vee$
be a line and
$p_L^{\circ}\colon X_L^{\circ}\to L$
be the morphism defined
by the pencil.

{\rm 1.}
Let $u\in X_L^{\circ}$ be a closed point.
Then, $u$ is (resp.\ at most) an isolated
characteristic point of 
$p_L^{\circ}\colon X_L^{\circ}\to L$
if and only if (resp.\ either)
$u$ is an isolated point of
(resp.\ or not contained in)
the intersection
$X_L^{\circ} \cap
{\mathbf P}(\widetilde C)\subset
X\times_{\mathbf P}Q$.

Further, if $u$ is at most an isolated
characteristic point of 
$p_L^{\circ}\colon X_L^{\circ}\to L$,
we have an equality
\begin{equation}
(C,dp_L^{\circ})_{T^*X,u}
=
({\mathbf P}(\widetilde C),X_L^{\circ})_{
X\times_{\mathbf P}Q,u}
\label{eqAL}
\end{equation}
of the intersection numbers.

{\rm 2.}
Assume that $k$ is algebraically closed
and 
that $C$ is irreducible.
Suppose that $p^\vee\colon 
X\times_{\mathbf P}Q
\to {\mathbf P}^\vee$ is generically 
finite on ${\mathbf P}(\widetilde C)$
and let $\xi\in {\mathbf P}(\widetilde C)$
and $\eta\in \Delta
=\overline{p^\vee({\mathbf P}(\widetilde C))}
\subset {\mathbf P}^\vee$
denote the generic points.
Then, there exists a smooth dense
open subscheme $\Delta^{\circ}
\subset \Delta$ satisfying
the following condition:

For a line $L\subset {\mathbf P}^\vee$
meeting $\Delta$ transversally
at $H\in \Delta^{\circ}$ and 
for an isolated $C$-characteristic point $u$
of $p_L^{\circ}\colon 
X_L^{\circ}\to L$ such that $(u,H)\in 
{\mathbf P}(\widetilde C)^{\circ}
={\mathbf P}(\widetilde C)\times
_\Delta\Delta^{\circ}$, 
the both sides of {\rm (\ref{eqAL})}
are equal to the inseparable
degree $[\xi:\eta]_{\rm insep}$.
\end{lm}

Since ${\mathbf P}(\widetilde C)
\subset
X\times_{\mathbf P}Q$ 
is of codimension $n$,
the intersection product in the right hand
side of (\ref{eqAL}) is defined.

\proof{
1.
The first assertion follows
immediately from
Lemma \ref{lmAL}.
We show the equality (\ref{eqAL}).
Let $\tilde p_L^{\circ}\colon 
{\mathbf P}_L^{\circ}=
{\mathbf P}\sm {\mathbf A}_L
\to L$ denote the morphism
$p_L^{\circ}$ defined with
$X$ replaced by ${\mathbf P}$.
Let $\omega$ be a basis of
$T^*L$ on the image $p_L(u)$
and let $\tilde p_L^{\circ*}\omega$
denote the section
$X_L^{\circ} \to 
X_L^{\circ} \times_{\mathbf P}T^*{\mathbf P}$
defined on a neighborhood of $u$ by
the pull-back of $\omega$
by $\tilde p_L^{\circ}$.
Then, we have
\begin{equation}
(C,dp_L^{\circ})_{T^*X,x}
=
(C,p_L^{\circ *}\omega)_{T^*X,x}
=
(\widetilde C,\tilde
p_L^{\circ *}\omega)_
{X\times_{\mathbf P}T^*{\mathbf P},x}
=
({\mathbf P}(\widetilde C),
\overline{\widetilde p_L^{\circ *}\omega})_{
X\times_{\mathbf P}Q,x}.
\label{eqCpL}
\end{equation}

The graph ${\mathbf P}_L^{\circ} \to 
{\mathbf P}\times L$
of $\tilde p_L^{\circ}$
is a regular immersion of codimension 1
and defines an exact sequence
\begin{equation}
\begin{CD}
0\to
{\mathbf P}_L^{\circ}\times_Q
T^*_Q
({\mathbf P}\times {\mathbf P}^\vee)
@>>> 
({\mathbf P}_L^{\circ}\times_{\mathbf P}T^*{\mathbf P})
\times_{{\mathbf P}_L^{\circ}}
({\mathbf P}_L^{\circ}\times_LT^*L)
@>>>T^*{\mathbf P}_L^{\circ} \to0\end{CD}
\label{eqTPL}
\end{equation}
of vector bundles on ${\mathbf P}_L^{\circ}$.
Hence the left arrow induces
an isomorphism
${\mathbf P}_L^{\circ}\times_Q
T^*_Q
({\mathbf P}\times {\mathbf P}^\vee)
\to
{\mathbf P}_L^{\circ}\times_LT^*L$
to the second factor.
Since 
$T^*_Q
({\mathbf P}\times {\mathbf P}^\vee)
\to 
Q\times_{\mathbf P}T^*{\mathbf P}$
is the universal sub line bundle
on
$Q={\mathbf P}(T^*{\mathbf P})$,
the morphism
${\mathbf P}_L^{\circ}\times_LT^*L
\to
T^*{\mathbf P}_L^{\circ}$
also defines the restriction
of the universal sub line bundle
on ${\mathbf P}_L^{\circ}\subset Q$.
Hence
the image of the section
$\overline{\widetilde p_L^{\circ *}\omega}
\colon X_L^{\circ} \to 
X \times_{\mathbf P}Q=
{\mathbf P}(X\times_{\mathbf P}T^*{\mathbf P})$ equals 
$X_L^{\circ} \subset
X\times_{\mathbf P}Q$
and (\ref{eqCpL}) implies (\ref{eqAL}).

2.
Let $\Delta^{\circ}\subset \Delta$
be a smooth dense open subscheme 
such that the base change
${\mathbf P}(\widetilde C)^{\circ}
=
{\mathbf P}(\widetilde C)
\times_\Delta\Delta^{\circ}
\to \Delta^{\circ}$
is the composition of
a finite flat radicial morphism
of degree $[\xi:\eta]_{\rm insep}$
and a finite \'etale morphism.
Then, for a line $L$ and $(u,H)\in 
{\mathbf P}(\widetilde C)^{\circ}$ as in the assumption,
the intersection number
$({\mathbf P}(\widetilde C),X_L^{\circ})_{
X\times_{\mathbf P}Q,u}$
equals the degree of
the localization at $u$
of the fiber of the finite flat morphism
${\mathbf P}(\widetilde C)^{\circ}\to
\Delta^{\circ}$
at $H\in \Delta^{\circ}$
and is equal to
$[\xi:\eta]_{\rm insep}.$
\qed}
\medskip

We give a condition 
for a function on isolated characteristic
points to be given 
as an intersection number.

\begin{df}\label{dfflZ}
Let $k$ be an algebraically closed field.
%of characteristic $p\geqq 0$.
Let $f\colon Z\to S$ be a {\em quasi-finite} morphism
of schemes of finite type 
over $k$.
We say that a function $\varphi
\colon Z(k)\to {\mathbf Q}$
on the set of closed points
is {\em flat over} $S$
if for every closed point $x\in Z$,
there exists a commutative diagram
$$\xymatrix
{U\ar[r]\ar[dr]_g&
V\times_SZ\ar[r]^h
\ar[d]&
Z\ar[d]^f\\
& V\ar[r]&S
}$$
%with cartesian square
satisfying the following conditions
{\rm (1)--(4):}

{\rm (1)} $V\to S$ is an \'etale neighborhood
of  $s=f(x)$.

{\rm (2)} $U\subset V\times_SZ$ is an open neighborhood
of $x$.

{\rm (3)} $g\colon
U\to V$ is finite
and $g^{-1}(s)=\{x\}$.

{\rm (4)} The function $g_*\varphi$ on $V(k)$
defined by $g_*\varphi(t)=
\sum_{z\in g^{-1}(t)}\varphi(hz)$
is {\em constant}.
\end{df}

\begin{lm}\label{lmfla}
Let $f\colon Z\to S$ be a {\em quasi-finite} morphism
of schemes of finite type over
an algebraically closed field $k$.

{\rm 1.}
Let $\varphi
\colon Z\to {\mathbf Q}$ be
a function.
If $\varphi$ is flat over $S$
in the sense of Definition {\rm \ref{df11}},
its restriction on $Z(k)$ is flat over $S$
in the sense of Definition {\rm \ref{dfflZ}}.

{\rm 2.}
Let $\varphi
\colon Z(k)\to {\mathbf Q}$
be a function flat over $S$
in the sense of Definition {\rm \ref{dfflZ}}.
Then, it is constructible on $Z(k)$.
If $\varphi=0$ on $U(k)\subset
Z(k)$ for a dense open subset
$U\subset Z$,
we have
$\varphi=0$ on $Z(k)$.
\end{lm}

\proof{1.
Let the notation be as in 
Definition \ref{dfflZ}.
Then, by Lemma \ref{lmscZ}.4,
$g_*\varphi$ is locally constant
and the assertion follows.

2.
If $Z\to S$ is finite, surjective and
radiciel, then $\varphi$ is locally
constant.
The constructibility follows from this
by devissage.

Since the second assertion is
\'etale local on $Z$,
we may assume $Z$ is finite over $S$.
Then $\varphi$ on $Z(k)$
is determined by its restriction
to a dense open subset by
the condition (4) in Definition \ref{dfflZ}.
\qed}

\begin{df}\label{dffla}
Let $X$ be a smooth scheme
over $k$ and
$C\subset T^*X$
be a closed conical subset.

{\rm 1.}
We say that $\varphi$ is
a {\em function on 
isolated $C$-characteristic points}
if a number $\varphi(f,u)\in {\mathbf Q}$
%Z}[\frac 1p]$ if $p>0$
%(resp.\ $\varphi(f,u)\in {\mathbf Z}$ if $p=0$)
is defined
%Assume that every irreducible component of $X$
%and of $C$ are of dimension $n$.
%Assume that 
for every 
morphism $f\colon U\to Y$
over $k$
defined on an open subscheme $U\subset X$
to a smooth curve $Y$
with at most 
an isolated $C$-characteristic point $u\in U$
%Z}[\frac 1p]$ if $p>0$
%(resp.\ $\varphi(f,u)\in {\mathbf Z}$ if $p=0$)
and if we have $\varphi(f,u)=0$
if $u\in U$ is not
an isolated $C$-characteristic point.

{\rm 2.}
Let $\varphi$ be
a function on 
isolated $C$-characteristic points.
We say that $\varphi$ is {\em flat}
if for every commutative diagram
\begin{equation}
\xymatrix{
Z\ar[r]^{\subset} &
U\ar[rr]^f\ar[rd]
\ar[d]_{{\rm pr}_1}&
&Y\ar[ld]^g\\
&X&S&}
\label{eqfla}
\end{equation}
of schemes over $k$
satisfying the conditions 
{\rm (1)--(5)} below,
the function 
$\varphi_f$ on $Z(k)$
defined by $\varphi_f(u)=\varphi(f_s,u)$ 
for the base change
$f_s\colon U_s
\to Y_s$ of $f$ at $s={\rm pr}_2(u)\in S$
is {\rm flat} in the sense of 
Definition {\rm \ref{dfflZ}:}

{\rm (1)} $S$ is smooth over $k$.

{\rm (2)} $g\colon Y\to S$
is a smooth curve.

{\rm (3)} $U$ is an open subscheme
of $X\times S$.

{\rm (4)} $Z\subset U$
is a closed subset {\em quasi-finite}
over $S$.

{\rm (5)} The pair
$({\rm pr}_1,f)$
is $C$-transversal
on the complement $U\sm Z$.
\end{df}

Note that 
the morphisms
$f_s\colon U_s\to Y_s$
are $C$-transversal except
possibly at at most isolated
characteristic points at $Z_s$
under the conditions
(1)--(5) on the diagram (\ref{eqfla})
by Lemma \ref{lmChf}.1.

\begin{pr}\label{prfla}
Let $X$ be a smooth scheme
over
an algebraically closed field
$k$ of characteristic $p\geqq 0$
and
$C\subset T^*X$
be a closed conical subset.
Assume that every irreducible
component of $X$
and of $C=\bigcup_aC_a$ are of dimension $n$.
Let $\varphi$ be
a ${\mathbf Z}[\frac 1p]$-valued
(resp.\ ${\mathbf Z}$-valued)
function on isolated $C$-characteristic points
if $p>0$
(resp.\ if $p=0$).

{\rm 1.}
The following conditions
are equivalent:

{\rm (1)}
$\varphi$ is flat
in the sense of Definition {\rm \ref{dffla}.2}.
If $g\colon Y\to Z$ is an \'etale morphism
of smooth curves over $k$,
we have
$\varphi(f,u)=\varphi(gf,u)$.

{\rm (2)}
There exists a
${\mathbf Z}[\frac 1p]$-linear
(resp.\ ${\mathbf Z}$-linear)
combination
$A=\sum_am_aC_a$ satisfying
\begin{equation}
\varphi(f,u)=(A,df)_{T^*X,u}
\label{eqafu}
\end{equation}
for every  $f\colon U\to Y$
with at most
an isolated $C$-characteristic point $u\in U$.

Further $A$ in {\rm (2)} is unique.
Further $A$ is independent of $C$
in the sense that if $C'\supset C$
is a closed conical subset
such that
every irreducible
component of $C'=\bigcup_bC'_b$ is of dimension $n$,
then the linear combination
$A'=\sum_bm'_bC'_b$ satisfying
{\rm (\ref{eqafu})}
for every  $f\colon U\to Y$
with at most
an isolated $C'$-characteristic point $u\in U$
equals $A$.

{\rm 2.}
Let $h\colon W\to X$ be an \'etale
morphism of smooth schemes
over $k$ and
let $\psi$ be
a ${\mathbf Z}[\frac 1p]$-valued
(resp.\ ${\mathbf Z}$-valued)
function on isolated $h^* C$-characteristic points
if $p>0$
(resp.\ if $p=0$).
We assume that 
$\varphi$ and $\psi$ satisfy
the equivalent condition in {\rm 1.}\
and let $A$ and $A'$
be the linear combinations
of irreducible components
of $C$ and of $h^* C$
satisfying {\rm (\ref{eqafu})}
for $\varphi$ and $\psi$
respectively.
If $\varphi(f,hv)=\psi(fh,v)$
for isolated characteristic points
$v$ of $fh\colon W\times_XU\to U\to Y$,
we have $A'=h^*A$.
\end{pr}

\proof{
1.
(2)$\Rightarrow$(1):
Let the notation be as in
(\ref{eqfla}).
We show that
the function
$\varphi_f$ on $Z(k)$
defined by
$\varphi_f(u)=(A,df_s)_{T^*X,u}$
is flat over $S$ in the sense of Definition
\ref{dffla}.2.
We may assume  $A=C_a=C$.
Since the assertion is local,
we may assume that
$\Omega^1_{Y/S}$
is free of rank $1$
and the section
$df\colon U\to T^*X$
defined by a basis
is globally defined on $U$.
Since $T^*X$ is regular,
the ${\cal O}_{T^*X}$-module
${\cal O}_C$ is of finite tor-dimension.
Since $U\to S$ is flat,
the complex of 
${\cal O}_U$-modules 
${\cal O}_C
\otimes^L_{{\cal O}_{T^*X}}
{\cal O}_U$ defined
as the pull-back by $df$
is of finite tor-dimension
as a complex of 
${\cal O}_S$-modules.
Hence the function
$\varphi_f$ on $Z(k)$
is flat over $S$ by Lemma \ref{lmA}.1
and Lemma \ref{lmfla}.1.

If $g\colon Y\to Z$ is an \'etale morphism
of smooth curves over $k$,
we have
$(A,df)_{T^*X,u}
=(A,dgf)_{T^*X,u}$.

(1)$\Rightarrow$(2):
Since the question is local,
we may assume $X$ is affine.
We take a closed immersion
$X\to {\mathbf A}^n
={\rm Spec}\ k[T_1,\ldots,T_n]
\subset
{\mathbf P}^n$.
Let $E=\Gamma({\mathbf P}^n,
{\cal O}(1))$ and let
${\cal L}$ be the pull-back to
$X$ of ${\cal O}(1)$.
After replacing $E\to 
\Gamma(X,{\cal L})$
by $S^dE\to 
\Gamma(X,{\cal L}^{\otimes d})$
for $d\geqq 3$,
we may assume that
the conditions (E) and (C)
before and after Lemma \ref{lmloc}
are satisfied.
We may identify
$E$ with the $k$-linear subspace
of $k[T_1,\ldots,T_n]$
consisting of polynomials of degree $\leqq 1$.
Similarly,
$S^dE$ is identified
with the $k$-linear subspace
of $k[T_1,\ldots,T_n]$
consisting of polynomials of degree $\leqq d$.

We consider the universal family 
as in Section \ref{sspl}.
Define an open subset
\begin{equation}
(X\times {\mathbf G})^\triangledown
\subset
(X\times {\mathbf G})^{\circ}
\label{eqXGt}
\end{equation}
of $(X\times {\mathbf G})^{\circ}
\subset (X\times {\mathbf G})'$
(\ref{eqXGo}) 
to be the largest open subset such that
${\mathbf Z}(\widetilde C)$
defined by the left cartesian square in
\begin{equation}
\begin{CD}
{\mathbf Z}(\widetilde C)
@>>>
(X\times {\mathbf G})^\triangledown
@>f>> {\mathbf D}
@>>> {\mathbf G}\\
@VVV@VVV@VVV\\
{\mathbf P}(\widetilde C)
@>>> X\times_{\mathbf P}Q
@>>> {\mathbf P}^{\vee}
\end{CD}
\label{eqZPC}
\end{equation}
is quasi-finite over
${\mathbf G}$.
By Lemma \ref{lmPC}.1,
the complement
$(X\times {\mathbf G})^\triangledown
\sm
{\mathbf Z}(\widetilde C)$
is the largest open subset
where the pair of
${\rm pr}_1\colon
(X\times {\mathbf G})^\triangledown
\to X$
and 
$f\colon
(X\times {\mathbf G})^\triangledown
\to {\mathbf D}$
is $C$-transversal.
Further by Lemma \ref{lmAL},
for the pair $(u,L)\in (X\times {\mathbf G})^{\circ}$
of a line $L\subset {\mathbf P}^\vee$
and $u\in X_L^{\circ}$,
the pair $(u,L)\in (X\times {\mathbf G})^{\circ}$ is a
point of ${\mathbf Z}(\widetilde C)$
(resp.\ of $(X\times {\mathbf G})^\triangledown$)
if and only if
$u\in X_L^{\circ}$
is (resp.\ at most)
an isolated $C$-characteristic point
of $p_L^{\circ} \colon
X_L^{\circ} \to L$. %with respect to $C$.

We consider the diagram
\begin{equation}
\xymatrix{{\mathbf Z}(\widetilde C)
\ar[r]^{\!\!\!\!\!\!\!\!\subset}&
(X\times {\mathbf G})^\triangledown
\ar[rr]^f\ar[rd]_{p^\triangledown}
\ar[d]&&
{\mathbf D}\ar[ld]^g\\
&X&{\mathbf G}}
\label{eqMfg}
\end{equation}
as (\ref{eqfla}).
For a point $L$ of ${\mathbf G}$,
the fiber of $f\colon 
(X\times {\mathbf G})^\triangledown
\to {\mathbf D}$
is a restriction of $p_L^{\circ}\colon 
X_L^{\circ} \to L$.
By the condition (1),
the function $\varphi_f$ on 
${\mathbf Z}(\widetilde C)$
is flat.
Hence, there exists a dense
open subscheme ${\mathbf Z}(\widetilde C)^{\circ}
\subset {\mathbf Z}(\widetilde C)$
where the function
$\varphi_f$
is locally constant by Lemma \ref{lmfla}.2.

For each irreducible component
$C_a$ of $C=\bigcup_aC_a$,
the function $\varphi_f$
is constant on a dense open subscheme
${\mathbf Z}(\widetilde C_a)^{\circ}
=
{\mathbf Z}(\widetilde C_a)
\cap {\mathbf Z}(\widetilde C)^{\circ}
\subset
{\mathbf Z}(\widetilde C_a)$.
Define a number $\varphi_a
\in {\mathbf Z}[\frac1p]$ to be
the value of $\varphi_f$
on ${\mathbf Z}(\widetilde C_a)^{\circ}$.
The restriction 
${\mathbf P}(\widetilde C_a)
\to D_a=
\overline {p^\vee({\mathbf P}(\widetilde C_a))}$
of
$p^\vee\colon X\times_{\mathbf P}Q
\to {\mathbf P}^\vee$
is generically finite
by Corollary \ref{corrad}.1
since $E$ is assumed to satisfy the conditions 
(E) and (C) before 
and after Lemma {\rm \ref{lmloc}}.
Let $\xi_a\in 
{\mathbf P}(\widetilde C_a)$
and $\eta_a\in D_a$
be the generic points.
We define
\begin{equation}
A=
\sum_a
\dfrac{\varphi_a}
{[\xi_a:\eta_a]_{\rm insep}}[C_a].
\label{eqChE}
\end{equation}
Since the inseparable degree
$[\xi_a:\eta_a]_{\rm insep}$ is a power of $p$
if $p>0$ (resp. is $1$ if $p=0$),
the coefficients
in $A$ are in ${\mathbf Z}[\frac1p]$
(resp. in ${\mathbf Z}$).

We show that $A$
satisfies (\ref{eqafu})
for morphisms defined by pencils.
Let $L\subset {\mathbf P}^\vee$ be a line
and
$u\in X_L^{\circ}$ be at most an isolated
characteristic point of
$p_L^{\circ}\colon 
X_L^{\circ}\to L$.
If $u\in X_L^{\circ}$ is not an isolated
characteristic point of
$p_L^{\circ}$,
then the both sides of (\ref{eqafu}) are $0$.
Assume $u\in X_L^{\circ}$ is an isolated
characteristic point of
$p_L^{\circ}$.
Then we have 
$(u,L)\in {\mathbf Z}(\widetilde C)$
by Lemma \ref{lmApL}.1.

Shrinking ${\mathbf Z}(\widetilde C)^{\circ}$
if necessary,
we may assume that
${\mathbf Z}(\widetilde C)^{\circ}
=
\coprod_a
{\mathbf Z}(\widetilde C_a)^{\circ}$.
If $(u,L)\in
{\mathbf Z}(\widetilde C_a)^{\circ}$,
the left hand side of (\ref{eqafu})
is $\varphi_a$
by the definition of $\varphi_a$.
By Lemma \ref{lmApL}.2,
the left hand side of (\ref{eqafu})
is $\dfrac{\varphi_a}{[\xi_a:\eta_a]_{\rm insep}}
\cdot [\xi_a:\eta_a]_{\rm insep}
=\varphi_a$.
Hence the equality (\ref{eqafu})
holds on the dense open subset
${\mathbf Z}(\widetilde C)^{\circ}
\subset
{\mathbf Z}(\widetilde C)$.
This also proves the uniqueness of $A$
and implies the independence of $C$.

The left hand side
of (\ref{eqafu}) is 
the function $\varphi_f$
on ${\mathbf Z}(\widetilde C)$
and is constructible and
flat over ${\mathbf G}$
by assumption.
The right hand side
of (\ref{eqafu}) is also a function on 
${\mathbf Z}(\widetilde C)$
flat over ${\mathbf G}$
as is proved in (2)$\Rightarrow$(1).
Hence the equality (\ref{eqafu})
holds for every $(u,L)
\in {\mathbf Z}(\widetilde C)$
by Lemma \ref{lmfla}.2.

We show that the linear combination
$A$ defined for $X\to {\mathbf A}^n
\subset {\mathbf P}^n
={\mathbf P}(E^\vee)$
equals that defined for
$X\to {\mathbf P}(S^dE^\vee)$
for $d\geqq 2$.
Since $E\subset
k[T_1,\ldots,T_n]$ consisting of
polynomials of degree $\leqq 1$
is canonically identified
with a subspace of $S^dE\subset
k[T_1,\ldots,T_n]$ consisting of
polynomials of degree $\leqq d$,
the uniqueness of $A$
implies the independence of $d$.

We show that the equality (\ref{eqafu})
holds for every morphism
$f\colon U\to Y$
with at most
an isolated $C$-characteristic point $u\in U$.
By taking an \'etale morphism
to ${\mathbf A}^1$
defined on a neighborhood
of $f(u)\subset Y$,
we may assume $Y={\mathbf A}^1$
and $f$ is defined by
a function on $U\subset X$.
We may assume that
$f$ is defined by a ratio of
polynomials in $k[T_1,\ldots,T_n]$ 
of degree $d\geqq 1$.
In other words,
$f$ equals a morphism
defined by a pencil
$L$ and the assertion is proved.

2.
We may assume $X$ is affine
and take an immersion $X\to
{\mathbf P}^n$ as in
the proof of 1.
Let $C'_b\subset h^*C$ be an irreducible
component 
and let $C_a\subset C$ be the closure of
its image.
We take a closed point
$(u,L)\in {\mathbf Z}(\widetilde C_a)^\circ$
as in the notation in the proof of 1.\
such that $u=h(v)$ for 
a point $v\in W$.
Then, since $\psi(p_Lh,v)
=\varphi(p_L,u)$,
the coefficient of
$C'_b$ in $A'$ equals
that of
$C_a$ in $A$.
\qed}

\subsection{Characteristic cycle}\label{ssCC}

We state and prove the existence
of characteristic cycle
satisfying the Milnor formula.

\begin{thm}[{\rm cf.{} \cite[Principe p.\ 7]{bp}}]\label{thmM}
Let $X$ be a smooth scheme over 
a perfect field $k$
of characteristic $p>0$ (resp.\ $p=0$)
and ${\cal F}$ be a constructible
complex of $\Lambda$-modules 
of finite tor-dimension on $X$.
Let $C=\bigcup_aC_a$ be a closed conical subset
of the cotangent bundle $T^*X$
such that ${\cal F}$
is micro-supported on $C$.
Assume that every irreducible component of
$X$ and  every irreducible component $C_a$
of $C$ are of dimension $n$.
Then, there exists a unique 
${\mathbf Z}[\frac1p]$-linear 
(resp.\ ${\mathbf Z}$-linear) combination
$CC_C{\cal F}=
\sum_a m_a[C_a]$ satisfying the following
condition:

For every \'etale morphism
$j\colon W\to X$,
every morphism $f\colon W\to Y$
to a smooth curve and 
every at most {\em isolated $C$-characteristic
point} $u\in W$ of $f$, % with respect to $C$,
we have
\begin{equation}
-\dim{\rm tot}\ \phi_u
(j^*{\cal F},f)
=
(CC_C {\cal F},df)_{T^*W,u}.
\label{eqMil}
\end{equation}
Further,
the linear combination
$CC_C{\cal F}$
is independent of $C$
on which ${\cal F}$ is micro-supported.
\end{thm}%\medskip

We will give a proof
by Beilinson
of the fact that
the characteristic cycle 
has ${\mathbf Z}$-coefficients
in Section \ref{sZ}.
This is a generalization of 
the Hasse-Arf theorem
\cite{CL} in the case
$\dim X\leqq 1$. 

\proof{
We may assume $k$ is algebraically
closed by replacing $k$ by
an algebraic closure.
Let $\Lambda_0$ denote
the residue field of the
finite local ring $\Lambda$
and set ${\cal F}_0=
{\cal F}\otimes_\Lambda^L\Lambda_0$.
Then, we have
$\dim{\rm tot} \phi_u
(j^*{\cal F},f)
=
\dim{\rm tot} \phi_u
(j^*{\cal F}_0,f)$
%by Lemma \ref{lmctf}.2
and ${\cal F}_0$ is
micro-supported on $C$
if and only if
and ${\cal F}$ is
micro-supported on $C$
by Lemma \ref{lmmc}.7.
Thus, we may assume
$\Lambda$ is a field.

We regard the left hand side
of (\ref{eqMil})
as a function $\varphi$ on
isolated $C$-characteristic points
in the sense of
Definition \ref{dffla}.1.
In fact,
if the pair
of $j\colon W\to X$ and $f\colon
W\to Y$ is $C$-transversal,
then $f\colon W\to Y$ 
is universally locally acyclic
relatively to $j^*{\cal F}$
and the left hand side of
(\ref{eqMil}) is 0.
If $g\colon Y\to Z$
is an \'etale morphism,
the function $\varphi$ satisfies 
the condition 
$\varphi(f,u)=\varphi(gf,u)$ in (1)
in Proposition \ref{prfla}.1.
If $h\colon W\to X$
is an \'etale morphism,
it also satisfies the condition 
$\varphi(f,hv)=\psi(fh,v)$
in Proposition \ref{prfla}.2.
Thus by Proposition \ref{prfla},
it suffices to show that
$\varphi$ is flat in the sense
of Definition \ref{dffla}.2.

Let the notation be as in (\ref{eqfla})
and we apply 
Proposition \ref{prMsc}.
The morphism $f\colon
U\to Y$
is locally acyclic relatively to
the pull-back of ${\cal F}$
on the complement of
$Z$ by Lemma \ref{lmPC}.2.
The projection
${\rm pr}_2\colon
U\to S$
is locally acyclic relatively to
the pull-back of ${\cal F}$
by the generic universal local acyclicity
\cite[Th\'eor\`eme 2.13]{TF}.
Since $Z$ is quasi-finite over $S$,
the function
$\varphi_f$ is constructible and flat
over $S$ by 
Proposition \ref{prMsc}
and Lemma \ref{lmfla}.1.
\qed}

\begin{df}\label{dfCC}
We define the {\em characteristic cycle}
$CC {\cal F}$ 
to be $CC_C{\cal F}$
independent of $C$
on which ${\cal F}$ is micro-supported.
\end{df}

The Milnor formula \cite{Milnor}
and (\ref{eqMil})
imply that for the constant sheaf $\Lambda$,
we have
$$CC \Lambda
=(-1)^n\cdot [T^*_XX].$$
Thus, the formula (\ref{eqMil})
is a generalization of
the Milnor formula
proved by Deligne in \cite{Milnor}
and shall be also called a Milnor formula.
We will give more examples 
and properties of characteristic cycles 
in the rest of this
subsection and in Section \ref{spb}.
%\subsection{Elementary properties of characteristic cycles}\label{ssep}
We keep assuming that $k$ is
perfect and
$X$ is smooth of dimension $n$ over $k$.
Let ${\cal F}$ be a constructible
complex of $\Lambda$-modules of finite tor-dimension on $X$.

\begin{lm}\label{lmele}
{\rm 1.}
If ${\cal F}$
is locally constant,
we have 
\begin{equation}
CC {\cal F}=(-1)^n
{\rm rank}\ {\cal F}\cdot [T^*_XX].
\end{equation}

{\rm 2.}
For an \'etale morphism
$j\colon U\to X$,
we have
\begin{equation}
CC j^*{\cal F}
=j^*CC {\cal F}.
\end{equation}

{\rm 3.}
Assume that $\dim X=1$
and let $U\subset X$
be a dense open subscheme
where ${\cal F}$ is locally constant.
For $x\in X\sm U$,
let $\bar \eta_x$
denote a geometric generic
point of the strict localization
at a geometric point $\bar x$
above $x$ and let 
\begin{equation}
a_x({\cal F})=
{\rm rank}\ {\cal F}_{\bar \eta_x}
-
{\rm rank}\ {\cal F}_{\bar x}
+
{\rm Sw}_x {\cal F}_{\bar \eta_x}
\label{eqaxF}
\end{equation}
be the Artin conductor.
Then, we have
\begin{equation}
CC {\cal F}=
-\Bigl(
{\rm rank}\ {\cal F}\cdot [T^*_XX]
+
\sum_{x\in X\!\sm \!U}a_x({\cal F})
\cdot [T^*_xX]
\Bigr)
\label{eqdim1}
\end{equation}
\end{lm}

\proof{
1.
It follows from the Milnor formula \cite{Milnor}.
It will also follow immediately
from the compatibility of
the characteristic cycles
with smooth pull-back
Proposition \ref{prsm*}.

2.
Since the characterization (\ref{eqMil})
is an \'etale local condition,
the assertion follows.

3.
By Lemma \ref{lmlcst}.2,
it suffices to determine the coefficients.
For the $0$-section $T^*_XX$,
it follows from 1 and 2.
For the fibers,
since 
$\dim{\rm tot}_x\phi_x
({\cal F},{\rm id})
=
a_x({\cal F})$,
it follows from the Milnor formula (\ref{eqMil})
for the identity $X\to X$.
\qed}
\medskip

For surfaces,
the characteristic cycle
is studied in \cite{surface}.

\begin{df}\label{dfi*A}
Let $i\colon X\to Y$
be a closed immersion 
of smooth schemes over $k$
and let 
\begin{equation}
\begin{CD}
T^{*}X@<<< X\times_{Y}T^{*}Y
@>>> T^{*}Y
\end{CD}
\label{eqi*A}
\end{equation}
be the canonical morphisms.
Let $C\subset T^*X$
be a closed conical subset.
Assume that every irreducible component of $X$
and every irreducible component 
$C_a$ of $C=\bigcup_aC_a$
are of dimension $n$
and 
that every irreducible component of $Y$
is of dimension $m$.
Then, for a linear combination
$A=\sum_am_a[C_a]$,
we define $i_*A$ to be 
$(-1)^{n-m}$-{\em times} the push-forward
by the second arrow
$X\times_{Y}T^{*}Y\to T^{*}Y$
in {\rm (\ref{eqi*A})}
of the pull-back of $A$ by the first arrow
$X\times_{Y}T^{*}Y\to T^{*}X$
in the sense of intersection theory.
\end{df}

\begin{lm}\label{lmRf}
Let ${\cal F}$ be a constructible complex
of $\Lambda$-modules 
of finite tor-dimension on $X$.

{\rm 1.}
For a distinguished triangle
$\to {\cal F}'\to {\cal F}\to {\cal F}''\to$
in $D_{\rm ctf}(X,\Lambda)$,
we have
\begin{equation}
CC {\cal F}=
CC {\cal F}'+
CC {\cal F}''.
\label{eqadd}
\end{equation}

{\rm 2.}
For a closed immersion 
$i\colon X\to Y$
of smooth schemes over $k$,
we have
\begin{equation}
CC i_*{\cal F}=
i_*CC {\cal F}.
\label{eqi*F}
\end{equation}

{\rm 3.}
For a morphism $f\colon X\to Y$
of separated 
smooth schemes of finite type over $k$, 
we have
\begin{equation}
CC Rf_*{\cal F}=
CC Rf_!{\cal F}.
\label{eqRf}
\end{equation}

{\rm 4.}
We have
\begin{equation}
CC D_X{\cal F}=
CC {\cal F}.
\label{eqChD}
\end{equation}
\end{lm}

\proof{
1.
By the characterization
of characteristic cycle
by the Milnor formula (\ref{eqMil}),
it follows from the additivity
of the total dimension.

2.
Let $C\subset T^*X$ be the singular
support of ${\cal F}$.
Then,
$i_*{\cal F}$ is micro-supported
on $i_{\circ}C\subset T^*Y$
by Lemma \ref{lmmc}.6.
Let $Y\to {\mathbf P}$
be an immersion satisfying the condition
(E) before Lemma \ref{lmloc}
and the condition (C) before
Proposition {\rm \ref{prwi}}
for $i_{\circ}C$.
Then, by the description of
the characteristic cycle
$CC {\cal F}
=CC_C^E {\cal F}$
in (\ref{eqChE}),
it follows from the canonical isomorphism
$\phi({\cal F},p_L^{\circ} \circ i)
\to \phi(i_*{\cal F},p_L^{\circ} )$
for the morphism
$p_L^{\circ} \colon Y_L^{\circ}\to L$
defined by a pencil $L$.

3.
By 1,
it follows from \cite{La}.

4. 
We have $SS{\cal F}=SSD_X{\cal F}$
by Corollary \ref{corssD}.
By 2 and Lemma \ref{lmele}.2,
we may assume $X$ is projective
as in the proof of Corollary \ref{corssD}.
Let $C=SS{\cal F}=SSD_X{\cal F}$
be the singular support.
Let $X\to {\mathbf P}$
be a closed immersion
satisfying the condition
(E) and (C) before 
and after Lemma \ref{lmloc}.
Then, for a point
$(u,L)$
in the dense open subset
${\mathbf P}(\widetilde C)^{\circ}
\subset
{\mathbf P}(\widetilde C)
\subset X\times_{\mathbf P}Q$
and $v=p_L(u)\in L$,
we have
$\dim{\rm tot}\phi_u({\cal F},p_L^{\circ})
=
a_v(Rp^\vee_*p^*{\cal F})|_L$
and similarly for $D_X{\cal F}$.
Since
$a_v(Rp^\vee_*p^*{\cal F})|_L
=
a_vD_L(Rp^\vee_*p^*{\cal F})|_L
=
a_v(Rp^\vee_*p^*D_X{\cal F})|_L$,
it follows from the description of
the characteristic cycle
$CC {\cal F}
=CC_C^E {\cal F}$
in (\ref{eqChE}).
\qed}
\medskip

For the residue field $\Lambda_0$ of
$\Lambda$,
a constructible complex
${\cal F}$ of $\Lambda$-modules
of finite tor-dimension on $X$
is a perverse sheaf if and only if
${\cal F}\otimes^L_\Lambda\Lambda_0$
is a perverse sheaf.

\begin{pr}\label{prperv}
Assume ${\cal F}$ is a perverse sheaf on $X$.

{\rm 1.} {\rm (\cite[Question p.\ 7]{bp})}
We have
\begin{equation}
CC {\cal F}\geqq0
\label{eqperv}
\end{equation}.

{\rm 2.}
The support of
$CC {\cal F}$
equals $SS{\cal F}$.
\end{pr}

\proof{
By the description of
the characteristic cycle
$CC {\cal F}
=CC_C^E {\cal F}$
in (\ref{eqChE}),
it follows from Lemma \ref{lmRn}.
\qed}%\medskip

\begin{cor}\label{corRj}
Let $j\colon U=X\sm D\to X$
be the open immersion of
the complement of a Cartier divisor.
Then, for a perverse sheaf
${\cal F}$ of $\Lambda$-modules on $U$,
we have
\begin{equation}
SS Rj_*{\cal F}=
SS j_!{\cal F}.
\label{eqSSRj}
\end{equation}
\end{cor}

\proof{
Since the open immersion
$j\colon U=X\sm D\to X$ is
an affine morphism,
$Rj_*{\cal F}$ and 
$j_!{\cal F}$ are perverse sheaves on $X$
by \cite[Corollaire 4.1.10]{BBD}.
Hence, it follows from
Lemma \ref{lmRf}.3 and Proposition \ref{prperv}.2.
\qed}
\medskip

We show the compatibility with
smooth pull-back.

\begin{df}\label{dfCsm}
Let $h\colon W\to X$ be 
a smooth morphism of
smooth schemes over
a perfect field $k$ and let
$C\subset T^*X$ be a closed conical
subset.
Assume that every irreducible
component of $X$ and
every irreducible
component of $C$
are of dimension $n$
and that
every irreducible
component of $W$ is of dimension $m$.
Let 
\begin{equation}
\begin{CD}
T^*W@<<< W\times_XT^*X
@>>> T^*X
\end{CD}
\label{eqhsm}
\end{equation}
be the canonical morphisms.
Then, for a linear combination
$A=\sum_am_a[C_a]$
of irreducible components of
$C=\bigcup_aC_a$,
we define $h^!A$ to be 
$(-1)^{n-m}$-{\em times} the push-forward
by the first arrow
$W\times_XT^*X\to T^*W$
in {\rm (\ref{eqh!A})}
of the pull-back of $A$ by the second arrow
$W\times_XT^*X\to T^{*}X$
in the sense of intersection theory.
\end{df}

\begin{pr}\label{prsm*}
Let $h\colon W\to X$ be a smooth morphism of 
smooth schemes over a perfect field $k$ and let
$C\subset T^*X$ be a closed conical
subset.
Assume that every irreducible
component of $X$ and
every irreducible
component of $C$
are of dimension $n$
and that
every irreducible
component of $W$ is of dimension $m$.

Let ${\cal F}$ be
a constructible complex of
$\Lambda$-modules on $X$
of finite tor-dimension
micro-supported on
$C\subset T^*X$.
Then, we have 
\begin{equation}
CC h^*{\cal F}=
h^!CC {\cal F}.
\label{eqsm*}
\end{equation}
\end{pr}

\proof{
Since the assertion is \'etale local
on $W$, we may assume
$W=X\times {\mathbf A}^n$
and $h$ is the first projection.
Then, this is the case
where $Y={\mathbf A}^n$
and ${\cal G}=\Lambda$
of \cite[Theorem 3.6.2]{notes},
which is proved 
using the Thom-Sebastiani formula
\cite{TS}
without
using the results in the rest of
this article.
\qed}

\subsection{Integrality}\label{sZ}

In this section, we give a proof
by Beilinson
of the integrality of characteristic cycles.

\begin{thm}[Deligne]\label{thmZ}
The coefficients
of $CC{\cal F}$ are integers.
\end{thm}

This is a generalization of 
the Hasse-Arf theorem
\cite{CL} in the case
$\dim X\leqq 1$. 
We will deduce Theorem from
the Milnor formula
and the following Proposition.

\begin{pr}[{\rm \cite[Proposition 4.12]{Be}}]\label{prT}
Let $X$ be a smooth scheme
of dimension $n$
over an algebraically closed
field $k$
and $C\subset T^*X$
be a closed irreducible conical subset
of dimension $n$.
Let $u\in X$ be a closed
point and
let $(u,\omega)\in C$
be a closed smooth point of
$C\subset T^*X$ regarded as
a reduced closed subscheme.

Then,
there exists a function
$f$ defined on a neighborhood of $u$
such that the section $df$
of $T^*X$ meets $C$ transversely
at $(u,\omega)$
if one of the following conditions {\rm (1)}
and {\rm (2)} is satisfied:

{\rm (1)}
The characteristic $p$ of $k$ is
different from $2$.

{\rm (2)}
Let $T$ be the tangent space
$T_{(u,\omega)}C$ of $C$ at 
the smooth point $(u,\omega)$.
Then, 
there exists a function
$g$ defined on a neighborhood of $u$
such that the section $dg$
of $T^*X$ meets $C$ at 
$(u,\omega)$
and that
the dimension of
the intersection
$(dg)_*(T_uX)\cap T\subset
T_{(u,\omega)}(T^*X)$
is even.
\end{pr}

\proof{%[Proof of Proposition {\rm \ref{prT}}]{
Take a local coordinate
$x_1,\ldots,x_n$ of $X$ at $u$
and 
write $\omega=\sum_ia_idx_i$
as a $k$-linear combination
of the basis $dx_1,\ldots,dx_n$ of
the cotangent space $T^*_uX$.
Then, the function $g=\sum_ia_ix_i$
satisfies $dg(u)=\omega$.
We will modify 
$g$ as $f=g+\sum_{i,j}b_{ij}x_ix_j$
to find a function $f$ satisfying
the required condition.

Let $V=T_uX$
and
$W=T_{(u,\omega)}(T^*X)$
denote the tangent spaces
at $u\in X$ and at
$(u,\omega)\in T^*X$ respectively.
The tangent space $W$
is decomposed as the direct sum of
the image $(dg)_*(V)$ of
the morphism $(dg)_*\colon V\to W$
defined by the section $dg\colon X\to T^*X$
with the tangent space
$T_\omega(T^*_uX)$ of the fiber
$T^*_uX\subset T^*X$.
The latter $T_\omega(T^*_uX)$
is naturally identified
with the cotangent space
$T^*_uX\subset T^*X$ 
that is the dual
$V^\vee$ of $V=T_uX$.
Further identifying $V$ with
its image $(dg)_*(V)$, we identify
$W=V\oplus V^\vee$.

The required transversality condition means that
the intersection 
$(df)_*(V)\cap T\subset W$ is $0$.
Since $df=dg+\sum_{i,j}(b_{ij}+b_{ji})x_idx_j$,
we have $(df)_*=(dg)_*+(B+B^\vee)$,
where $B\colon V\to V^\vee$
denotes the bilinear form on $V$
defined by the matrix $(b_{ij})$
with respect to the basis $(dx_i)$.
Consequently, under the identification
$W=V\oplus V^\vee$
above, the image of
$(df)_*\colon V\to W$
is identified with the graph $\Gamma$ of
$A=B+B^\vee\colon V\to V^\vee$.
Thus, the assertion is
a consequence of the
following lemma on linear algebra.
\qed}

\begin{lm}\label{lmWV}
Let $V$ be a $k$-vector space
of finite dimension
and $V^\vee$ be the dual.

{\rm 1.}
Let $T\subset W=V\oplus V^\vee$
be a linear subspace of
$\dim T=\dim V$
and set
$V_1=V\cap T$.
Then there exists
a direct sum decomposition
$V=V_1\oplus V_2$
satisfying the following property:

For a non-degenerate bilinear form
$A_1\colon V_1\to V_1^\vee$,
we extend it as $A=A_1\oplus 0 
\colon V=V_1\oplus V_2\to 
V^\vee=V_1^\vee\oplus V_2^\vee$
and let $\Gamma\subset
W=V\oplus V^\vee$ denote the graph.
Then we have
$\Gamma\cap T=0$.
% and $\Gamma\cap V^\vee=0$.

{\rm 2.}
Either if $\dim V$ is even or
if the characteristic $p$ of $k$ is
different from $2$,
there exists a bilinear form
$B\colon V\to V^\vee$ such that
$A=B+B^\vee$ is non-degenerate.
\end{lm}

\proof{1.
Let $\bar T$ denote the
image of the morphism $T\to 
W\to W/V=V^\vee$
and $\bar T^\perp\subset V$
be the orthogonal subspace.
Then by the assumption $\dim T=\dim V$,
the subspaces 
$V_1=V\cap T$
and $\bar T^\perp\subset V$
have the same dimension.
Hence, there exists a
direct sum decomposition
$V=V_1\oplus V_2$ satisfying
$V_2\cap \bar T^\perp=0$.
Since $V_2=(V_1^\vee)^\perp$,
we have $V_1^\vee+\bar T=V^\vee$.

By the assumption that
$A$ is non-degenerate,
we have $\Gamma+T\supset V_1+V_1^\vee+V_2$.
Hence, further
by $V_1^\vee+\bar T=V^\vee$,
we obtain
$\Gamma+T=W$.
Thus $\dim \Gamma=\dim V=\dim T$
implies 
$\Gamma\cap T=0$.

2.
If $p\neq 2$, it suffices
to take a non-degenerate symmetric 
bilinear form $B$.
If $p=2$ and if $\dim V$ is even,
it suffices to take
a non-degenerate alternating 
bilinear form $A$
and to write $A=B+B^\vee$
by taking a symplectic basis.
\qed}

\proof[Proof of Theorem {\rm \ref{thmZ}}
{\rm (Beilinson)}]{
We may assume $k$ is algebraically closed.
Write $C=\bigcup_aC_a$
as the union of irreducible components.
We show that each
coefficient $m_a$ in
$CC{\cal F}=\sum_am_a[C_a]$ is
an integer.
Let $(u,\omega)$ be
a smooth point of $C_a$
not contained in any other 
irreducible component
$C_b, (b\neq a)$
of $C$.

If one of the conditions
(1) and (2) in Proposition \ref{prT}
is satisfied, 
for $f$ as in Proposition \ref{prT},
the coefficient
$m_a$ equals
$(CC{\cal F},df)_{T^*X,u}$
since 
$(C_a,df)_{T^*X,u}=1$ and
$(C_b,df)_{T^*X,u}=0$ for $b\neq a$.
By the Hasse-Arf theorem
\cite{CL},
the total dimension
$\dim{\rm tot} \phi_u({\cal F},f)$
is an integer and
the Milnor formula
$(CC{\cal F},df)_{T^*X,u}=
-\dim{\rm tot} \phi_u({\cal F},f)$
(\ref{eqMil})
implies that the coefficient
$m_a=(CC{\cal F},df)_{T^*X,u}$
is an integer in this case.

In the exceptional case in $p=2$,
we take $X\times {\mathbf A}^1$
and the pull-back ${\rm pr}_1^*{\cal F}$.
If the original $g$ does
not satisfies the condition (2)
in Proposition \ref{prT}
at a smooth point $(u,\omega)$
of $C_a$,
then the composition
${\rm pr}_1^*g$ satisfies
the condition
at the smooth point $((u,0),{\rm pr}_1^*\omega)$
of ${\rm pr}_1^oC_a
\subset T^*(X\times {\mathbf A}^1)$.
Thus, the assertion follows from
the compatibility Proposition \ref{prsm*}
of characteristic cycle with smooth pull-back.
\qed}

\section{Characteristic class}\label{scc}

\subsection{Cycle classes
on projective space bundles}\label{ssCh}

In this preliminary subsection,
we recall some basic facts
on the Chow groups of
${\mathbf P}^n$-bundles.
In this section,
we assume that $X$ is a
scheme of finite type over a field $k$.
We can replace this assumption
by some condition
which assures necessary properties
on Chow groups.

To describe the Chow group
of a projective space bundle, 
we introduce some notation.
Let $A=\bigoplus_iA_i$
be a graded module
and let $c_q, q=0,\ldots, n+1$
be endomorphisms of
$A$ of degree $-q$
sending $A_i$ to $A_{i-q}$.
We assume that $c_0$ is the identity.
We formally set $f=\sum_{q=0}^{n+1}
c_qh^{n+1-q}$ and 
define a graded module
\begin{equation}
A[h]/(f)=A^{\oplus n+1}=
\bigoplus_{q=0}^n
Ah^q
\label{eqAhf}
\end{equation}
where $A_ih^q$
has degree $i+n-q$.
We define an endomorphism
$h$ of $A[h]/(f)$ of degree $-1$
by sending
on $ah^q$ to $ah^{q+1}$
for $a\in A$ and $q<n$
and sending $a\cdot h^n$ to
$-(\sum_{q=1}^{n+1}c_qa\cdot h^{n+1-q})
=(h^{n+1}-f)\cdot a$.

\begin{lm}\label{lmCHP}
Let $X$ be a scheme of
finite type over $k$.
Let $E$ be a vector bundle
over $X$ of rank $n+1$
and let $p\colon {\mathbf P}(E)\to X$
be the associated ${\mathbf P}^n$-bundle.
Let $c_q(E)$ denote the $q$-th Chern classes of
$E$ and let $c_h(E)$ denote $\sum_{q=0}^{n+1}
c_q(E)h^{n+1-q}$.

{\rm 1 (\cite[Theorem 3.3 (b)]{Ful}).}
The morphism
\begin{equation}
{\rm CH}_\bullet (X)[h]/(c_h(E))
\to
{\rm CH}_\bullet ({\mathbf P}(E))
\label{eqCHP}
\end{equation}
sending $ah^q$
to $c({\cal O}(1))^q\cap p^*a$
is an isomorphism
of graded modules.

{\rm 2.}
Let $i\geqq 0$ be an integer.
Then, the inverse of the degree $i$-part
\begin{equation}
\bigoplus_{q=0}^n
{\rm CH}_{i-n+q}(X)
=
\bigl({\rm CH}_\bullet(X)[h]/(c_h(E))\bigr)_i
\to 
{\rm CH}_i ({\mathbf P}(E))
\label{eqCHPn}
\end{equation}
of the isomorphism {\rm (\ref{eqCHP})}
sends $b\in {\rm CH}_i ({\mathbf P}(E))$
to $c(E)\cap p_*(c_1({\cal O}(-1))^{-1}b)$.
\end{lm}

\proof{
2.
By the isomorphism
(\ref{eqCHP}),
we identify ${\rm CH}_{\bullet}({\mathbf P}(E))$
and its degree $i$-part 
with ${\rm CH}_\bullet (X)[h]/
(c_h(E))$ and $\bigoplus_{q=0}^n
{\rm CH}_{i-n+q}(X)h^q$.
By the projection formula, the morphism
$p_*\colon {\rm CH}_{\bullet}({\mathbf P}(E))
\to {\rm CH}_{\bullet}(X)$
is the morphism
${\rm CH}_\bullet (X)[h]/
(c_h(E))\to {\rm CH}_\bullet (X)$
taking the coefficient of $h^n$.

We may assume $b=ah^q$ for $q=0,\ldots,n$
and $a\in {\rm CH}_{i-n+q}(X)$.
The difference 
$c(E)(1-h)^{-1}$ $h^q
-c_h(E)(1-h)^{-1}
=\sum_{i=0}^{n+1}c_i(E)(h^q-h^{n+1-i})/(1-h)$
is of degree $\leqq n$ in $h$
and the top term is $h^n$.
Since $c_h(E)$ acts as $0$ on 
${\rm CH}_\bullet (X)[h]/(c_h(E))$,
the assertion follows.
\qed}

\begin{lm}\label{lmincl}
Let $i\colon F\to E$ be an injection
of vector bundles on $X$
of ${\rm rank}\ E=n+1\geqq
{\rm rank}\ F=m+1$
and let $i\colon 
{\mathbf P}(F)
\to
{\mathbf P}(E)$
also denote the induced morphism.
Let $E'$ denote the cokernel.

{\rm 1.}
The diagram 
\begin{equation}
\begin{CD}
{\rm CH}_\bullet({\mathbf P}(F))
@<{i^*}<<
{\rm CH}_\bullet({\mathbf P}(E))\\
@AAA@AAA\\
{\rm CH}_\bullet(X)[h]/(c_h(F))
@<{\rm can}<<
{\rm CH}_\bullet(X)[h]/(c_h(E))
\end{CD}
\label{eqincl}
\end{equation}
is commutative.

{\rm 2.}
The diagram 
\begin{equation}
\begin{CD}
{\rm CH}_\bullet({\mathbf P}(F))
@>{i_*}>>
{\rm CH}_\bullet({\mathbf P}(E))\\
@AAA@AAA\\
{\rm CH}_\bullet(X)[h]/(c_h(F))
@>{c_h(E')\cap }>>
{\rm CH}_\bullet(X)[h]/(c_h(E))
\end{CD}
\label{eqincl*}
\end{equation}
is commutative.
\end{lm}

\proof{
1.
Since the pull-back of ${\cal O}(1)$
by ${\mathbf P}(F)\to {\mathbf P}(E)$
is also
${\cal O}(1)$ on ${\mathbf P}(F)$,
the assertion follows from the definition.

2.
By 1, it suffices to
show that the endomorphism $i_*i^*$
of
${\rm CH}_\bullet({\mathbf P}(E))$
is identified with the multiplication by $c_h(E')$
on ${\rm CH}_\bullet(X)[h]/(c_h(E))$.
By the self-intersection formula,
the endomorphism $i_*i^*$
equals the action of
the top Chern class
$c_{n-m}
(T_{{\mathbf P}(F)}{\mathbf P}(E))$
of the normal bundle.

By the exact sequence
$0\to \Omega_{{\mathbf P}(E)/X}
\to E^\vee\otimes {\cal O}(-1)\to {\cal O}_{{\mathbf P}(E)}\to 0$
and the corresponding one for
${\mathbf P}(F)$,
we obtain an isomorphism
$N_{{\mathbf P}(F)/{\mathbf P}(E)}
\to 
E^{\prime \vee}\otimes {\cal O}(-1)$
for the conormal sheaf.
Hence, the top Chern class
$c_{n-m} (T_{{\mathbf P}(F)}{\mathbf P}(E))$
equals 
$c_{n-m}(E'\otimes {\cal O}(1))$
and is identified with $c_h(E')$.
\qed}

\begin{lm}\label{lmEF1}
Assume
$E=F\oplus {\mathbf A}^1_X$.
Let $i\colon 
{\mathbf P}(F)
\to
{\mathbf P}(E)$ denote the injection
and let $p\colon 
{\mathbf P}(E)\to X$ be the projection.

{\rm 1.}
The morphism
\begin{equation}
i^*\oplus p_*\colon
{\rm CH}_\bullet({\mathbf P}(E))
\to
{\rm CH}_\bullet({\mathbf P}(F))
\oplus
{\rm CH}_\bullet(X)
\label{eqQZ}
\end{equation}
is an isomorphism.
If $X$ is irreducible
of dimension $d\leqq n$
and $i_x\colon x\to X$ is
the immersion of a smooth
$k$-rational point, 
the second projection
$p_*\colon
{\rm CH}_d({\mathbf P}(E))
\to
{\rm CH}_d(X)=
{\mathbf Z}$
equals the pull-back $i_x^*
\colon
{\rm CH}_d({\mathbf P}(E))
\to
{\rm CH}_0({\mathbf P}^n)=
{\mathbf Z}$.

{\rm 2.}
The morphism
\begin{equation}
 i_*+ p^*\colon
{\rm CH}_\bullet({\mathbf P}(F))
\oplus
{\rm CH}_\bullet(X)
\to
{\rm CH}_\bullet({\mathbf P}(E))
\label{eqQZ2}
\end{equation}
is an isomorphism.
The composition 
${\rm CH}_\bullet ({\mathbf P}(E))
\to {\rm CH}_{\bullet}(X)$
of the inverse of
the isomorphism {\rm (\ref{eqQZ2})}
with the second projection
equals the pull-back by
the $0$-section
$s\colon X\to F\subset {\mathbf P}(E)$.
\end{lm}

\proof{1.
We identify
${\rm CH}_\bullet({\mathbf P}(E))
=
{\rm CH}_\bullet(X)[h]/(c_h(E))$
and 
${\rm CH}_\bullet({\mathbf P}(F))
=
{\rm CH}_\bullet(X)[h]/(c_h(F))$.
Then, $c_h(E)=c_h(F)h$
and the morphism
$i^*\colon
{\rm CH}_\bullet({\mathbf P}(E))
\to
{\rm CH}_\bullet({\mathbf P}(F))$
is identified with the surjection
${\rm CH}_\bullet(X)[h]/(c_h(E))
\to
{\rm CH}_\bullet(X)[h]/(c_h(F))$
sending $h^n$ to $0$
by Lemma \ref{lmincl}.1.
Since $p_*\colon
{\rm CH}_\bullet({\mathbf P}(E))
=
{\rm CH}_\bullet(X)[h]/(c_h(E))
\to
{\rm CH}_\bullet(X)$
is the morphism taking
the coefficient of
$h^n$,
the morphism (\ref{eqQZ}) is
an isomorphism.

The second assertion follows from
the commutative diagram
$$\begin{CD}
{\rm CH}_d({\mathbf P}(E))
@>{i_x^*}>>
{\rm CH}_0({\mathbf P}^n)\\
@V{p_*}VV@VV{\deg}V\\
{\rm CH}_d(X)
@>{i_x^*}>>
{\rm CH}_0(x).
\end{CD}$$

2.
The morphism
$i_*\colon
{\rm CH}_\bullet({\mathbf P}(F))
\to
{\rm CH}_\bullet({\mathbf P}(E))$
is identified with the multiplication
$h\cap \colon
{\rm CH}_\bullet(X)[h]/(c_h(F))
\to
{\rm CH}_\bullet(X)[h]/(c_h(E))$
by Lemma \ref{lmincl}.2.
Hence
the morphism (\ref{eqQZ2}) is
an isomorphism.

Since $s^*i_*=0$ and
$s^*p^*$ is the identity of
${\rm CH}_\bullet(X)$,
the second assertion follows
from the isomorphism (\ref{eqQZ2}).
\qed}

\begin{lm}\label{lmsurj}
Let $\theta\colon E\to F$ be a surjection
of vector bundles on $X$
of ${\rm rank}\ E=n+1\geqq
{\rm rank}\ F=m+1$ and
let $K$ denote the kernel.
Let $\pi\colon 
{\mathbf P}(E)'
\to
{\mathbf P}(E)$
be the blow-up at
${\mathbf P}(K)
\subset
{\mathbf P}(F)$
and let $\theta'\colon 
{\mathbf P}(E)'
\to
{\mathbf P}(F)$
denote the morphism induced by $\theta$.

{\rm 1.}
Let 
$c_h(K)^{-1}\cap 
\colon
{\rm CH}_\bullet(X)[h]/(c_h(E))
\to
{\rm CH}_\bullet(X)[h]/(c_h(F))$
denote the composition 
of the first projection
with the inverse of the isomorphism
$$(c_h(K)\cap)+ {\rm can}\colon
{\rm CH}_\bullet(X)[h]/(c_h(F))
\oplus
{\rm CH}_\bullet(X)[h]/(c_h(K))
\to
{\rm CH}_\bullet(X)[h]/(c_h(E)).$$
Then, the diagram 
\begin{equation}
\begin{CD}
{\rm CH}_\bullet({\mathbf P}(E))
@>{\theta'_*\pi^*}>>
{\rm CH}_\bullet({\mathbf P}(F))\\
@AAA@AAA\\
{\rm CH}_\bullet(X)[h]/(c_h(E))
@>{c_h(K)^{-1}\cap }>>
{\rm CH}_\bullet(X)[h]/(c_h(F))
\end{CD}
\label{eqsurj*}
\end{equation}
is commutative.

{\rm 2.}
The diagram 
\begin{equation}
\begin{CD}
{\rm CH}_\bullet({\mathbf P}(E))
@<{\pi_*\theta^{\prime *}}<<
{\rm CH}_\bullet({\mathbf P}(F))\\
@AAA@AAA\\
{\rm CH}_\bullet(X)[h]/(c_h(E))
@<{\rm can}<<
{\rm CH}_\bullet(X)[h]/(c_h(F))
\end{CD}
\label{eqsurj}
\end{equation}
is commutative.
\end{lm}

\proof{
Let $A={\rm CH}_\bullet(X)$, $c_h(F)=f$ and $c_h(E)=f\cdot g$
and identify 
${\rm CH}_\bullet({\mathbf P}(E))
=A[h]/(f\cdot g)$ and
${\rm CH}_\bullet({\mathbf P}(F))
=A[h']/(f)$.
Let $L\subset {\mathbf P}(F)\times_XF$
be the universal sub line bundle
and let
$V\subset {\mathbf P}(F)\times_XE$
be its pull-back by the base change of $\theta$.
Then, 
${\mathbf P}(E)'
\to
{\mathbf P}(F)$
is canonically identified with
the ${\mathbf P}^{n-m}$-bundle
${\mathbf P}(V)$
by Lemma \ref{lmFEL}.
Thus, we identify
${\rm CH}_\bullet({\mathbf P}(E)')
=A[h']/(f)[h]/((h-h')\cdot g)$.
Let $e$ denote the polynomial
in $h$ and $h'$ defined
by $f(h)-f(h')=(h-h')e$.

We have a cartesian diagram
$$\begin{CD}
{\mathbf P}(K)
\times 
{\mathbf P}(F)
@>>>
{\mathbf P}(E)'\\
@VVV@VV{\pi}V\\
{\mathbf P}(K)
@>>>
{\mathbf P}(E).
\end{CD}$$
We show that this induces a cocartesian
diagram
\begin{equation}
\begin{CD}
A[h]/(g(h))[h']/(f(h'))
@>{(h-h')\times}>>
A[h']/(f(h'))[h]/((h-h')\cdot g(h))\\
@A{e\times}AA@AAA\\
A[h]/(g(h))
@>{f(h)\times}>>
A[h]/(f(h)\cdot g(h)).
\end{CD}
\label{eqPE'}
\end{equation}
on the Chow groups. 
The diagram is cocartesian by
\cite[Proposition 6.7 (e)]{Ful}.
By Lemma \ref{lmincl}.2,
the lower horizontal arrow is
the multiplication by $f(h)$.
The descriptions of
the left vertical arrow
and the upper horizontal arrow
follow similarly from
the excess intersection formula
and the self-intersection formula
respectively.
Since $f(h')=0$ 
and hence $f(h)=(h-h')e$ in
$A[h']/(f(h'))[h]/((h-h')\cdot g(h))$,
the right vertical arrow
is well-defined as
the canonical morphism.

1.
The morphism $\theta'_*\pi^*$
is the composition of
the right vertical arrow in (\ref{eqPE'})
with the morphism
$A[h']/(f(h'))[h]/((h-h')\cdot g(h))\to
A[h']/(f(h'))$
taking the coefficient of $h^{n-m}$.
For $i< n-m$,
we have $\theta'_*\pi^* h^i=0$.
For $i\geqq 0$,
we have $h^ig(h)=
h^{\prime i}g(h)$ in
$A[h']/(f(h'))[h]/((h-h')\cdot g(h))$
and hence $\theta'_*\pi^* 
h^ig(h)=h^{\prime i}$.
Thus, 
the commutativity of
(\ref{eqsurj*}) is proved.

2.
By the cocartesian diagram
(\ref{eqPE'}),
we identify 
$A[h']/(f(h'))[h]/((h-h')\cdot g(h))$
with the amalgamated sum.
Then, the morphism
$\pi_*\colon
A[h']/(f(h'))[h]/((h-h')\cdot g(h))
\to
A[h]/(f(h)\cdot g(h))$
is induced by the
identity of $A[h]/(f(h)\cdot g(h))$
and 
the morphism
$A[h]/(g(h))[h']/(f(h'))
\to 
A[h]/(g(h))$
taking the coefficient of $h^{\prime m}$.
The morphism $\pi_*\theta^{\prime *}$
is the composition of
the canonical morphism
$A[h']/(f(h'))\to
A[h']/(f(h'))[h]/((h-h')$
with the above morphism.
For $0\leqq i\leqq m$,
we have $h^{\prime i}
=h^i+(h'-h)(\deg <m)$.
Hence, we have
$\pi_*\theta^{\prime *}h^{\prime i}
=h^i$.
Thus, 
the commutativity of
(\ref{eqsurj}) is proved.
\qed}

\begin{lm}\label{lmtheta}
Let $\theta\colon E\to F$ be
a surjection of vector bundles
of rank $n\geqq m$ over $X$
and 
let $K={\rm Ker}(\theta
\colon E\to F)$ be the kernel of
$\theta$.

Let $A=\sum_am_aC_a$
be a ${\mathbf Z}$-linear combination
of irreducible closed
conical subsets $C_a\subset E$
of dimension $i\leqq m$
and define
$\bar A
=
\sum_am_a\bar C_a
\in {\rm CH}_i(
{\mathbf P}(E\oplus {\mathbf A}^1_X))
=
{\rm CH}_{\bullet\leqq i}(X)$.
Assume that the intersection
$C\cap K$ of the union $C=\bigcup_aC_a$
with $K$
is a subset of the $0$-section.
Define a linear combination $\theta(A)$ 
of closed conical subsets of $F$ to be
the direct image of $A$ by $\theta\colon
E\to F$ which is finite on $C$.
Then,
we have
\begin{equation}
[\overline {\theta(A)}]=
c_h(K)^{-1}\cap
[\bar A]
\label{eqPhC}
\end{equation}
in ${\rm CH}_i(
{\mathbf P}(F\oplus {\mathbf A}^1_X))
=
{\rm CH}_{\bullet\leqq i}(X)$.
\end{lm}

\proof{
Let $\pi\colon
{\mathbf P}(E\oplus {\mathbf A}^1_X)'
\to
{\mathbf P}(E\oplus {\mathbf A}^1_X)$
and $\theta'\colon
{\mathbf P}(E\oplus {\mathbf A}^1_X)'
\to
{\mathbf P}(F\oplus {\mathbf A}^1_X)$
be as in Lemma \ref{lmsurj}.
The assumption that
$C\cap K$ 
is a subset of the $0$-section
means that ${\mathbf P}(K)$ does not
meet the closure 
$\bar C\to {\mathbf P}(E\oplus {\mathbf A}^1_X)$.
Hence we have $[\overline {\theta(A)}]
= \theta'_*\pi^* ([\bar A])$.
Thus,  the assertion follows
by Lemma \ref{lmsurj}.1.
\qed}

\subsection{Characteristic class
%and its functorial properties
}\label{ssccc}

In this subsection,
$k$ denotes a perfect field of
characteristic $p\geqq 0$
and $\Lambda$ denotes
a finite field of characteristic
$\neq p$.

Let $X$ be a scheme of
finite type over $k$
and $i\colon X\to M$
be a closed immersion 
to a smooth scheme $M$ 
of dimension $n$ over $k$.
We identify
${\rm CH}_n({\mathbf P}(X\times_MT^*M
\oplus {\mathbf A}^1_X))$
with
${\rm CH}_\bullet(X)
=\bigoplus_{i=0}^n
{\rm CH}_i(X)$
by the canonical isomorphism
\begin{equation}
{\rm CH}_\bullet(X)
=\bigoplus_{i=0}^n
{\rm CH}_i(X)
\to 
{\rm CH}_n({\mathbf P}(X\times_MT^*M
\oplus {\mathbf A}^1_X))
\label{eqCCH}
\end{equation}
sending $(a_i)_i$
to $\sum_ip^*a_ih^i$
where $p\colon 
{\mathbf P}(X\times_MT^*M
\oplus {\mathbf A}^1_X)
\to X$ is the canonical projection
and 
$h=c_1(L^\vee)$
be the first Chern class of the dual 
$L^\vee$ of
the universal sub line bundle
$L\subset
{\mathbf P}(X\times_MT^*M
\oplus {\mathbf A}^1_X)\times_X
(X\times_MT^*M
\oplus {\mathbf A}^1_X)$.
The inverse mapping sends $a$
to $c(T^*X)\cap p_*(c(L)^{-1}
\cap a)$.

\begin{lm}\label{lmwell}
Let $X$ be a scheme of
finite type over $k$
and ${\cal F}$ be a constructible
complex of $\Lambda$-modules on $X$.
Let $i\colon X\to M$
be a closed immersion 
to a smooth scheme $M$ 
of dimension $n$ over $k$.
Then, the class of the closure
$\overline{CCi_*{\cal F}}
={\mathbf P}(
CCi_*{\cal F}\oplus {\mathbf A}^1_X)$
regarded as an element of
${\rm CH}_n(X\times_MT^*X
\oplus {\mathbf A}^1_X)=
{\rm CH}_\bullet (X)$
is independent of the
choice of $M$ or $i$.
\end{lm}

\proof{
Let $j\colon X\to N$ be
another closed immersion
to a smooth scheme $N$ over $k$.
By considering the product
and the projections,
we may assume that there exists
a smooth morphism
$f\colon M\to N$ compatible
with the immersions $i$ and $j$.
It suffices to show that
$CCi_*{\cal F}$ meets properly
the image of the injection 
$X\times_NT^*N\to
X\times_MT^*M$
and that 
$CCj_*{\cal F}$
is the pull-back of
$CCi_*{\cal F}$
by Lemma \ref{lmincl}.1.
Since the claim is \'etale local on $X$,
we may assume that
there exists a section $s\colon N\to M$
of $f$ such that $i=s\circ j$.
Then the claim follows from
$CCi_*{\cal F}=s_*CCj_*{\cal F}$
Lemma \ref{lmRf}.2.
\qed}

Following the construction
in \cite[A.4]{Gi} in an analytic context,
we make the definition below.

\begin{df}\label{dfccX}
Let $X$ be a scheme of
finite type over a perfect field $k$
and let $\Lambda$ be a
finite field of characteristic invertible in $k$.

{\rm 1.}
We say that $X$ is {\em embeddable}
if there exists a closed immersion 
$i\colon X\to M$
to a smooth scheme $M$ 
over $k$.

{\rm 2. (cf.\ \cite[A.4]{Gi})}
Assume that $X$ is embeddable and
let ${\cal F}$ be a constructible
complex of $\Lambda$-modules
on $X$.
Let $i\colon X\to M$
be a closed immersion 
to a smooth scheme $M$ 
of dimension $n$ over $k$.
We call the class 
\begin{equation}
cc_X{\cal F}
=
\overline{CCi_*{\cal F}}
={\mathbf P}(
CCi_*{\cal F}\oplus {\mathbf A}^1_X)
\in
{\rm CH}_{n}({\mathbf P}(X\times_MT^*M
\oplus {\mathbf A}^1_X))
=
{\rm CH}_\bullet(X)
\label{eqccXF}
\end{equation}
of the closure of the characteristic cycle
$CCi_*{\cal F}$
regarded as an element of
${\rm CH}_\bullet(X)$
the {\em characteristic class}
of ${\cal F}$.

Let $K(X,\Lambda)
=K(X)$ denote the Grothendieck
group of the triangulated category
of constructible complexes
of $\Lambda$-modules on $X$.
We define a linear morphism
\begin{equation}
cc_X\colon
K(X,\Lambda)
\to
{\rm CH}_\bullet(X).
\label{eqch}
\end{equation}
to be that sending the class of 
${\cal F}$ to the characteristic class
$cc_X{\cal F}$.

{\rm 3.}
Assume that $X$ is embeddable and
let $j\colon U\to X$
be an open immersion.
Let 
${\cal F}$ be a locally constant
constructible sheaf of $\Lambda$-modules
on $U$.
Then, we define the 
(total) {\em Swan class}
${\rm Sw}_X{\cal F}
\in {\rm CH}_\bullet (X\sm U)$
to be the class of
the closure of the difference
$CC j_!{\cal F}
-{\rm rank}\ {\cal F}\cdot
CC j_!\Lambda_U$.
\end{df}

If $X$ is quasi-projective,
then $X$ is embeddable.

\begin{cn}\label{cncc}
Let $X$ be an embeddable
scheme of finite type over $k$.
Let $a\colon X\to {\rm Spec}\ k$
be the canonical morphism
and define ${\cal K}_X=Ra^!\Lambda$.

{\rm 1.}
The cohomology class 
$cl(cc_{X,0}{\cal F})
\in H^0(X,{\cal K}_X)$ 
of the dimension $0$-part of
the characteristic class
equals the characteristic class
$C({\cal F})$ defined in 
{\rm \cite{AS}}.

{\rm 2.}
The dimension $0$-part
${\rm Sw}_{X,0}{\cal F}
\in {\rm CH}_0(X\sm U)$
of the total Swan class
equals the Swan class
defined in {\rm \cite{KS}}.
\end{cn}

\begin{lm}\label{lmccd}
{\rm 1.}
The dimension $0$-part
$cc_{X,0}{\cal F}$
is the intersection
product $(CC{\cal F},T^*_XX)_{T^*X}
\in {\rm CH}_0(X)$
with the $0$-section.

{\rm 2.}
Assume that $X$ is irreducible
of dimension $d$.
For a constructible
complex ${\cal F}$
of $\Lambda$-modules
on $X$,
let ${\rm rank}^\circ {\cal F}$
denote the rank of its
restriction on
a dense open subset of $X$
where ${\cal F}$ is locally constant.
Then, the dimension $d$-part
$cc_{X,d}{\cal F}\in
{\rm CH}_d(X)={\mathbf Z}$
is $(-1)^d{\rm rank}^\circ {\cal F}$.
\end{lm}

\proof{
1.
It follows from Lemma \ref{lmEF1}.1.

2.
After shrinking $X$,
we may assume that
${\cal F}$ is locally constant
and $X$ is smooth over $k$.
Then $CC{\cal F}$
is
$(-1)^d{\rm rank}^\circ {\cal F}$-times
the $0$-section
of $T^*X$ and the assertion follows.
\qed}
\medskip

We discuss some functorial
properties of the morphism
$cc_X$.
For $k={\mathbf C}$,
the commutativity of the diagram
(\ref{eqCHf}) below in analytic context
is proved for
projective morphisms in \cite[Theorem 9.4]{Gi}.
This implies \cite[Theorem A.6]{Gi}
a commutative diagram
\begin{equation}
\xymatrix{
K(X,\Lambda)\ar[rr]^{cc_X}
\ar[d]&&
{\rm CH}_\bullet(X)\\
F(X)\ar[rru]_{(-1)^\bullet c_X}&& }
\label{eqcX}
\end{equation}
for the MacPherson-Chern class
$c_X\colon F(X)\to {\rm CH}_\bullet(X)$
\cite{Mac} with multiplied $(-1)^i$ 
on degree $i$ and the canonical morphism
$K(X,\Lambda)\to F(X)$
to the module of ${\mathbf Z}$-valued
constructible functions
defined by the rank of stalks.

However in positive characteristic,
one cannot expect
to have a functoriality
\begin{equation}
\begin{CD}
K(X,\Lambda)
@>{cc_X}>>
{\rm CH}_\bullet(X)\\
@V{f_*}VV@VV{f_!}V\\
K(Y,\Lambda)
@>{cc_Y}>>
{\rm CH}_\bullet(Y)
\end{CD}
\label{eqCHf}
\end{equation}
for proper morphism
$f\colon X\to Y$ over $k$
in full generality (cf.\ \cite[Note 87$_1$]{RS})
except for the dimension $0$-part,
as the following counterexample shows.

\begin{ex}\label{exctr}
We identify ${\rm CH}_\bullet({\mathbf P}^n)$
with $\bigoplus_{i=0}^n{\mathbf Z}$.

{\rm 1.}
The characteristic class $cc_{{\mathbf P}^n}
\Lambda=
(-1)^n[T^*_{{\mathbf P}^n}
{\mathbf P}^n]=
(-1)^n
c(\Omega^1_{{\mathbf P}^n})=
((-1)^i\binom {n+1}{i+1})_{\dim i}$
is non-trivial on every degree.
The endomorphism
$f_*$ of ${\rm CH}_\bullet({\mathbf P}^n)$
induced by the Frobenius morphism
$f\colon {\mathbf P}^n\to {\mathbf P}^n$
is the multiplication by
$p^i$ on the dimension $i$-part.
Since $f_*\Lambda=\Lambda$,
the diagram 
{\rm (\ref{eqCHf})}
is {\em not} commutative for
the Frobenius morphism
$f\colon {\mathbf P}^n\to {\mathbf P}^n$
except possibly for the dimension $0$-part.

{\rm 2.}
Let $j\colon 
{\mathbf A}^1
={\rm Spec}\ k[x]
\to {\mathbf P}^1$
be the open immersion
and let $i\colon 
{\rm Spec}\ k
\to {\mathbf P}^1$
be the closed immersion
of the complement.
Let ${\cal L}$ on ${\mathbf A}^1$ be
the locally constant constructible
sheaves of rank $1$ defined by
the Artin-Schreier equation
$t^p-t=x$ and define
${\cal F}=j_!{\cal L}
\oplus \Lambda$
and 
${\cal G}=j_!\Lambda^{\oplus 2}$.

Then, we have
$CC{\cal F}
=
CC{\cal G}
=-2(T^*_{{\mathbf P}^1}{\mathbf P}^1
+T^*_\infty{\mathbf P}^1)$
and hence
$cc_{{\mathbf P}^1}
{\cal F}
=
cc_{{\mathbf P}^1}
{\cal G}
=(-2,2)$.
Since
$cc_{{\rm Spec} k}
i^*{\cal F}=1\neq
cc_{{\rm Spec} k}
i^*{\cal G}=0$,
there exists {\em no}
right vertical arrow that makes the diagram
\begin{equation}
\begin{CD}
K(X,\Lambda)
@>{cc_X}>>
{\rm CH}_\bullet(X)\\
@V{h^*}VV@VV{?}V\\
K(W,\Lambda)
@>{cc_W}>>
{\rm CH}_\bullet(W)
\end{CD}
\label{eqCHh}
\end{equation}
commutative
for the immersion $\infty\to {\mathbf P}^1$.
\end{ex}

We consider the functoriality
with respect to push-forward.
By Lemma \ref{lmccd}.1,
the commutative diagram
(\ref{eqCHf}) for 
proper smooth morphism $f\colon
X\to {\rm Spec}\ k$
means that the index formula
$(CC{\cal F},T^*_XX)=\chi(X_{\bar k},{\cal F})$
holds for
constructible complexes ${\cal F}$ on $X$.

\begin{lm}\label{lmclim}
The diagram {\rm (\ref{eqCHf})}
is commutative if $f\colon X\to Y$
is a closed immersion
of embeddable schemes
of finite type over $k$.
\end{lm}

\proof{
Let $i\colon Y\to M$
be a closed immersion to
a smooth scheme over $k$.
Then, since both $cc_X{\cal F}$
and $cc_Yf_*{\cal F}$
are defined by
$CC(if)_*{\cal F}$,
the assertion follows.
\qed}
\medskip

Let $f\colon X\to Y$
be a proper morphism
of embeddable schemes
of finite type over $k$.
As in the proof of Corollary
\ref{corcnsm},
we obtain a commutative diagram
\begin{equation}
\begin{CD}
X@>i>> M\\
@VfVV @VV{g}V\\
Y@>j>> N.
\end{CD}
\label{eqXYMN}
\end{equation} where
$M$ and $N$ are smooth,
the right vertical arrow is smooth
and the horizontal arrows are
closed immersions.
We consider the morphisms
\begin{equation}
\begin{CD}
Y\times_NT^*N
@<{\tilde f}<<
X\times_NT^*N
@>h>>
X\times_NT^*M
\end{CD}
\label{eqfh}
\end{equation}
where the left arrow is induced
by $f$ and the right arrow is
the canonical injection.

Let $C
\subset X\times_MT^*M$
be a closed conical subset.
Assume that every irreducible 
component of $M$ and of $C=\bigcup_aC_a$
is of dimension $m$.
Assume also that every irreducible 
component of $N$ 
is of dimension $\leqq n$
and that every irreducible 
component of
the closed conical subset
$f_\circ C=\tilde f(h^{-1}C)
\subset Y\times_NT^*N$
is of dimension $\leqq n$.
Then for a linear combination
$A=\sum_am_aC_a$ of
irreducible components
of $C$,
the push-forward $f_!A=
\sum_a
m_a\tilde f_*h^!C_a$
is defined as a linear combination of
irreducible components of dimension $n$
of $f_\circ C$.

\begin{lm}\label{lmcnpr}
Let $f\colon X\to Y$
be a proper morphism
of embeddable schemes
of finite type over $k$
and let {\rm (\ref{eqXYMN})}
be a commutative diagram
of schemes over $k$.
Assume that
$M$ and $N$ are smooth
of dimension $m$ and $n$
respectively,
that the right vertical arrow is smooth
and that the horizontal arrows are
closed immersions.

Let ${\cal F}$
be a constructible complex of
$\Lambda$-modules on $X$ 
and $C=SS{\cal F}
\subset X\times_MT^*M$ 
be the singular support of
${\cal F}$.
Assume that 
$f_\circ C=\tilde f(h^{-1}C)
\subset Y\times_NT^*N$
is of dimension $\leqq n$
and that we have
$CCR(gi)_*{\cal F}=f_!CCi_*{\cal F}$.
Then, we have
\begin{equation}
cc_YRf_*{\cal F}
=f_*cc_X{\cal F}.
\label{eqprcc}
\end{equation}
\end{lm}

\proof{
We consider the morphisms
\begin{equation}
\begin{CD}
{\mathbf P}(Y\times_NT^*N
\oplus {\mathbf A}^1_Y)
@<\bar f<<
{\mathbf P}(X\times_NT^*N
\oplus {\mathbf A}^1_X)
@>\bar h>>
{\mathbf P}(X\times_MT^*M
\oplus {\mathbf A}^1_X)
\end{CD}
\label{eqfhb}
\end{equation}
extending (\ref{eqfh})
and let $\bar f_!\colon
{\rm CH}_m({\mathbf P}(X\times_MT^*M
\oplus {\mathbf A}^1_X))
\to
{\rm CH}_n({\mathbf P}(Y\times_NT^*N
\oplus {\mathbf A}^1_Y))$
 denote the composition 
$\bar f_*\bar h^!$.
By Lemma \ref{lmincl}.2,
the diagram\begin{equation}
\begin{CD}
{\rm CH}_n({\mathbf P}(X\times_MT^*M
\oplus {\mathbf A}^1_X))
@>>>
{\rm CH}_\bullet(X)
\\
@V{\bar f_!}VV@VV{f_*}V\\
{\rm CH}_m({\mathbf P}(Y\times_NT^*N
\oplus {\mathbf A}^1_Y))
@>>>
{\rm CH}_\bullet(Y)
\end{CD}
\label{eqCCpr}
\end{equation}
is commutative.

Let $A=CC{\cal F}$.
By the assumption
$\dim f_\circ C\leqq n$,
a cycle $\overline{f_!A}$ is 
defined as that
supported on the closure of
$f_\circ C$ in
${\mathbf P}(Y\times_NT^*N
\oplus {\mathbf A}^1_Y)$
by taking the closure of
the cycle $f_!A$.
Further, we have
$\overline{f_!A}=
\bar f_!\bar A$.
Thus, by the assumption
$CCj_*Rf_*{\cal F}=f_!CCi_*{\cal F}$
and (\ref{eqCCpr}),
we obtain
$\overline{CCj_*Rf_*{\cal F}}
=\bar f_!\overline{CCi_*{\cal F}}
=f_*\overline{CCi_*{\cal F}}$.
Hence, the assertion follows.
\qed}
\medskip

\section{Pull-back of characteristic cycle
and the index formula}\label{spb}

We prove that the construction
of characteristic cycles
is compatible with the pull-back
by properly transversal morphisms
in Section \ref{sspb}.
We will derive from this
an index formula for the Euler number
at the end of Section \ref{ssREP}.

In this section,
$k$ denotes a field
of characteristic $\geqq0$.
We assume that
irreducible components of
a smooth scheme over $k$
have the same dimension unless otherwise
stated.
We also assume that
irreducible components of
a closed conical subset 
of the cotangent of
a smooth scheme over $k$
have the same dimension
as the base scheme.

\subsection{Pull-back of characteristic cycle}\label{sspb}

In this subsection,
we assume that
$X$ and $W$ are
smooth schemes over a field $k$.
We assume that
every irreducible component
of $X$ (resp.\ of $W$) is of dimension 
$n$ (resp.\ $m$) and that
every irreducible component of
a closed conical subset
$C\subset T^*X$ is of dimension $n$.
We assume that
a constructible complex
${\cal F}$ of $\Lambda$-modules
on $X$ is of finite tor-dimension.

\begin{df}\label{dfCt}
Let $X$ and $W$ be 
smooth schemes over
a field $k$ and let
$C\subset T^*X$ be a closed conical
subset.
Assume that every irreducible
component of $X$ and
every irreducible
component of $C$
are of dimension $n$
and that
every irreducible
component of $W$ is of dimension $m$.

{\rm 1.}
We say that a $C$-transversal
morphism $h\colon W\to X$ over $k$
is {\em properly $C$-transversal}
if every irreducible component of
$h^*C=W\times_XC$
is of dimension $m$.

{\rm 2.}
Let $h\colon W\to X$ be 
a properly $C$-transversal morphism
and let 
\begin{equation}
\begin{CD}
T^*W@<<< W\times_XT^*X
@>>> T^*X
\end{CD}
\label{eqh!A}
\end{equation}
be the canonical morphisms.
Then, for a linear combination
$A=\sum_am_a[C_a]$
of irreducible components of
$C=\bigcup_aC_a$,
we define $h^!A$ to be 
$(-1)^{n-m}$-{\em times} the push-forward
by the first arrow
$W\times_XT^*X\to T^*W$
in {\rm (\ref{eqh!A})}
of the pull-back of $A$ by the second arrow
$W\times_XT^*X\to T^{*}X$
in the sense of intersection theory.
\end{df}

\begin{lm}
\label{lmdimP}
Let $h\colon W\to X$ be a morphism
of smooth schemes over $k$
and $C\subset T^*X$ be a closed
conical subset.
Assume that every irreducible component
of $X$ is of dimension $n$ and
that every irreducible component
of $W$ is of dimension $m=n-c$.
Let $\dim_{h(W)}C$ denote
the minimum of $\dim C\cap T^*U$
where $U$ runs through
open neighborhoods of the image $h(W)$.

{\rm 1.}
For $h^*C=W\times_XC$,
we have
$\dim h^*C\geqq \dim C_{h(W)}-c$.

{\rm 2.}
Let $g\colon V\to W$ be a morphism
of smooth schemes over $k$.
Assume that every irreducible component
of $C$ is of dimension $n$ and
that every irreducible component
of $V$ is of dimension $l=m-c'$.
Then, 
the following conditions are equivalent:

{\rm (1)}
$h$ is properly $C$-transversal
on a neighborhood of
$g(V)\subset W$ and 
$g\colon V\to W$ is properly 
$h^{\circ}C$-transversal.

{\rm (2)}
The composition $h\circ g\colon V\to X$
is properly $C$-transversal.
\end{lm}

\proof{
1.
If $h$ is smooth, we have
$\dim h^*C= \dim_{h(W)} C-c$.
Hence, it suffices to consider the
case where $h$ is a regular immersion
of codimension $c$.
Then it follows from
\cite[Chap.\ 0 Proposition (16.3.1)]{EGA4}
%(cf.\ Proposition \ref{prpr}.2).

2.
By Lemma \ref{lmCh}.3,
it suffices to verify the conditions
on the dimension.
The implication (1)$\Rightarrow$(2)
is clear. We show
(2)$\Rightarrow$(1).
By 1, we have
$\dim g^*(h^*C)\geqq \dim_{g(V)} h^*C-c'
\geqq n-c-c'$.
Hence (2) implies
the equalities
$\dim_{g(V)} h^*C= n-c$ and
$\dim g^*(h^*C)=\dim_{g(V)} h^*C-c'$.
Thus we have
(2)$\Rightarrow$(1).
\qed}
\medskip

If $h\colon W\to X$ is 
properly $C$-transversal,
then every irreducible component of
$h^{\circ}C$
is of dimension $\dim W$
by Lemma \ref{lmnc}.
A smooth morphism
$h\colon W\to X$ is properly $C$-transversal
by Lemma \ref{lmCh}.1.

\begin{ex}
{\rm (\cite[Example 2.18]{nonlog})}
Assume that $k$ is a perfect
field of characteristic $p>2$.
Let $X={\mathbf A}^2
={\rm Spec}\ k[x,y]
\supset
U={\mathbf G}_m\times
{\mathbf A}^1
={\rm Spec}\ k[x^{\pm1},y]$.
Let ${\cal G}$ be a locally free sheaf
of rank $1$ on $U$
defined by the Artin-Schreier equation
$t^p-t=y/x^p$.
Then, the singular support
$C=SSj_!{\cal G}$
for the open immersion $j\colon
U\to X$ equals the union
of the $0$-section
$T^*_XX$ with the line bundle
$\langle dy\rangle_D$
on the $y$-axis $D=X\sm U$
spanned by the section $dy$.
Hence, the immersion
$D\to X$ is $C$-transversal
but is {\em not} properly
$C$-transversal.
\end{ex}

%as in Definition \ref{dfiv}.

%\section{Characteristic classes and
%Radon transform}\label{sCR}

%\subsection{Radon transform and pull-back}\label{ssRR}
In the rest of this section,
we assume that $k$ is perfect.

\begin{pr}[{Beilinson}]\label{prh!}
Let ${\mathbf P}={\mathbf P}^n$
be a projective space 
and let ${\mathbf P}^\vee$
be the dual projective space.
Let ${\cal G}$ be
a constructible complex of
$\Lambda$-modules
on ${\mathbf P}^\vee$
and let ${\cal F}$ denote
the naive inverse Radon transform
$R{\bm p}_*{\bm p}^{\vee *}{\cal G}$.
Let $C^\vee
\subset T^*{\mathbf P}^\vee$
be a closed conical subset 
such that every irreducible
component is of dimension $n$
and let
$C={\bm p}_{\circ}
{\bm p}^{\vee \circ}C^\vee
\subset T^*{\mathbf P}$.
Assume that ${\cal G}$ is
micro-supported on $C^\vee\subset 
T^*{\mathbf P}^\vee$.

{\rm 1.}
We have
\begin{equation}
{\mathbf P}(CC {\cal F})=
{\mathbf P}({\bm p}_!CC {\bm p}^{\vee*}{\cal G})
=
{\mathbf P}({\bm p}_!{\bm p}^{\vee!}CC {\cal G}).
\label{eqp!}
\end{equation}

{\rm 2.}
Let $X$ be a smooth subscheme
of ${\mathbf P}$
and assume that the immersion
$h\colon X\to {\mathbf P}$
is {\em properly $C$-transversal}.
Then, we have
\begin{equation}
CC h^*{\cal F}=
h^!CC{\cal F}.
\label{eqh!}
\end{equation}
\end{pr}

\proof{
First, we prove
\begin{equation}
{\mathbf P}(CC Rp_*p^{\vee *}{\cal G})=
{\mathbf P}(p_!CC p^{\vee*}{\cal G})
\label{eqph!}
\end{equation}
for properly $C$-transversal
immersion
$h\colon X\to {\mathbf P}$.
We may assume that $k$ is algebraically closed.
Both $CC Rp_*p^{\vee *}{\cal G}$
and $p_!CC p^{\vee*}{\cal G}$
are supported on 
$p_{\circ}p^{\vee \circ}C^\vee
=h^{\circ}C$ by Corollary \ref{corhCt}.2.
By Lemma \ref{lmh!}
and the assumption that
$h$ is properly $C$-transversal,
every irreducible component
of $h^\circ C$ is of dimension $X$.
Hence it suffices to show the equality
\begin{equation}
(CC Rp_*p^{\vee *}{\cal G},df)_u
=(p_!CC p^{\vee*}{\cal G},df)_u
\label{eqMilR}
\end{equation}
for smooth
morphisms $f\colon U\to {\mathbf A}^1$
defined on open subschemes $U\subset X$
with at most an isolated $C$-characteristic point $u$.

By the Milnor formula,
the left hand side of (\ref{eqMilR})
equals
$-\dim{\rm tot}\
\phi_u(Rp_*p^{\vee *}{\cal G},f)$.
By Corollary \ref{corfp},
there exist at most finitely many
isolated $p^{\vee \circ}C$-characteristic points of
$f p\colon U\times_{\mathbf P}Q
\to {\mathbf A}^1$.
%with respect to .
Hence, 
the right hand side of (\ref{eqMilR})
equals
$\sum_v
(CC p^{\vee*}{\cal G},d(f p))_v$
where $v$ runs through
isolated $p^{\vee \circ}C$-characteristic points of
$fp$.
%with respect to $p^{\vee \circ}C$.
Further by the Milnor formula,
this equals
$-\sum_v\dim{\rm tot}\
\phi_v(p^{\vee *}{\cal G},fp)$.
Thus, the equality (\ref{eqMilR})
follows from the isomorphism
$$\phi_u(Rp_*p^{\vee *}{\cal G},f)
\to
R\Gamma(Q\times_Xu,
\phi(p^{\vee *}{\cal G},fp))
\to
\bigoplus_v
\phi_v(p^{\vee *}{\cal G},fp).$$

1.
For the first equality,
it suffices to take
$X={\mathbf P}$ in (\ref{eqph!}).
The second follows from
Proposition \ref{prsm*}.

2.
By proper base change theorem,
we have an isomorphism
$h^*{\cal F}
\to Rp_*p^{\vee *}{\cal G}$.
Hence by (\ref{eqp!})
and (\ref{eqph!}),
we have
${\mathbf P}(CCh^*{\cal F})
={\mathbf P}(h^!CC{\cal F})$.
By the assumption
that the immersion
$h$ is properly $C$-transversal,
the coefficients
of the $0$-section
$T^*_XX$ in $CCh^*{\cal F}$
and in $h^!CC{\cal F}$
are both equal to $(-1)^{\dim X}
{\rm rank}\ {\cal F}$.
Thus the assertion follows.
\qed}
\medskip

For a linear combination
$A=\sum_am_aC_a$
of irreducible closed  conical subsets
of dimension $n$,
we define the Legendre transform
$LA=(-1)^{n-1}{\bm p}_!{\bm p}^!A$.
This is also a linear combination
of irreducible closed  conical subsets
of dimension $n$
by Lemma \ref{lmh!}.
Since the definition of ${\bm p}^!A$
involves the sign $(-1)^{n-1}$,
that of the Legendre transform
does not involve sign and
we have
${\mathbf P}(L(A))=
\sum_{a;C_a\not\subset 
T^*_{\mathbf P}{\mathbf P}}
m_a{\mathbf P}(C_a^\vee)$.

\begin{cor}[Beilinson]\label{corCCR}
Let ${\cal F}$
be a constructible complex
of $\Lambda$-modules on ${\mathbf P}$.
Then, for the Radon transform
$R{\cal F}$, we have
\begin{equation}
{\mathbf P}(CCR{\cal F})=
{\mathbf P}(LCC{\cal F}).
\label{eqCCR}
\end{equation}
\end{cor}

We will remove ${\mathbf P}$
in (\ref{eqCCR})
in Corollary \ref{corLR}.

\proof{
By Proposition \ref{prh!}.1,
we have 
$${\mathbf P}(CCR{\mathcal F})
=
{\mathbf P}((-1)^{n-1}{\bm p}^\vee_!
CC{\bm p}^*{\mathcal F})
=
{\mathbf P}((-1)^{n-1}{\bm p}^\vee_!
{\bm p}^!
CC{\mathcal F})
=
{\mathbf P}
(LCC{\mathcal F}).$$
\qed}
%\medskip

\begin{thm}\label{thmi*}
Let $X$ and $W$ be 
smooth schemes over a perfect field $k$ and let
$C\subset T^*X$ be a closed conical
subset.
Assume that every irreducible
component of $X$ and
every irreducible
component of $C$
are of dimension $n$
and that
every irreducible
component of $W$ is of dimension $m$.

Let ${\cal F}$ be
a constructible complex of
$\Lambda$-modules on $X$
of finite tor-dimension
micro-supported on
$C\subset T^*X$
and
let $h\colon W\to X$ be a 
properly $C$-transversal morphism.
Then, we have 
\begin{equation}
CC h^*{\cal F}=
h^!CC {\cal F}.
\label{eqii}
\end{equation}
\end{thm}

\proof[Proof {\rm (Beilinson)}]{
Let $h\colon W\to X$
be a properly $C$-transversal morphism.
Since $h$ is decomposed
as the composition
of the graph $W\to W\times X$
and the projection $W\times X\to X$
and since it is proved for
smooth morphisms in Proposition \ref{prsm*},
it is sufficient to show
the case where $h$ is an immersion.

First, we show the case
where $X$ is the projective space 
${\mathbf P}={\mathbf P}^n$.
The case where ${\cal F}=R{\bm p}_*
{\bm p}^{\vee*}{\cal G}$
is the naive inverse Radon transform
has been proved in Proposition
\ref{prh!}.2.
Let ${\cal F}$ be a constructible
complex on ${\mathbf P}^n$.
Since ${\cal F}$ is isomorphic
to $R^\vee R{\cal F}$ upto
constant sheaf
and the assertion is clear for
the constant sheaf,
it follows in the case
$h$ is an immersion  
to ${\mathbf P}$.

We show the general case.
Since the assertion is local,
we may assume that $X$ is
affine and take an immersion
$i\colon X\to {\mathbf P}$.
Further, we may assume that
there is a smooth subscheme
$V\subset {\mathbf P}$
such that $X\cap V=W$
and that the intersection
is transversal.
Then, the immersion $\tilde h\colon
V\to {\mathbf P}$
is $i_{\circ}C$-transversal
on a neighborhood of $W$.
Hence, it follows from the case
where $h$ is an immersion  
to ${\mathbf P}$.
\qed}
\smallskip

We study the compatibility
of characteristic classes with pull-back.
Let $F\to E$ be
an injection of vector bundles
over a scheme $W$ of finite type over $k$
and let $p\colon E\to E/F$
be the canonical surjection
Let $C\subset E$
be a closed conical subset
such that
the intersection
$C\cap F$
is a subset of the $0$-section.
Then, for a linear combination
$A=\sum_am_aC_a$ of irreducible
components of $C=\bigcup_aC_a$,
the intersection theory defines
a cycle $p^!p_*A$ supported
on $p^{-1}(p(C))=C+F
\subset E$.

\begin{df}\label{dfprim}
Let $X$ be an embeddable scheme of
finite type over $k$
and let $i\colon X\to M$
be a closed immersion to
a smooth scheme over $k$.
Let $C\subset X\times_MT^*M$
be a closed conical subset.
Let $h\colon W\to X$
be a regular immersion
of codimension $c$.

{\rm 1.}
We say that $h$ is 
{\em properly $C$-transversal}
if the following conditions are
satisfied:
The canonical morphism
$T^*_WX\to W\times_MT^*M$
is an injection of vector bundles on $W$.
For every irreducible
components $C_a$ of $C=\bigcup_aC_a$,
the pull-back $W\times_XC_a\subset C_a$
is of codimension $c$.
The intersection 
$h^*C\cap T^*_WX$
for $h^*C=W\times_XC
\subset W\times_MT^*M$
is a subset of the $0$-section.

{\rm 2.}
Assume that $h\colon W\to X$
is properly $C$-transversal.
We regard the pull-back
$h^*C=
W\times_XC$ as a subset
of $W\times_MT^*M$
and the conormal bundle
$T^*_WX$ as a sub vector bundle
of $W\times_MT^*M$.
Then, we define a closed conical
subset 
$h^!C\subset W\times_MT^*M$
to be the sum
$h^*C+T^*_WX$.

For a linear combination
$A=\sum_am_aC_a$ of
irreducible components of
$C=\bigcup_aC_a$,
let $h^*C_a$ denote the
pull-back of $C_a$ by
the regular immersion
$W\times_MT^*M\to 
X\times_MT^*M$ in the sense of
intersection theory
and let $p\colon W\times_MT^*M
\to (W\times_MT^*M)/T^*_WX$
be the canonical surjection.
Then,
we define
$h^!A=(-1)^c\sum_am_a
p^!p_*h^*C_a$.
\end{df}

\begin{pr}\label{prprim}
Let $X$ be an embeddable scheme of
finite type over $k$
and let $i\colon X\to M$
be a closed immersion to
a smooth scheme over $k$.
Let ${\cal F}$
be a constructible complex of
$\Lambda$-modules on $X$ 
and $C=SS{\cal F}
\subset X\times_MT^*M$ 
be the singular support of
${\cal F}$.

Let $h\colon W\to X$
be a {\em properly
$C$-transversal} closed regular immersion
of codimension $c$.
Then, we have
\begin{equation}
CC(i\circ h)_*h^*{\cal F}
=h^!CCi_*{\cal F},
\label{eqprim}
\end{equation}
\begin{equation}
cc_Wh^*{\cal F}
=
(-1)^c
c(T^*_WX)^{-1}\cap h^!cc_X{\cal F}.
\label{eqprim2}
\end{equation}
\end{pr}

\proof{
We show (\ref{eqprim}).
Since the assertion is local on
$W$, we may assume that
there exists a cartesian diagram
$$\begin{CD}
W@>h>> X\\
@VVV@VViV\\
N@>j>>M
\end{CD}$$
where the lower
horizontal arrow
$j\colon N\to M$ 
is a regular immersion
of codimension $c$
of smooth schemes over $k$.
Further, we may assume that
 the immersion $j\colon N\to M$ 
is properly $C$-transversal.
Then, by Theorem \ref{thmi*},
we have $CCj^*i_*{\cal F}
=j^!CCi_*{\cal F}$.
Since
$(i\circ h)_*h^*{\cal F}
=
j_*j^*i_*{\cal F}$, 
we obtain
$CC(i\circ h)_*h^*{\cal F}
=
CCj_*j^*i_*{\cal F}
=
j_*CCj^*i_*{\cal F}
=
j_*j^!CCi_*{\cal F}
=
h^!CCi_*{\cal F}$
by Lemma \ref{lmRf}.2.

The equality (\ref{eqprim2})
follows from (\ref{eqprim}) and
Lemma \ref{lmtheta}
since the definition of $h^!CCi^*{\cal F}$
involves the sign $(-1)^c$.
\qed}

\begin{cor}\label{corcnsm}
Let $f\colon X\to Y$
be a smooth morphism of 
relative dimension $d$ of
embeddable schemes 
of finite type over $k$
and let $T^*X_{/Y}$ 
and $c(T^*X_{/Y})$ denote
the relative cotangent bundle
and its total Chern class.
Then,
the diagram
\begin{equation}
\begin{CD}
K(Y,\Lambda)
@>{cc_Y}>>
{\rm CH}_\bullet(Y)\\
@V{f^*}VV@VV{(-1)^dc(T^*X_{/Y})\cap f^!}V\\
K(X,\Lambda)
@>{cc_X}>>
{\rm CH}_\bullet(X)
\end{CD}
\label{eqCHhf}
\end{equation}
is commutative.
In particular, for $Y={\rm Spec}\ k$
and for a geometrically constant
sheaf ${\cal F}$ on $X$ of dimension $n$,
we have
\begin{equation}
cc_X{\cal F}
=
(-1)^n
{\rm rank}\ {\cal F}\cdot
c(T^*X).
\label{eqY=k}
\end{equation}
\end{cor}

\proof{
Let $i\colon X\to M$
and $j\colon Y\to N$
be closed immersions
to smooth schemes over $k$.
After replacing $i\colon X\to M$
by $(i,fj)\colon X\to M\times N$,
we obtain a commutative diagram
\begin{equation}
\xymatrix{
X\ar[rdd]_f\ar[rd]^{\!\!\! h}\ar[rrd]^i&&\\
&Y'%=Y\times_NM
\ar[r]_{ j'}\ar[d]^{g'}&M\ar[d]^g\\
&Y\ar[r]^j&N
}
\label{eqXYMNg}
\end{equation}
where $g\colon M\to N$ is smooth
and the square is cartesian.
Since $g$ is smooth
and $j'_*g^{\prime *}{\cal F}=
g^*j_*{\cal F}$,
we have
$CCj'_*g^{\prime *}{\cal F}=
g^!CCj_*{\cal F}$
by Theorem \ref{thmi*}.
Let $C=SSj_*{\cal F}
\subset Y\times_NT^*N$
be the singular support
and let $C'=g^oC
\subset Y'\times_MT^*M$.
Then, since $g\colon M\to N$ is smooth,
the closed immersion
$h\colon X\to Y'$
is properly $C'$-transversal.
Hence by Proposition \ref{prprim},
we have
$CCi_*h^*g^{\prime *}{\cal F}=
h^!CCj'_*g^{\prime *}{\cal F}=
h^!g^!CCj_*{\cal F}$.
Namely, we have
$CCi_*f^*{\cal F}=
f^!CCj_*{\cal F}$.
Since the definition of
the right hand side
involves the sign $(-1)^d$,
we obtain
\begin{equation}
cc_Xf^*{\cal F}=
c(T^*X_{/Y})\cap
(-1)^df^!cc_Y{\cal F}
\label{eqccTXY}
\end{equation}
by Lemma \ref{lmincl}.2.
Thus the assertion follows.
\qed}
%\medskip

\subsection{Radon transform and 
the index formula}\label{ssREP}

Assume $X={\mathbf P}^n$.
We identify 
$Q={\mathbf P}(T^*{\mathbf P})$
and let
${\bm p}\colon Q\to {\mathbf P}$ and
${\bm p}^\vee\colon Q\to {\mathbf P}^\vee$
be the projections.
The Radon transforms
$R=R{\bm p}^\vee_!
{\bm p}^*[n-1]$ and 
$R^\vee=R{\bm p}_!
{\bm p}^{\vee*}[n-1](n-1)$
define morphisms
\begin{equation}
R\colon
K({\mathbf P},\Lambda)
\to
K({\mathbf P}^\vee,\Lambda),
\quad
R^\vee\colon
K({\mathbf P}^\vee,\Lambda)
\to
K({\mathbf P},\Lambda).
\label{dfKR}
\end{equation}
Define also morphisms $\chi\colon 
K({\mathbf P},\Lambda)
\to {\mathbf Z}$
and $\chi\colon 
K({\mathbf P},\Lambda)
\to {\mathbf Z}$
by $\chi{\cal F}=
\chi({\mathbf P}_{\bar k},{\cal F})$
and
by $\chi{\cal G}=
\chi({\mathbf P}^\vee_{\bar k},{\cal G})$.

\begin{lm}\label{lmPn}
Let $n\geqq 1$ be an integer
and ${\mathbf P}={\mathbf P}^n$.

{\rm 1.}
For $a=0,\ldots,n$,
let ${\mathbf P}^a
\subset {\mathbf P}
={\mathbf P}^n$
denote a linear subspace
of dimension $a$.
Then, the images
$cc_{{\mathbf P}^n}(
\Lambda_{{\mathbf P}^a}[a])$
for $a=0,\ldots,n$
do not depend on
the choice of linear subspaces
and 
form a basis of
a free ${\mathbf Z}$-module
${\rm CH}_\bullet ({\mathbf P})$.

{\rm 2.}
The diagram
\begin{equation}
\begin{CD}
K({\mathbf P},\Lambda)
@>{\chi}>>
{\mathbf Z}
\\
@VRVV@VV{(-1)^{n-1}n\times}V\\
K({\mathbf P}^\vee,\Lambda)
@>{\chi}>>
{\mathbf Z}
\end{CD}
\label{eqRchi}
\end{equation}
and that with $R$ replaced
by $R^\vee$ and with
${\mathbf P}$ and
${\mathbf P}^\vee$ switched
are commutative.

{\rm 3.}
Assume $k$ is algebraically closed.
The composition $R^\vee R\colon
K({\mathbf P},\Lambda)\to
K({\mathbf P},\Lambda)$
maps ${\cal F}$ to
${\cal F}+(n-1)\chi{\cal F}\Lambda$.
If $n\neq 1$
and if $f\colon
K({\mathbf P},\Lambda)
\to {\mathbf Z}$
is a morphism satisfying
$fR^\vee R=n^2f$,
then we have
$f=f([\Lambda_{{\mathbf P}^0}])\cdot
\chi$.
\end{lm}

\proof{1.
It follows from
$CC
\Lambda_{{\mathbf P}^a}[a]=
T^*_{{\mathbf P}^a}{\mathbf P}^n$.

2.
For constructible complexes ${\cal F}$
on ${\mathbf P}$,
we have
$\chi R{\cal F}
=
(-1)^{n-1}\chi(Q_{\bar k},{\bm p}^*{\cal F})
=
(-1)^{n-1}\chi 
R{\bm p}_*{\bm p}^*{\cal F}
=
(-1)^{n-1}\chi
R{\bm p}_*{\bm p}^*{\Lambda}\otimes{\cal F}$
by the projection formula.
Hence the assertion follows from
$R^q{\bm p}_*{\bm p}^*{\Lambda}
=\Lambda(-q/2)$
for $0\leqq q\leqq 2(n-1)$ even
and $=0$ for otherwise.

3.
Since $R^\vee R{\cal F}$
is isomorphic to ${\cal F}$
up to geometrically constant sheaves,
we have
$R^\vee R{\cal F}-{\cal F}
={\rm rank}^\circ
(R^\vee R{\cal F}-{\cal F})\cdot
\Lambda$.
By taking $\chi$ and applying 2,
we obtain
$(n^2-1)$
$\chi{\cal F}
=
{\rm rank}^\circ
(R^\vee R{\cal F}-{\cal F})\cdot
\chi\Lambda$.
Since
$\chi\Lambda=n+1$,
we obtain
${\rm rank}^\circ
(R^\vee R{\cal F}-{\cal F})$
$=
(n-1)\chi{\cal F}$
and the first assertion follows.

If $f$ satisfies the condition,
similarly we obtain
$(n^2-1)
f({\cal F})
=
(n-1)\chi{\cal F}\cdot
f(\Lambda)$.
Thus, $f$ is a constant multiple of
$\chi$.
Since $\chi\Lambda_{{\mathbf P}^0}=1$,
the constant is given by
$f(\Lambda_{{\mathbf P}^0})$.
\qed}
\smallskip

We define the Legendre transform
\begin{equation}
L\colon
{\rm CH}_\bullet ({\mathbf P})
\to
{\rm CH}_\bullet ({\mathbf P}^\vee)
\label{eqbarL}
\end{equation}
by $L={\bm p}^\vee_*{\bm p}^*$
and $L^\vee={\bm p}_*{\bm p}^{\vee*}$
for the projections
${\bm p}\colon Q\to {\mathbf P}$ and
${\bm p}^\vee\colon Q\to {\mathbf P}^\vee$.

\begin{pr}[Beilinson]\label{prRR}
Let $n\geqq 1$ be an integer
and ${\mathbf P}={\mathbf P}^n$.

{\rm 1.}
The diagram
\begin{equation}
\begin{CD}
K({\mathbf P},\Lambda)
@>{cc_{\mathbf P}}>>
{\rm CH}_\bullet ({\mathbf P})
\\
@VRVV@VVLV\\
K({\mathbf P}^\vee,\Lambda)
@>{cc_{{\mathbf P}^\vee}}>>
{\rm CH}_\bullet ({\mathbf P}^\vee)
\end{CD}
\label{eqR}
\end{equation}
and that with $R$ and $L$ replaced
by $R^\vee$ and $L^\vee$ and with
${\mathbf P}$ and
${\mathbf P}^\vee$ switched
are commutative.

{\rm 2.}
The diagram
\begin{equation}
\xymatrix{
K({\mathbf P},\Lambda)\ar[rrd]_{\chi}
\ar[rr]^{\!\!\!\!\!\!\!
cc_{\mathbf P}}&&
{\rm CH}_\bullet ({\mathbf P})
\ar[d]^{\deg}\\
&&{\mathbf Z}}
\label{eqchP}
\end{equation}
is commutative.
\end{pr}

\proof{
We may assume that $k$ is algebraically closed.
We prove the assertions by 
induction on $n$.
If $n=1$, the projections
${\bm p}\colon Q\to {\mathbf P}$
and
${\bm p}^\vee\colon Q\to {\mathbf P}^\vee$
are isomorphisms
and the assertion 1 is obvious.
Since $\deg cc_{\mathbf P}{\cal F}
=
(CC{\cal F},T^*_{\mathbf P}{\mathbf P})_
{T^*{\mathbf P}}$,
the assertion 2 for $n=1$ is nothing but
the Grothendieck-Ogg-Shafarevich
formula for ${\mathbf P}={\mathbf P}^1$.

We show that the assertion 2 for $n-1\geqq 1$
implies the assertion 1 for $n$.
We show the commutativity of
the diagram (\ref{eqR})
by using the direct sum decomposition
\begin{equation}
\begin{CD}
{\rm CH}_\bullet({\mathbf P}^\vee)
=
{\rm CH}_n({\mathbf P}(T^*{\mathbf P}^\vee
\oplus {\mathbf A}^1_{{\mathbf P}^\vee}))
@>>>
{\rm CH}_{n-1}({\mathbf P}(T^*{\mathbf P}^\vee))
\oplus
{\rm CH}_n({\mathbf P}^\vee)\\
@.=
{\rm CH}_{n-1}(Q)
\oplus {\mathbf Z}
\end{CD}
\label{eqCHQZ}
\end{equation}
(\ref{eqQZ}).
The compositions with the
first projection $
{\rm CH}_\bullet({\mathbf P}^\vee)
\to {\rm CH}_{n-1}(Q)$
are equal by Corollary \ref{corCCR}.
We show that the compositions with the
second projection ${\rm pr}_2\colon
{\rm CH}_\bullet({\mathbf P}^\vee)
\to {\mathbf Z}$
induced by the projection
${\mathbf P}(T^*{\mathbf P}^\vee
\oplus {\mathbf A}^1_{{\mathbf P}^\vee})
\to 
{\mathbf P}^\vee$
are the same.

Let ${\cal F}$ be a constructible
complex of $\Lambda$-modules
on ${\mathbf P}$
and $C=SS{\cal F}$ 
be the singular support of ${\cal F}$.
Let $H\subset {\mathbf P}$
be a hyperplane
such that the immersion
$h\colon H\to {\mathbf P}$
is properly $C$-transversal
and let $i\colon {\rm Spec}\ k
\to {\mathbf P}^\vee$ 
be the immersion of the $k$-rational
point of ${\mathbf P}^\vee$ 
corresponding to $H$.

By the hypothesis of induction,
we have
$\deg cc_Hh^*{\cal F}=
\chi h^*{\cal F}$.
By Proposition \ref{prprim},
we have
$cc_Hh^*{\cal F}=
-h^!cc_{\mathbf P}{\cal F}$.
Hence by the commutative diagram
\begin{equation}
\begin{CD}
{\rm CH}_\bullet ({\mathbf P})
@>{h^!}>>
{\rm CH}_\bullet (H)
\\
@VLVV@VV{\deg}V\\
{\rm CH}_\bullet ({\mathbf P}^\vee)
@>{i^!}>>{\mathbf Z}
\end{CD}
\label{eqH}
\end{equation}
and by the last assertion in
Lemma \ref{lmCHP}.2,
we obtain
${\rm pr}_2L
cc_{\mathbf P}{\cal F}
=i^!L
cc_{\mathbf P}{\cal F}
=-\deg cc_Hh^*{\cal F}
=-\chi h^*{\cal F}$.
We also have
${\rm pr}_2
cc_{{\mathbf P}^\vee} R{\cal F}
=
(-1)^n
{\rm rank}^\circ R{\cal F}
=
-\chi h^*{\cal F}$
by Lemma \ref{lmccd}
in the notation ${\rm rank}^\circ$
defined there.
Hence the assertion 1 follows.

We show that the assertion 1 for $n\geqq 2$
implies the assertion 2 for $n$.
By the commutative diagrams
(\ref{eqR}),
the endomorphism $R^\vee R$
of $K({\mathbf P},\Lambda)$
preserves the kernel
$K({\mathbf P},\Lambda)^0$  of 
$cc_{\mathbf P}
\colon K({\mathbf P},\Lambda)
\to 
{\rm CH}_\bullet({\mathbf P})$.
We show that $R^\vee R$ acts
on $K({\mathbf P},\Lambda)^0$ 
as the identity.
Since $R^\vee R{\cal F}$
is isomorphic to ${\cal F}$
up to constant sheaf,
we have
$R^\vee R{\cal F}-{\cal F}
=
{\rm rank}^\circ(R^\vee R{\cal F}-{\cal F})
\cdot \Lambda$.
Since 
${\rm rank}^\circ=(-1)^n
{\rm pr}_2cc_{\mathbf P}$
by Lemma \ref{lmccd},
we have
${\rm rank}^\circ {\cal F}=0$
for ${\cal F}\in K({\mathbf P},\Lambda)^0$.
Further by the commutative diagrams
(\ref{eqR}),
we also have
${\rm rank}^\circ R^\vee R{\cal F}=0$
for ${\cal F}\in K({\mathbf P},\Lambda)^0$.
Hence
$R^\vee R$ acts
as the identity on
$K({\mathbf P},\Lambda)^0$.

On the other hand,
by Lemma \ref{lmPn}.3,
we have
$\chi R^\vee R=n^2\chi $.
Since $n^2>1$,
$\chi$ annihilates
$K({\mathbf P},\Lambda)^0=
{\rm Ker}\ cc_{\mathbf P}$
and 
$\chi$ induces a unique morphism
${\rm CH}_\bullet({\mathbf P})
\to {\mathbf Z}$ by Lemma \ref{lmPn}.1.

For the classes of
${\cal F}=\Lambda_{{\mathbf P}^a}[a]$,
we have
$\chi {\cal F}=(CC{\cal F},T^*_{\mathbf P}
{\mathbf P})_{T^*{\mathbf P}}
=\deg cc_{\mathbf P}{\cal F}$.
Thus, the diagram (\ref{eqchP}) is commutative
by Lemma \ref{lmPn}.1.
\qed}

\begin{cor}[Beilinson]\label{corLR}
Let ${\cal F}$
be a constructible complex
of $\Lambda$-modules on ${\mathbf P}$.
Then, for the Radon transform
$R{\cal F}$, we have
\begin{equation}
CCR{\cal F}=
LCC{\cal F}.
\label{eqCCL}
\end{equation}
\end{cor}

\proof{
Except for the coefficient
of the $0$-section, it is proved in
Corollary \ref{corCCR}.
Since the coefficient
of the $0$-section
is given by 
${\rm pr}_2\colon
{\rm CH}_\bullet({\mathbf P}^\vee)
\to {\mathbf Z}$,
it follows from
Proposition \ref{prRR}.1.
\qed}

Corollary \ref{corLR} means that
\cite[Conjecture 2.2]{notes} holds
for $p^\vee\colon Q\to {\mathbf P}^\vee$
and $p^*{\cal F}$
for ${\cal F}$ on ${\mathbf P}$ 
by Lemma \ref{lmh!}.

We state and prove
the index formula
for the Euler-Poincar\'e characteristic.

\begin{thm}\label{thmEP}
Let $X$ be a projective smooth
variety over an algebraically closed
field.
Then, we have
\begin{equation}
\chi(X,{\cal F})
=
(CC {\cal F},
T^{*}_{X}X)_{T^{*}X}.
\label{eqEP}
\end{equation}
\end{thm}

\proof[Proof {\rm (Beilinson)}]{
Since $X$ is assumed projective,
it follows from Lemma \ref{lmclim}
and Proposition \ref{prRR}.2.
\qed}

\subsection{Characteristic cycle
and ramification}\label{ssramc}

We briefly recall the definition 
of the characteristic cycle
\cite[Definition 3.5]{nonlog}
in the strongly non-degenerate case.
We use the notation in Section \ref{ssram}.

Let $X$ be a smooth scheme
of dimension $n$ over a perfect 
field $k$ and
$U=X\sm D$ be
the complement of a divisor $D$
with normal crossings.
Let $j\colon U\to X$ be
the open immersion and
${\cal G}$ be a locally constant
constructible sheaf of free $\Lambda$-modules
on $U$.
Assume that the ramification of
${\cal G}$ along $D$ is strongly non-degenerate.

Assume first that $R=D$.
Then, the locally constant sheaf
${\cal G}$ on $U$ is tamely ramified along $D$
and we define a linear combination by
\begin{equation}
C( j_!{\cal G})
=(-1)^n{\rm rank}\ {\cal G}\cdot
\sum_I[T^*_{D_I}X]
\label{eqtame}
\end{equation}
where $T^*_{D_I}X$
denotes the conormal bundle
of the intersection $D_I=\bigcap_ID_i$
for sets of indices
$I\subset \{1,\ldots,m\}$.

Assume $R=\sum_ir_iD_i>
D=\sum_iD_i$. 
For each irreducible component,
we have a decomposition by 
characters 
$V=\bigoplus_{\chi\colon
{\rm Gr}^{r_i}G_{K_i}\to 
{\mathbf F}_p}
\chi^{\oplus m(\chi)}$.
Further, the characteristic form
of each character $\chi$
appearing in the decomposition
defines a sub line bundle $L_\chi$
of the pull-back $D_\chi\times_XT^*X$
of the cotangent bundle 
to a finite covering $\pi_\chi
\colon D_\chi\to D_i$
by the non-degenerate assumption.
Then we define
\begin{equation}
C(j_!{\cal G})
=(-1)^n\Bigl(
{\rm rank}\ {\cal G}\cdot
[T^*_XX]
+
\sum_i
\sum_{\chi}
\dfrac{r_i\cdot m(\chi)}{[D_\chi\colon D_i]}
\pi_*[L_\chi]\Bigr).
\label{eqwild}
\end{equation}
In the general strongly non-degenerate case,
we define $C(j_!{\cal G})$ by the additivity
and \'etale descent.

\begin{thm}\label{thmram}
Let $X$ be a smooth scheme
of dimension $n$ over a perfect 
field $k$ and
let $j\colon U\to X$ be the
open immersion of
the complement
$U=X\sm D$ of
a divisor $D$ with simple normal crossings.
Let ${\cal G}$
be a locally constant constructible
sheaf of free $\Lambda$-modules
on $U$ such that
the ramification along $D$
is {\em strongly non-degenerate}.
Then we have
\begin{equation}
CC j_!{\cal G}
=
C(j_!{\cal G}),
\label{eqCCj}
\end{equation}
In other words, 
$C(j_!{\cal G})$ defined by
{\rm (\ref{eqtame})},
{\rm (\ref{eqwild})},  the additivity
and by \'etale descent
using
ramification theory
satisfies the Milnor formula
{\rm (\ref{eqMil})}.
\end{thm}

Theorem \ref{thmram} is
proved for dimension $\leqq 1$
in Lemma \ref{lmele}.1 and 2.
Recall that Theorem \ref{thmram} is
proved in dimension $2$
in \cite[Proposition 3.20]{surface}
using a global argument, as in \cite{bp}.
The tamely ramified case
of Theorem \ref{thmram} has been proved
in \cite{Yang} by a different method.
Theorem \ref{thmram} gives
an affirmative answer to
\cite[Conjecture 3.16]{nonlog}.

\proof{
It suffices to show the
equality of the coefficients
for each irreducible component
of $C=S(j_!{\cal F})$ by Proposition \ref{prnd}.
By the additivity of characteristic cycles
and the compatibility with \'etale pull-back,
it suffices to show the
tamely ramified case
and the totally wildly ramified
case separately.

First, we prove the tamely ramified case.
It suffices to determine the
coefficients of $[T^*_{D_I}X]$
by induction on the number $d$ of
elements of $I$.
By Proposition \ref{prsm*},
we may assume $X={\mathbf A}^n$
and $D$ is the complement of
${\mathbf G}_m^n$.
If $n=0$ and $X$ consists of a
single point, we have
$CC {\cal F}=
{\rm rank}\ {\cal F}\cdot [T^*_XX]$.
If $n\geqq 1$,
it follows from Lemma \ref{lmtame}.

We prove the totally wildly ramified case.
If $\dim X=0$, it is proved above.
It follows from (\ref{eqdim1}) if
$\dim X\leqq 1$.
It suffices to compare
the coefficients in
$CC j_!{\cal G}$
and $C(j_!{\cal G})$
assuming $\dim X\geqq 2$.
Since the assertion is
\'etale local, we may 
assume that $C(j_!{\cal G})$
has a unique irreducible
component different from the $0$-section.
By Theorem \ref{thmi*}
and Proposition \cite[Proposition 3.8]{nonlog},
for every properly $C$-transversal immersion $i\colon W\to X$
of a smooth curve,
we have 
$i^!CC j_!{\cal G}=
CC i^* j_!{\cal G}$
and $i^!C(j_!{\cal G})
=C(i^*j_!{\cal G})$
respectively.
Hence it is reduced to the case
$\dim X=1$.
\qed}

\section{${\cal F}$-transversality and singular support}\label{sFt}

We introduce and study a notion of
${\cal F}$-transversality
for a morphism $h\colon W\to X$
with respect to a constructible
complex ${\cal F}$ on $X$
in Section \ref{ssFt}
using a canonical morphism (\ref{eqprpg})
defined and studied in Section \ref{ssFcan}.
We study the relation between
the ${\cal F}$-transverality
and the local acyclicity in Section \ref{ssFla}
and deduce a relation
with the singular support in Section \ref{ssFs}.

\subsection{A canonical morphism}\label{ssFcan}

Let $k$ be a field
and $\Lambda$ be a finite
local ring such that
the residue characteristic
$\ell$ is invertible in $k$.
For a separated morphism $h\colon W\to X$
of schemes of finite type
over $k$,
the functor
$Rh^!\colon D^b_c(X,\Lambda)
\to D^b_c(W,\Lambda)$
is defined as the adjoint
of
$Rh_!\colon D^b_c(W,\Lambda)
\to D^b_c(X,\Lambda)$.
One should be able to define
the functors $Rh_!,Rh^!$
without assuming separatedness
using cohomological descent
but we will not go further in this
direction.

Let ${\cal F}$ and ${\cal G}$ 
be constructible complexes
of $\Lambda$-modules
on $X$ and on $W$
respectively
and assume that
${\cal G}$ is of finite tor-dimension.
Then, we have an isomorphism
\begin{equation}
{\cal F}\otimes^L_\Lambda
Rh_!{\cal G} 
\to Rh_!(h^*{\cal F}\otimes^L_\Lambda{\cal G} )
\label{eqprj}
\end{equation}
of projection formula
\cite[(4.9.1)]{Rapport}.
\if{This is defined as the adjoint
$h^*{\cal F}\otimes^L_\Lambda
h^*Rh_!{\cal G} 
\to h^*{\cal F}\otimes^L_\Lambda{\cal G}$
of the morphism induced
by the adjunction
$h^*Rh_!{\cal G} 
\to {\cal G}$
if $h$ is proper.
It is defined as the inverse of
the isomorphism
$h^*{\cal F}\otimes^L_\Lambda
h^*Rh_!{\cal G} 
\gets h^*{\cal F}\otimes^L_\Lambda{\cal G}$
if $h$ is an open immersion.}\fi

\begin{df}\label{dfAB}
Let $h\colon W\to X$
be a morphism of finite
type of schemes.
Let ${\cal F}$ and ${\cal E}$ be
constructible complexes
of $\Lambda$-modules on $X$
and let ${\cal G}$ be
a constructible complex
of $\Lambda$-modules 
of finite tor-dimension on $W$.
We say that morphisms
\begin{equation}
{\cal F}
\otimes^L
Rh_!{\cal G}
\to {\cal E},\quad
h^*{\cal F}
\otimes^L
{\cal G}
\to Rh^!{\cal E}
\label{eqAB}
\end{equation}
correspond to each other
if the first one
is the composition 
%of the inverse 
of the isomorphism {\rm (\ref{eqprj})}
with the adjoint
$Rh_!(h^*{\cal F}
\otimes^L
{\cal G})
\to
{\cal E}$
of the second one.
\end{df}

Since {\rm (\ref{eqprj})}
is an isomorphism,
the correspondence
{\rm (\ref{eqAB})}
is a one-to-one correspondence.
If $h\colon W\to X$
is a smooth separated morphism
of relative dimension $d$,
the morphism
${\cal F}
\otimes^L
Rh_!\Lambda(d)[2d]
\to {\cal F}$
induced by
the trace mapping
$Rh_!\Lambda(d)[2d]\to \Lambda$
corresponds to the canonical isomorphism
\begin{equation}
h^*{\cal F}
\otimes^L
\Lambda(d)[2d]
\to Rh^!{\cal F}
\label{eqPD}
\end{equation}
of Poincar\'e duality
\cite[Th\'eor\`eme 3.2.5]{DP}.

Recall that a morphism
${\cal F}\otimes^L_\Lambda
Rh_!{\cal G} 
\to {\cal E}$
corresponds to a morphism
${\cal F}\to
R{\cal H}om(
Rh_!{\cal G},{\cal E})$
bijectively.
Similarly
a morphism
$h^*{\cal F}\otimes^L_\Lambda
{\cal G} 
\to Rh^!{\cal E}$
corresponds to a morphism
$h^*{\cal F}\to
R{\cal H}om({\cal G},Rh^!{\cal E})$
bijectively.

\begin{lm}\label{lmAB}
Let $h\colon W\to X$
be a morphism of finite
type of schemes.
Let ${\cal F}$ and ${\cal E}$ be
constructible complexes
of $\Lambda$-modules on $X$
and let ${\cal G}$ be
a constructible complex
of $\Lambda$-modules 
of finite tor-dimension on $W$.
Let
\begin{equation}
{\cal F}
\otimes^L
Rh_!{\cal G}
\to {\cal E},\quad
h^*{\cal F}
\otimes^L
{\cal G}
\to Rh^!{\cal E}
\label{eqAB0}
\end{equation}
be morphisms corresponding
to each other.
Then,
the morphism
\begin{equation}
h^*{\cal F}
\to R{\cal H}om({\cal G},Rh^!{\cal E})
\label{eqAB1}
\end{equation}
induced by the second one
equals the adjoint of
the composition
\begin{equation}
{\cal F}
\to R{\cal H}om(Rh_!{\cal G},{\cal E})
\to Rh_*
R{\cal H}om({\cal G},Rh^!{\cal E})
\label{eqAB2}
\end{equation}
of that induced by the first
one and the inverse of the
canonical isomorphism of adjunction.
\end{lm}

\proof{
We consider the diagram
\begin{equation}
\xymatrix{
{\cal F}
\ar[d]\ar[rd]&
\\
Rh_*
R{\cal H}om
({\cal G},
h^*{\cal F}\otimes^L{\cal G})
\ar[d]\ar[r]
&
R{\cal H}om
(Rh_!{\cal G},
Rh_!(h^*{\cal F}\otimes^L{\cal G}))
\ar[d]\\
Rh_*
R{\cal H}om
({\cal G},
Rh^!{\cal E})
\ar[rd]\ar[r]
&
R{\cal H}om
(Rh_!{\cal G},
Rh_!Rh^!{\cal E})
\ar[d]\\
&
R{\cal H}om
(Rh_!{\cal G},
{\cal E}).}
\label{eqdgAB}
\end{equation}
The top vertical arrow is
the adjoint of $h^*{\cal F}
\to 
R{\cal H}om
({\cal G},
h^*{\cal F}\otimes^L{\cal G})$
induced by the identity of
$h^*{\cal F}\otimes^L{\cal G}$
and the horizontal arrows
are defined by the functoriality
of $Rh_!$.
The upper slant arrow is induced
by the isomorphism of projection formula 
(\ref{eqprj}) and
the upper triangle is commutative.
The middle vertical arrows are
induced by the second morphism
of (\ref{eqAB0})
and the middle square is commutative.
The lower slant arrow
is the isomorphism of adjunction
and the lower vertical arrow
is induced by the adjunction map.
Thus, the diagram (\ref{eqdgAB})
is commutative.

The composition of 
the left column in
the diagram (\ref{eqdgAB})
is the adjoint of 
(\ref{eqAB1}).
Since the composition 
of the upper slant arrow
and the right column 
is induced by
the first morphism of (\ref{eqAB0}),
the assertion follows.
\qed}
\medskip

Let ${\cal F}$ and ${\cal G}$ 
be constructible complexes
of $\Lambda$-modules
on $X$
and assume that
${\cal G}$ is of finite tor-dimension.
We define a canonical morphism
\begin{equation}
c_{h,{\cal F},{\cal G}}\colon
h^*{\cal F}\otimes^L_\Lambda
Rh^!{\cal G} 
\to Rh^!({\cal F}\otimes^L_\Lambda{\cal G} )
\label{eqFG}
\end{equation}
to be that corresponding
to the morphism
${\cal F}
\otimes^L_\Lambda
Rh_!Rh^!{\cal G} \to
{\cal F}
\otimes^L_\Lambda
{\cal G}$ induced by the adjunction
$Rh_!Rh^!{\cal G} \to
{\cal G}$.

\begin{lm}\label{lmhRh}
Let $h\colon W\to X$
be a separated
morphism of schemes.
Let ${\cal F}$ and ${\cal G}$ 
be constructible complexes
of $\Lambda$-modules
on $X$
and assume that
${\cal G}$ is of finite tor-dimension.

{\rm 1.}
Let $g\colon V\to W$
be a separated
morphism of schemes
and $f\colon V\to X$
be the composition.
Then the morphisms {\rm (\ref{eqFG})}
for $h,g$ and $f$ form
a commutative diagram
\begin{equation}
\xymatrix{
f^*{\cal F}\otimes^L Rf^!{\cal G}
\ar[d]\ar[rrrr]^{c_{f,{\cal F},{\cal G}}}&&&&
Rf^!({\cal F}\otimes^L{\cal G})\ar[d]\\
g^*h^*{\cal F}\otimes^L Rg^!Rh^!{\cal G}
\ar[rr]^{c_{g,h^*{\cal F},Rh^!{\cal G}}}&&
Rg^!(h^*{\cal F}\otimes^L Rh^!{\cal G})
\ar[rr]^{Rg^!(c_{h,{\cal F},{\cal G}})}&&
Rg^!Rh^!({\cal F}\otimes^L{\cal G}).}
\label{eqhgf}
\end{equation}
The vertical arrows are the
canonical isomorphisms.

{\rm 2.}
Let ${\cal E}$ 
be another constructible complex
of $\Lambda$-modules
on $X$
and assume that ${\cal F}$
is of finite tor-dimension.
Then the morphisms {\rm (\ref{eqFG})}
form a commutative diagram
\begin{equation}
\xymatrix{
h^*{\cal E}\otimes^L h^*
{\cal F}\otimes^L Rh^!{\cal G}
\ar[d]
\ar[rr]
^{{\rm id}\otimes c_{h,{\cal F},{\cal G}}}
&&
h^*{\cal E}\otimes^L Rh^!
({\cal F}\otimes^L {\cal G})
\ar[d]^{c_{h,{\cal E},{\cal F}\otimes{\cal G}}}
\\
h^*({\cal E}\otimes^L
{\cal F})\otimes^L Rh^!{\cal G}
\ar[rr]^{c_{h,{\cal E}\otimes{\cal F},{\cal G}}}&&
Rh^!({\cal E}\otimes^L
{\cal F}\otimes^L{\cal G}).}
\label{eqEFG}
\end{equation}

{\rm 3.}
Assume that $h$ is an immersion.
Then the morphism {\rm (\ref{eqFG})}
is induced by the restriction of
\begin{equation}
{\cal F}\otimes^L Rh_*Rh^!{\cal G}
=
{\cal F}\otimes^L 
R{\cal H}om (h_!\Lambda,
{\cal G})
\to 
Rh_*Rh^!(
{\cal F}\otimes^L {\cal G})
= 
R{\cal H}om (h_!\Lambda,
{\cal F}\otimes^L{\cal G}).
\label{eqhim}
\end{equation}
\end{lm}

\proof{
1.
It follows from the commutative 
diagram 
\begin{equation}
\xymatrix{
Rf_!(f^*{\cal F}\otimes^L Rf^!{\cal G})
\ar[d]\ar[rr]&&
{\cal F}\otimes^L
Rf_!Rf^!{\cal G}\ar[d]\\
Rh_!Rg_!(g^*h^*{\cal F}\otimes^L Rg^!Rh^!{\cal G})
\ar[r]&
Rh_!(h^*{\cal F}\otimes^L Rg_!Rg^!Rh^!{\cal G})
\ar[r]&
{\cal F}\otimes^L
Rh_!Rg_!Rg^!Rh^!{\cal G}.}
\label{eqhgpf}
\end{equation}
for the isomorphisms
of projection formula.

2.
We consider the diagram
$$\begin{CD}
Rh_!(h^*{\cal F}
\otimes^LRh^!{\cal G})
@>>>
Rh_!Rh^!({\cal F}
\otimes^L{\cal G})\\
@AAA@VVV\\
{\cal F}
\otimes^LRh_!Rh^!{\cal G}
@>>>
{\cal F}
\otimes^L{\cal G}.
\end{CD}$$
The right vertical arrow
is the adjunction morphism
and the bottom horizontal arrow
is induced by the adjunction morphism.
The left vertical arrow is
the isomorphism (\ref{eqprj})
of the projection formula 
and the top horizontal arrow is
induced by (\ref{eqFG}).
Since (\ref{eqFG}) is the adjoint
of the diagonal composition and is defined
by the commutativity of
the lower left triangle,
the diagram is commutative.
Tensoring ${\cal E}$
and applying the projection formula,
we obtain the adjoint of
(\ref{eqEFG}).

3.
We may assume that $h$ is a closed immersion.
Then, the composition of the
morphism
(\ref{eqhim}) with
the adjunction
$Rh_*Rh^!({\cal F}\otimes^L{\cal G})
\to {\cal F}\otimes^L{\cal G}$
is the same as
${\cal F}\otimes^LRh_*Rh^!{\cal G}
\to {\cal F}\otimes^L{\cal G}$
defining (\ref{eqFG}) and the assertion follows.
\qed}
\medskip

Assume ${\cal G}=\Lambda$
and we consider
the canonical morphism
\begin{equation}
c_{h,{\cal F}}
=c_{h,{\cal F},\Lambda}
\colon 
h^*{\cal F}\otimes^L_\Lambda Rh^!\Lambda
\to Rh^!{\cal F}.
\label{eqprpg}
\end{equation}
For another morphism
$g\colon V\to W$ and
the composition $f$,
The canonical morphisms
(\ref{eqFG})
form a commutative diagram
\begin{equation}
\xymatrix{
f^*{\cal F}
\otimes^L 
Rf^!\Lambda
\ar[r]^{c_{f,{\cal F}}}\ar[d]&
Rf^!{\cal F}
\ar[r]&
Rg^!Rh^!{\cal F}
\\
g^*h^*{\cal F}
\otimes^L 
Rg^!Rh^!\Lambda
\ar[rr]
^{c_{g,h^*{\cal F},Rh^!\Lambda}}
&&
Rg^!(h^*{\cal F}
\otimes^L 
Rh^!\Lambda)\ar[u]
_{Rg^!(c_{h,{\cal F}})}
\\
g^*h^*{\cal F}
\otimes^L 
Rg^!\Lambda
\otimes^L 
g^*Rh^!\Lambda
\ar[rr]^{c_{g,h^*{\cal F}}
\otimes{\rm id}}\ar[u]
^{{\rm id}\otimes
c_{g,Rh^!\Lambda}}&&
Rg^!h^*{\cal F}
\otimes^L 
g^*Rh^!\Lambda\ar[u]
_{c_{g,h^*{\cal F},Rh^!\Lambda}}
}
\label{eqhRh}
\end{equation}
by Lemma \ref{lmhRh}.1 and 2.

We give another description
of the morphism
$c_{h,{\cal F}}\colon
h^*{\cal F}\otimes^LRh^!\Lambda
\to Rh^!{\cal F}$
assuming further that ${\cal F}$
is of finite tor-dimension
and that $X$ is separated.
Recall that the dual $D_X{\cal F}$
is defined as $R{\cal H}om_\Lambda
({\cal F},{\cal K}_X)$ where
${\cal K}_X=Ra^!\Lambda$ 
for the structure morphism 
$a\colon X\to {\rm Spec}\ k$.
For a separated morphism $h\colon W\to X$
over $k$,
we have a canonical isomorphism
$Rh^!{\cal K}_X=
Rh^!Ra^!\Lambda
\to R(ah)^!\Lambda={\cal K}_W$.
The isomorphism of adjunction
$Rh_*D_W\to D_XRh_!$
induces an isomorphism
$D_Wh^*D_X\to Rh^!$.

\begin{lm}\label{lmAB2}
Assume that ${\cal F}$
is of finite tor-dimension and
we consider the morphisms
\begin{equation}
h^*D_X{\cal F}
\otimes^L
h^*{\cal F}
\otimes^L
Rh^!\Lambda
\to
h^*{\cal K}_X
\otimes^L
Rh^!\Lambda
\overset{c_{h,{\cal K}_X}}\longrightarrow
Rh^!{\cal K}_X\to {\cal K}_W
\label{eqFGDD}
\end{equation}
where the first arrow
is induced by the canonical morphism
$D_X{\cal F}
\otimes^L{\cal F}\to {\cal K}_X$.
Then, the composition
\begin{equation}
h^*{\cal F}
\otimes^L
Rh^!\Lambda
\to
D_Wh^*D_X{\cal F}
\to
Rh^!{\cal F}
\label{eqFGD}
\end{equation}
of the morphism induced
by {\rm (\ref{eqFGDD})}
with the canonical isomorphism
equals
the morphism 
$c_{h,{\cal F}}\colon
h^*{\cal F}
\otimes^L
Rh^!\Lambda
\to Rh^!{\cal F}$.
\end{lm}

\proof{
Let $h^*D_X{\cal F}
\to D_W(h^*{\cal F}
\otimes^LRh^!\Lambda)$
be the morphism induced by
(\ref{eqFGDD})
and let
\begin{equation}
D_X{\cal F}
\to 
Rh_*D_W(h^*{\cal F}
\otimes^LRh^!\Lambda)
\overset \simeq \to D_XRh_!(h^*{\cal F}
\otimes^LRh^!\Lambda)
\label{eq717}
\end{equation}
be its adjoint.
Then, the morphism (\ref{eqFGD})
is the adjoint of
the dual $Rh_!(h^*{\cal F}
\otimes^LRh^!\Lambda)
\to {\cal F}$ of (\ref{eq717})
since the isomorphism
$D_Wh^*D_X\to Rh^!$
is induced by
the isomorphism $Rh_*D_W\to D_XRh_!$.
The morphism
$D_X{\cal F}
\otimes^LRh_!(h^*{\cal F}
\otimes^LRh^!\Lambda)
\to {\cal K}_X$
corresponding to (\ref{eq717})
corresponds
to (\ref{eqFGDD})
by Lemma \ref{lmAB}.
Thus, it is induced by
the isomorphism
of the projection formula,
the canonical morphism
$D_X{\cal F}
\otimes^L{\cal F}\to {\cal K}_X$
and 
the adjunction morphism
$Rh_!h^!\Lambda\to \Lambda$.
Hence the morphism
(\ref{eqFGD}) equals
$c_{h,{\cal F}}$.
\qed}

\subsection{${\cal F}$-transversality}\label{ssFt}

\begin{df}\label{dfprpg}
Let $h\colon W\to X$ 
be a morphism of
schemes of finite type over $k$
and let ${\cal F}$ be a constructible complex
of $\Lambda$-modules on $X$.

In the case $h$ is {\em separated},
we say that $h$
is {\em ${\cal F}$-transversal}
if the canonical morphism 
\begin{equation*}
c_{h,{\cal F}}\colon
h^*{\cal F}\otimes^L_\Lambda Rh^!\Lambda
\to Rh^!{\cal F}
\leqno{\rm(\ref{eqprpg})}
\end{equation*}
%defined as {\rm (\ref{eqFG})} for ${\cal G}=\Lambda$
is an isomorphism.
For general $h$,
we say that $h$
is {\em ${\cal F}$-transversal}
if there exists an open covering
$W=\bigcup_iW_i$ consisting
of open subschemes
separated over $X$
such that
the restrictions $h|_{W_i}$
are ${\cal F}$-transversal.
\end{df}

Since the condition is local, a morphism
$h\colon W\to X$ is ${\cal F}$-transversal
if and only if
the restriction of $h$ to every open subset
of $W$ separated over $X$
is ${\cal F}$-transversal.

\begin{lm}\label{lmprpg}
Let ${\cal F}$ be a constructible complex
of $\Lambda$-modules on 
a scheme $X$ 
of finite type over $k$
and let $h\colon W\to X$
be a morphism of schemes
of finite type over $k$.

{\rm 1.}
If $h$ is smooth,
then $h\colon W\to X$
is ${\cal F}$-transversal.

{\rm 2.}
If ${\cal F}$ is locally constant,
then $h\colon W\to X$ is
${\cal F}$-transversal.

{\rm 3.}
Assume that
$h\colon W\to X$
is ${\cal F}$-transversal
and that 
$Rh^!\Lambda$ is isomorphic
to $\Lambda(c)[2c]$ for a
locally constant function $c$ on $W$.
Then for a morphism
$g\colon V\to W$
of schemes
of finite type over $k$,
the following conditions are equivalent:

{\rm (1)}
$g\colon V\to W$ 
is $h^*{\cal F}$-transversal.

{\rm (2)}
The composition
$f\colon V\to W$ 
is ${\cal F}$-transversal.

{\rm 4.}
Let $\Lambda_0$
be the residue field of $\Lambda$.
Assume that ${\cal F}$
is of finite tor-dimension
and set ${\cal F}_0=
{\cal F}\otimes^L_\Lambda\Lambda_0$,
Then,
$h\colon W\to X$ is 
${\cal F}$-transversal
if and only if
it is ${\cal F}_0$-transversal.

{\rm 5.}
Assume that ${\cal F}$ is a perverse sheaf.
Assume further that $h\colon
W\to X$ is locally of
complete intersection
of relative virtual dimension
$c$
and that 
$Rh^!\Lambda$ is isomorphic
to $\Lambda(c)[2c]$.
If $h\colon W\to X$ is
${\cal F}$-transversal,
then $h^*{\cal F}[c]$
and $Rh^!{\cal F}[-c]$
are perverse sheaves.

{\rm 6.}
Assume that ${\cal F}$
is of finite tor-dimension.
Then the following conditions are
equivalent

{\rm (1)}
$h$ is ${\cal F}$-transversal.

{\rm (2)}
The morphism
$h^*{\cal F}\otimes^L
Rh^!\Lambda
\to D_Wh^*D_X{\cal F}$
{\rm (\ref{eqFGD})}
induced by the morphism
$h^*D_X{\cal F}\otimes^L
h^*{\cal F}\otimes^L
Rh^!\Lambda
\to {\cal K}_W$
{\rm (\ref{eqFGDD})}
is an isomorphism.
\end{lm}

We show a converse of 2.\
in Corollary \ref{corfh}.

\proof{
By the remark after Definition \ref{dfprpg},
we may assume that morphisms
are separated.
We will omit to write
this remark in the sequel.

1. 
We may assume that
$h$ is smooth of relative dimension $d$.
We consider the canonical
isomorphism
$\Lambda(d)[2d]\to Rh^!\Lambda$
defined as the adjoint of
the trace mapping
$Rh_!\Lambda(d)[2d]\to \Lambda$.
Then, by the description
of the isomorphism
(\ref{eqPD}) loc.\ cit.,
the diagram
$$\begin{CD}
h^*{\cal F}\otimes^L
\Lambda(d)[2d]@>{\rm (\ref{eqPD})}>>
Rh^!{\cal F}\\
@VVV @|\\
h^*{\cal F}\otimes^L
Rh^!\Lambda
@>{c_{h,{\cal F}}}>>
Rh^!{\cal F}
\end{CD}$$
is commutative
and
the morphism $c_{h,{\cal F}}$
is an isomorphism.

2.
Since the assertion is \'etale local on $X$,
it is reduced to the case where ${\cal F}$
is constant by devissage.

3.
We consider the 
commutative diagram
$$\xymatrix{
f^*{\cal F}
\otimes^L 
Rf^!\Lambda
\ar[r]^{c_{f,{\cal F}}}\ar[d]&
Rf^!{\cal F}
\ar[r]&
Rg^!Rh^!{\cal F}
\\
g^*h^*{\cal F}
\otimes^L 
Rg^!Rh^!\Lambda
%\ar[rr]^{c_{g,h^*{\cal F},Rh^!\Lambda}}
&&
Rg^!(h^*{\cal F}
\otimes^L 
Rh^!\Lambda)\ar[u]
_{Rg^!(c_{h,{\cal F}})}
\\
g^*h^*{\cal F}
\otimes^L 
Rg^!\Lambda
\otimes^L 
g^*Rh^!\Lambda
\ar[rr]^{c_{g,h^*{\cal F}}
\otimes{\rm id}}\ar[u]
^{{\rm id}\otimes
c_{g,Rh^!\Lambda}}&&
Rg^!h^*{\cal F}
\otimes^L 
g^*Rh^!\Lambda\ar[u]
_{c_{g,h^*{\cal F},Rh^!\Lambda}}
.}
\leqno{\rm (\ref{eqhRh})}$$
Since 
$Rh^!\Lambda$ is assumed to be isomorphic
to $\Lambda(c)[2c]$ for a
locally constant function $c$ on $W$,
the lower vertical arrows are isomorphisms.
By the assumption that
$h\colon W\to X$
is ${\cal F}$-transversal,
the upper right vertical arrow
$Rg^!(c_{h,{\cal F}})$ is
an isomorphism.
The condition (1) means that
the top horizontal arrow
$c_{f,{\cal F}}$ is
an isomorphism
and the condition (2) is equivalent
to that the bottom horizontal arrow
$c_{g,h^*{\cal F}}
\otimes{\rm id}$ is
an isomorphism by 
the assumption
that 
$Rh^!\Lambda$ is isomorphic
to $\Lambda(c)[2c]$.
Hence the equivalence follows.

4.
Similarly as Lemma \ref{lmctf}.1,
the canonical morphism
$Rf^!{\cal F}\otimes_\Lambda^L
\Lambda_0
\to
Rf^!{\cal F}_0$
is an isomorphism.
Hence, similarly as the proof
of Lemma \ref{lmctf}.2,
the morphism $c_{h,{\cal F}}$
is an isomorphism
if and only if 
$c_{h,{\cal F}_0}$
is an isomorphism.

5. 
If $h$ is smooth, $h^*(c)[c]=Rh^![-c]$
is $t$-exact by \cite[4.2.4]{BBD}
and the assertion follows in this case.
Hence, we may assume
that $h$ is a regular immersion
of codimension $-c$.
Then, by \cite[Corollary 4.1.10]{BBD}
and induction on $-c$,
we have 
$h^*{\cal F}(c)[c]\in{^p}\!D_c^{[0,-c]}(W,\Lambda)$
and
$Rh^!{\cal F}[-c]\in{ ^p}\!D_c^{[c,0]}(W,\Lambda)$.
By the assumption,
we have an isomorphism
$h^*{\cal F}(c)[c]\to Rh^!{\cal F}[-c]$
and the assertion follows.

6.
It follows from Lemma \ref{lmAB2}.
\qed}
\medskip

\begin{pr}\label{prRhom}
Let ${\cal F}$ be a constructible complex
of $\Lambda$-modules of
finite tor-dimension on 
a scheme $X$ of finite type over $k$.

{\rm 1.}
Let $h\colon W\to X$
be a morphism of schemes 
of finite type over $k$
and assume that
$Rh^!\Lambda$ is isomorphic
to $\Lambda(c)[2c]$ for a
locally constant function $c$ on $W$.
Then, the following conditions are equivalent:

{\rm (1)}
The morphism $h\colon W\to X$ 
is ${\cal F}$-transversal.

{\rm (2)}
The morphism $h\colon W\to X$ 
is $D_X{\cal F}$-transversal.

\noindent
Further if $h\colon W\to X$ 
is ${\cal K}_X$-transversal,
they are 
equivalent to the following condition:

{\rm (3)}
The canonical morphism
$h^*R{\cal H}om_X({\cal F},{\cal K}_X)
\to 
R{\cal H}om_W(h^*{\cal F},h^*{\cal K}_X)$
is an isomorphism.

{\rm 2.}
Let ${\cal G}$ be a constructible complex
of $\Lambda$-modules on 
a scheme $X$ of finite type over $k$.
Then, the following conditions are equivalent:

{\rm (1)}
The diagonal morphism
$\delta\colon X\to X\times X$
is 
$R{\cal H}om_{X\times X}({\rm pr}_2^*{\cal F},
{\rm pr}_1^!{\cal G})$-transversal.

{\rm (2)}
The canonical morphism
${\cal G}\otimes^L
R{\cal H}om_X({\cal F},\Lambda)
\to
R{\cal H}om_X({\cal F},{\cal G})$
is an isomorphism.
\end{pr}

The assumptions in Proposition \ref{prRhom}.1
are satisfied if $X$ and $W$ are
smooth over $k$.
The canonical morphism
in Proposition \ref{prRhom}.1(3)
is an analogue of 
$h^*Sol\ \!{\cal M}\to Sol h^*{\cal M}$
for a ${\cal D}$-module ${\cal M}$.

\proof{
1.
The equivalence of (1) and (2)
follows from
Lemma \ref{lmprpg}.6,
the assumption
$Rh^!\Lambda
\simeq \Lambda(c)[2c]$ and
the isomorphism
${\rm id}\to D_WD_W$ of biduality
\cite[Th\'eor\`eme 4.3]{TF}.

The assumption that $h\colon W\to X$ 
is ${\cal K}_X$-transversal
means that the canonical morphism
$c_{h,{\cal K}_X}\colon
h^*{\cal K}_X\otimes^L Rh^!\Lambda
\to Rh^!{\cal K}_X={\cal K}_W$
is an isomorphism.
Hence by the definition of
(\ref{eqFGDD}),
the conditions (1) and (3)
are equivalent.

2.
We consider the commutative diagram
$$
\begin{CD}
\delta^*R{\cal H}om_X({\rm pr}_2^*{\cal F},
{\rm pr}_1^!{\cal G})
\otimes^L R\delta^!\Lambda
@>{c_{\delta,R{\cal H}om_X({\rm pr}_2^*{\cal F},
{\rm pr}_1^!{\cal G})}} >>
R\delta^!
R{\cal H}om_X({\rm pr}_2^*{\cal F},
{\rm pr}_1^!{\cal G})
\\
@AAA@VVV\\
{\cal G}\otimes^L
R{\cal H}om_X({\cal F},\Lambda)
@>>>
R{\cal H}om_X({\cal F},{\cal G})
\end{CD}$$
defined as follows.
The bottom horizontal arrow
is the canonical morphism
in the condition (2).
The canonical isomorphism
${\cal G}\boxtimes^L
D_X{\cal F}\to
R{\cal H}om_{X\times X}({\rm pr}_2^*{\cal F},
{\rm pr}_1^!{\cal G})$
\cite[(3.1.1)]{FL}
induces an isomorphism
${\cal G}\otimes^L
R{\cal H}om_X({\cal F},{\cal K}_X)
=\delta^*(
{\cal G}\boxtimes^L
D_X{\cal F})
\to 
\delta^*
R{\cal H}om_{X\times X}({\rm pr}_2^*{\cal F},
{\rm pr}_1^!{\cal G})$.
This together with the canonical isomorphism
$R\delta^!\Lambda \otimes {\cal K}_X
\to \Lambda$
defines the left vertical arrow.
The right vertical arrow is defined by
\cite[3.1.12.2]{DP}.
Since the condition (1) is equivalent
to that the top horizontal arrow
$c_{\delta,R{\cal H}om_X({\rm pr}_2^*{\cal F},
{\rm pr}_1^!{\cal G})}$
is an isomorphism,
the assertion follows.
\qed}

\begin{pr}\label{prpgj}
Let $$
\begin{CD}
W@>h>> X\\
@A{j'}AA@AAjA\\
V@>{h'}>>U
\end{CD}$$
be a cartesian diagram
of schemes 
of finite type over $k$
such that the vertical arrows
are open immersions.
Let ${\cal G}$ be a constructible complex
of $\Lambda$-modules on $U$.
We consider the conditions:

{\rm (1)}
The morphism $h\colon W\to X$ 
is $j_!{\cal G}$-transversal.

{\rm (2)}
The morphism $h'\colon V\to U$
is ${\cal G}$-transversal.

{\rm 1.}
The condition {\rm (1)}
implies {\rm (2)}.
Conversely,
if $Rh^!\Lambda$ is isomorphic to
$\Lambda(c)[2c]$ for a
locally constant function $c$
and if the canonical morphisms
\begin{equation}
j_!{\cal G}\to Rj_*{\cal G},\quad
j'_!h^{\prime *}{\cal G}\to 
Rj'_*h^{\prime *}{\cal G}
\label{eqj!j*}
\end{equation}
are isomorphisms,
the condition {\rm (2)}
implies {\rm (1)}.

{\rm 2.}
Assume that ${\cal G}$ is of finite
tor-dimension on $U$.
Then, the condition
{\rm (1)} is equivalent
to the combination of
{\rm (2)}
and the following condition:

{\rm (3)}
The base change morphism
\begin{equation}
\begin{CD}
h^*Rj_*R{\cal H}om({\cal G},{\cal K}_U)
@>>>
Rj'_*h^{\prime *}
R{\cal H}om({\cal G},{\cal K}_U)
\end{CD}
\label{eqbcg}
\end{equation}
is an isomorphism.
\end{pr}

\proof{1.
The implication {\rm (1)}$\Rightarrow
${\rm (2)} is clear.

We consider the commutative diagram
$$\xymatrix{
j'_!h^{\prime*}{\cal G}\otimes^L
Rh^!\Lambda
\ar[r]\ar[d]
&
h^*j_!{\cal G}\otimes^L
Rh^!\Lambda
\ar[rr]^{c_{h,j_!{\cal G}}}&&
Rh^!j_!{\cal G}\ar[d]\\
Rj'_*h^{\prime*}{\cal G}\otimes^L
Rh^!\Lambda\ar[dr]&
&&
Rh^!Rj_*{\cal G}\ar[d]\\
&
Rj'_*(h^{\prime*}{\cal G}\otimes^L
Rh^{\prime !}\Lambda)
\ar[rr]^{Rj'_*(c_{h',{\cal G}})}&&
Rj'_*Rh^{\prime!}{\cal G}
}$$
defined as follows.
The upper vertical arrows
are induced by the
canonical morphisms (\ref{eqj!j*}).
The top left horizontal arrow
is induced by the isomorphism
$j'_!h^{\prime*}\to 
h^*j_!$ and is an isomorphism.
The slant arrow is defined 
as the adjoint of the isomorphism
$j^{\prime*}
(Rj'_*h^{\prime*}{\cal G}\otimes^L
Rh^!\Lambda)\to 
Rj'_*(h^{\prime*}{\cal G}\otimes^L
Rh^{\prime !}\Lambda)$
and is an isomorphism
if the assumption on
$Rh^!\Lambda$ is satisfied.
The lower right vertical arrow
is the adjoint of the isomorphism
$j^*Rh_!\to Rh'_!j^{\prime *}$
and is an isomorphism.
Thus under the assumptions,
the implication {\rm (2)}$\Rightarrow
${\rm (1)} holds.

2.
Since the condition (1)
implies (2),
it suffices to show that
(\ref{eqbcg}) is an isomorphism
if and only if
the condition (3) for ${\cal F}=j_!{\cal G}$
in Proposition \ref{prRhom}.1 is satisfied,
assuming (2).
We consider the
commutative diagram
$$\begin{CD}
h^*R{\cal H}om_X(j_!{\cal G},{\cal K}_X)
@>{(3)}>>
R{\cal H}om_W(h^*j_!{\cal G},
h^*{\cal K}_X)
\\
@VVV@VVV\\
h^*Rj_*R{\cal H}om_U({\cal G},{\cal K}_U)
@.
R{\cal H}om_W(j'_!h^{\prime*}{\cal G},
h^*{\cal K}_X)
\\
@V{(\ref{eqbcg})}VV@VVV\\
Rj_*h^{\prime*}R{\cal H}om_U({\cal G},{\cal K}_U)
@>>>
Rj'_*R{\cal H}om_V(h^{\prime*}
{\cal G},h^{\prime*}{\cal K}_U)
\end{CD}$$
defined as follows.
The upper left and the
lower right vertical arrows are
the adjunction morphisms
and are isomorphisms.
The upper right vertical arrow
is induced by the isomorphism
$h^*j_!{\cal G}\to
j'_!h^{\prime*}{\cal G}$
and is an isomorphism.

The top horizontal arrow (3) is
the canonical morphism in
the condition (3) in Proposition 
\ref{prRhom}.1 for ${\cal F}=j_!{\cal G}$.
The lower one is induced by
that for ${\cal G}$
and is an isomorphism
if (2) is satisfied,
by Proposition \ref{prRhom}.1
(1)$\Rightarrow$(3).
Thus the assertion follows.
\qed}
%\medskip

\subsection{${\cal F}$-transversality
and local acyclicity}\label{ssFla}

In Proposition \ref{prla},
%and Corollary \ref{corlat},
$X$ and $Y$ denote arbitrary schemes
and $\Lambda$ denotes a finite local ring
such that the characteristic of the
residue field is invertible on $Y$.

\begin{pr}[{\rm \cite[Proposition 2.10]{app}}]
\label{prla}
Let 
\begin{equation}
\begin{CD}
X@<{j'}<< U\\
@VfVV @VV{f_V}V\\
Y@<j<<V
\end{CD}
\label{eqXYi}
\end{equation}
be a cartesian diagram of schemes
such that $j\colon V\to Y$
is an open immersion.
Let ${\cal F}$ 
be a complex of $\Lambda$-modules
on $X$
and let ${\cal G}$
be a complex of $\Lambda$-modules
on $V$.
Assume that 
$f$ is strongly locally acyclic relatively to ${\cal F}$.
Then, the base change morphism
\begin{equation}
{\cal F}\otimes f^*Rj_*{\cal G}
\to 
Rj'_*(j^{\prime*}{\cal F}
\otimes f_V^*{\cal G})
\label{eqla1}
\end{equation}
is an isomorphism.
\end{pr}

\begin{cor}\label{corlat}
Let $f\colon X\to Y$ be a 
{\em smooth} morphism
of schemes of finite type
over a field $k$,
$i\colon Z\to Y$
be an immersion
and 
\begin{equation}
\begin{CD}
W@>h>>X\\
@VgVV @VfVV\\
Z@>i>>Y
\end{CD}
\label{eqZiY}
\end{equation}
be a cartesian diagram of schemes.
Let ${\cal F}$
be a complex of $\Lambda$-modules
on $X$
and assume that
$f$ is strongly locally acyclic relatively to ${\cal F}$.
Then 
$h\colon W\to X$
is ${\cal F}$-transversal.
\end{cor}

\proof{
We may assume that $i\colon Z\to Y$
is a closed immersion.
Let $V=Y\sm Z$
and consider the cartesian diagram
\begin{equation}
\begin{CD}
W@>h>>X@<{j'}<< U\\
@VgVV @VfVV @VV{f_V}V\\
Z@>i>>Y@<j<<V.
\end{CD}
\label{eqZYV}
\end{equation}
By Proposition \ref{prla} applied
to the right square,
we obtain an isomorphism
${\cal F}\otimes f^*Rj_*\Lambda
\to 
Rj'_*j^{\prime*}{\cal F}$.
Since $f$ is smooth,
this induces an isomorphism
${\cal F}\otimes Rj'_*\Lambda
\to 
Rj'_*j^{\prime*}{\cal F}$
by smooth base change theorem
\cite{smbc}.
By the distinguished triangle
$\to h_*Rh^!\to{\rm id}\to
Rj'_*j^{\prime*}$,
we obtain an isomorphism
$c_{h,{\cal F}}\colon
h^*{\cal F}\otimes Rh^!\Lambda
\to 
Rh^!{\cal F}$.
\qed}

\begin{pr}\label{prF}
Let $f\colon X\to Y$ be a 
morphism of schemes of finite
type over a field $k$
and let ${\cal F}$
be a complex of $\Lambda$-modules
on $X$.
The morphism
$f$ is locally acyclic relatively to ${\cal F}$
if the following 
condition is satisfied:

Let $Y'$ and $Z$ be smooth
schemes 
over a finite extension of $k$
and 
\begin{equation}
\begin{CD}
X@<{p'}<< X'@<h<<W\\
@VfVV @VV{f'}V@VVgV\\
Y@<p<<Y'@<i<<Z
\end{CD}
\label{eqXYpi}
\end{equation}
be a cartesian diagram of
schemes where
$p\colon Y'\to Y$ is
proper and generically finite and
$i\colon Z\to Y'$ is a closed immersion.
Then, 
for the pull-back
${\cal F}'=p^{\prime *}{\cal F}$ on $X'$
the composition
\begin{equation}
h^*{\cal F}'\otimes g^*Ri^!\Lambda
\to 
h^*{\cal F}'\otimes Rh^!\Lambda
\overset{c_{h,{\cal F}'}}\longrightarrow 
Rh^!{\cal F}'
\label{eqF}
\end{equation}
where the first arrow
is induced by the base change morphism
is an isomorphism.
\end{pr}

Though
Proposition \ref{prF}
is stated for schemes of finite
type over a field $k$,
it can be generalized to
more general schemes
using more general versions  of
alteration.

\proof{
Let
\begin{equation}
\begin{CD}
X@<{p'}<< X'@<{j'}<<V\\
@VfVV @VV{f'}V@VV{f'_U}V\\
Y@<p<<Y'@<j<<U
\end{CD}
\label{eqXYUpi}
\end{equation}
be a cartesian diagram of
schemes where
$p\colon Y'\to Y$ is
proper and generically finite
and $j\colon U=Y\sm D
\to Y'$ is the open immersion
of the complement of
a divisor with simple normal crossings.
Then, by the assumption
and the relative property,
the canonical morphism
\begin{equation}
{\cal F}'\otimes f^{\prime*}Rj_*\Lambda
\to 
Rj'_*j^{\prime*}{\cal F}'
\label{eqY'U}
\end{equation}
is an isomorphism.
By projection formula
and the proper base change theorem,
(\ref{eqY'U}) shows that
the canonical morphism
\begin{equation}
{\cal F}\otimes f^*R(pj)_*\Lambda
\to 
R(p'j')_*(p'j')^*{\cal F}
\label{eqYU}
\end{equation}
is an isomorphism.

Let $s\gets t$ be a specialization
of geometric points of $Y$,
let $Y_{(s)}$ denote the strict localization
and let
\begin{equation}
\begin{CD}
X_s@>{i'_s}>>X\times_YY_{(s)}@<{j'_t}<< X_t\\
@V{f_s}VV@V{f_{(s)}}VV @VV{f_t}V\\
s@>{i_s}>>Y_{(s)}@<{j_t}<<t
\end{CD}
\label{eqXYtpi}
\end{equation}
be the cartesian diagram.
We assume that $t$ is
the spectrum of an algebraic
closure of the residue field
of the point of $Y$ below $t$.
It suffices to show
that the canonical morphism 
\begin{equation}
i^{\prime*}_s{\cal F}
\to 
i^{\prime*}_sRj'_{t*}j^{\prime*}_t{\cal F}
\label{eqijst}
\end{equation}
is an isomorphism.

By \cite[Theorem 4.1]{dJ},
we may write $t$ as a limit
$\varprojlim_\lambda
U_\lambda$
of the complements $U_\lambda=Y_\lambda\sm 
D_\lambda$,
in proper and generically
finite schemes $Y_\lambda$
over $Y$ 
smooth over finite extensions
of $k$ of divisors $D_\lambda\subset
Y_\lambda$ with simple normal crossings.
Then, as the limit of
(\ref{eqYU}), the canonical morphism
\begin{equation}
{\cal F}\otimes f_{(s)}^*
Rj_{t*}j^*_t\Lambda
\to 
Rj'_{t*}j^{\prime*}_t{\cal F}
\label{eqijst2}
\end{equation}
is an isomorphism.
Since
$\Lambda\to
i_s^*Rj_{t*}j^*_t\Lambda$
is an isomorphism,
the isomorphism
(\ref{eqijst2})
induces an isomorphism
(\ref{eqijst}).
\qed}

\begin{cor}\label{corF}
Let $f\colon X\to Y$ be a 
{\em smooth} morphism of 
schemes of finite type over
a field $k$
and let ${\cal F}$
be a complex of $\Lambda$-modules
on $X$.
Assume that
for every cartesian diagram 
{\rm (\ref{eqXYpi})} satisfying the condition
there,
the immersion $h\colon W\to X'$
%in the cartesian diagram {\rm (\ref{eqXYpi})}
is ${\cal F}'$-transversal
for the pull-back ${\cal F}'=p^{\prime*}{\cal F}$
by $p'\colon X'\to X$.
Then,
$f$ is locally acyclic relatively to ${\cal F}$.
\end{cor}

\proof{
Assume that a cartesian diagram 
{\rm (\ref{eqF})} satisfies the condition
there.
Since $f'$ is smooth,
the base change morphism
$g^*Ri^!\Lambda\to Rh^!\Lambda$
is an isomorphism
by the smooth base change theorem
\cite{smbc}.
Hence the isomorphism
(\ref{eqF}) means that
$h$ is ${\cal F}'$-transversal.
Thus it suffices to apply
Proposition \ref{prF}.
\qed}

\subsection{${\cal F}$-transversality and singular support}\label{ssFs}

Let $k$ be a field and
$\Lambda$ be a finite field
of characteristic invertible in $k$.

\begin{pr}\label{prfh}
Let $X$ be a smooth scheme
over a field $k$
and 
$\Lambda$ be a finite field
of characteristic invertible in $k$.
Let ${\cal F}$ be a constructible complex
of $\Lambda$-modules on $X$.
For a closed conical subset $C
\subset T^*X$,
we consider the following conditions:

{\rm (1)}
${\cal F}$ is micro-supported on $C$.

{\rm (2)}
The support of ${\cal F}$
is a subset of
the base $B=C\cap T^*_XX\subset X$ 
of $C$ and
every $C$-transversal
morphism $h\colon W\to X$
of smooth schemes
of finite type
is ${\cal F}$-transversal.

{\rm 1.}
If ${\cal F}$ is 
of finite tor-dimension,
then the condition {\rm (1)} implies {\rm (2)}.

{\rm 2.}
If $k$ is {\em perfect},
the condition {\rm (2)} implies {\rm (1)}.
\end{pr}

\proof{
1. 
By Lemma \ref{lmlcst}.3,
the condition that
${\cal F}$ is micro-supported on $C$
implies that the support of ${\cal F}$
is a subset of the base $B$ of $C$.

We show that a $C$-transversal morphism 
$h\colon W\to X$ of smooth schemes of finite type
is ${\cal F}$-transversal
assuming that
${\cal F}$ is micro-supported on $C$.
By applying Lemmas \ref{lmCh}.3 
and \ref{lmprpg}.3
to the graph 
$W\to W\times X\to X$ of $h$ and
by replacing $h\colon W\to X, C$
and ${\cal F}$ by 
$W\to W\times X, {\rm pr}_2^*C$ 
and ${\rm pr}_2^*{\cal F}$,
we may assume that $h$ is an immersion.
Since the assertion is local on $W$,
we may assume that
there exists a cartesian diagram
$$\begin{CD}
&W@>h>> X\\
&@VVV @VVfV\\
0=&y@>>>Y&={\mathbf A}^d
\end{CD}$$
where $f\colon X\to Y$ is smooth.
Since $h\colon W\to X$ is
$C$-transversal,
by Lemma \ref{lmCf}.1,
we may assume that $f\colon X\to Y$
is $C$-transversal.

Since ${\cal F}$ is assumed
to be micro-supported on $C$,
the morphism $f\colon X\to Y$ is locally
acyclic relatively to ${\cal F}$.
Hence the assertion follows from
Corollary \ref{corlat}.
%Lemma \ref{lmequiv}.1.

2.
%(2)$\Rightarrow$(1)
Let $h\colon W\to X$ and
$f\colon W\to Y$
be a $C$-transversal pair of
smooth morphisms.
It suffices to show that
$f\colon W\to Y$ is locally acyclic
relatively to $h^*{\cal F}$.
By replacing ${\cal F}$
and $C$ by $h^*{\cal F}$
and $h^\circ C$,
we may assume $W=X$.
Since the base $B$ of $C$
contains the support of ${\cal F}$,
after replacing $X$ by a neighborhood
of the support of ${\cal F}$,
we may assume that $f\colon X\to Y$
is smooth by Lemma \ref{lmCf}.6.

We show that the condition
in Corollary \ref{corF} is satisfied.
Let (\ref{eqF}) be a cartesian diagram
satisfying the condition there.
Since $f\colon X\to Y$
is smooth,
$Y'$ and $W$ are smooth over 
a finite extension of $k$
and since $k$ is assumed perfect,
the schemes $X'$ and $W$ are smooth over $k$.
By Lemma \ref{lmChf}.2,
the proper morphisms $p'\colon X'\to X$
and $p'h\colon W\to X$ are
$C$-transversal.
Thus, by the condition (2),
the morphisms $p'\colon X'\to X$
and $p'h\colon W\to X$ are
${\cal F}$-transversal.
Hence by Lemma \ref{lmprpg}.3,
$h\colon W\to X'$ is 
$p^{\prime *}{\cal F}$-transversal.
Thus by Corollary \ref{corF},
the morphism $f\colon X\to Y$ is
locally acyclic relatively to ${\cal F}$.
\qed}

\begin{cor}\label{corfh}
Assume that
$k$ is perfect
and let ${\cal F}$
be a constructible complex 
of $\Lambda$-modules 
of finite tor-dimension on $X$.
Then, 
the following conditions are equivalent:

{\rm (1)}
${\cal F}$ is locally constant.

{\rm (2)}
Every morphism $h\colon W\to X$
of finite type 
of smooth schemes
is ${\cal F}$-transversal.
\end{cor}

\proof{
The implication (1)$\Rightarrow$(2)
is Lemma \ref{lmprpg}.2.
By Lemma \ref{lmlcst},
the condition (1) is equivalent to
that ${\cal F}$ is micro-supported
on the $0$-section $T^*_XX$.
Hence, it follows from Proposition \ref{prfh}.
\qed}

\begin{cor}\label{corpSS}
Let ${\cal F}$ be a constructible
complex of $\Lambda$-modules
of finite tor-dimension on $X$
and let $C=SS{\cal F}$ denote the
singular support.
Then, for a properly $C$-transversal
morphism $h\colon W\to X$
of smooth schemes of
finite type over a perfect $k$,
we have 
\begin{equation}
SS h^*{\cal F}
=SS Rh^!{\cal F}=h^{\circ}SS{\cal F}.
\label{eqpSS}
\end{equation}
\end{cor}

\proof{
We may assume that $\Lambda$ is a field
by Corollary \ref{corss}.2.
By Proposition \ref{prfh}.1,
$h$ is ${\cal F}$-transversal. 
First, we assume that ${\cal F}$ is 
a perverse sheaf.
Then $h^*{\cal F}[\dim W-\dim X]$
is also a perverse sheaf
by Lemma \ref{lmprpg}.5
and $SS h^*{\cal F}$
is the support of
$CC h^*{\cal F}$ by Proposition \ref{prperv}.2.
Since
$CC h^*{\cal F}=h^!CC{\cal F}$ by Theorem \ref{thmi*}
and $CC{\cal F}\geqq 0$ by
Proposition \ref{prperv}.1, the assertion
follows in this case.
The general case 
is reduced to the case where
${\cal F}$ is a perverse sheaf by 
the equality (\ref{eqssu}) for 
${\cal F}$ and $h^*{\cal F}$.
\qed}

\end{document}